\documentclass[10pt]{amsart}
\usepackage{amsmath,amssymb,amsthm}
\usepackage{amsmath,amssymb,amsthm,amscd}
\numberwithin{equation}{section}

\DeclareMathOperator{\Aut}{Aut}
\DeclareMathOperator{\Eq}{Eq}
\DeclareMathOperator{\Autw}{\Aut_{{\mathfrak W}}}

\DeclareMathOperator{\autz}{\Aut_{{\mathfrak Z}}}
\DeclareMathOperator{\spec}{Spec}

\DeclareMathOperator{\Hom}{Hom}
\DeclareMathOperator{\Ext}{Ext}

\DeclareMathOperator{\image}{Im}

\DeclareMathOperator{\pr}{pr}

\def\W{{\mathfrak W}}

\def\F{{\mathfrak F}}

\def\W{{\mathfrak W}}

\def\cD{\mathcal D}

\def\cI{{\mathcal I}}
\def\cJ{{\mathcal J}}
\def\cK{{\mathcal K}}

\def\cO{{\mathcal O}}

\def\cU{{\mathcal U}}

\def\cW{{\mathcal W}}
\def\cX{{\mathcal X}}
\def\cY{{\mathcal Y}}
\def\cZ{{\mathcal Z}}

\def\bB{{\mathbf B}}

\def\bD{{\mathbf D}}

\def\bH{{\mathbf H}}

\def\bP{{\mathbf P}}
\def\bv{{\mathbf v}}

\def\bzero{{\mathbf 0}}

\def\bA{{\mathbf A}}

\def\bt{{\mathbf t}}

\def\bL{{\mathbf L}}

\def\bp{{\mathbf p}}

\def\bq{{\mathbf q}}
\def\kk{{\rm k}}

\def\bxi{{\boldsymbol \xi}}
\def\bleta{{\boldsymbol \eta}}
\def\bzero{{\boldsymbol 0}}

\def\bone{{\boldsymbol 1}}

\def\ZZ{{\mathbb Z}}

\def\CC{{\mathbb C}}

\def\WW{{\mathfrak W}}

\def\YY{{\mathfrak Y}}

\def\MM{{\mathfrak M}}

\def\Z{{\mathfrak Z}}

\def\AUT{\mathfrak{Aut}}

\newtheorem{prop}{Proposition}[section]
\newtheorem{theo}[prop]{Theorem}
\newtheorem{lemm}[prop]{Lemma}
\newtheorem{coro}[prop]{Corollary}

\newtheorem{defi}[prop]{Definition}

\def\begeq{\begin{equation}}
\def\endeq{\end{equation}}

\def\and{\quad{\rm and}\quad}

\def\adj{\sqcup}
\def\Ao{{{\mathbf A}^{\! 1}}}
\def\At{{\mathbf A}^{\! 2}}
\def\As{{\mathbf A}^{\! 3}}
\def\An{{{\mathbf A}^{\! n}}}

\def\Al{{\mathbf A}^{\! l}}
\def\Anl{{\mathbf A}^{\! n-l}}
\def\Almo{{\mathbf A}^{\! l-1}}
\def\Anpo{{\mathbf A}^{\!n+1}}
\def\Anmo{{\mathbf A}^{\!n-1}}

\def\Ampo{{\mathbf A}^{\! m+1}}
\def\Anpt{{\mathbf A}^{\! n+2}}

\def\bl{\bigl(}
\def\br{\bigr)}
\def\lbe{_{\beta}}

\def\bt{{\mathbf t}}

\def\bz{{\mathbf 0}}

\def\cn{C[n]}

\def\dual{^{\vee}}

\def\defeq{\triangleq}

\def\dn{D[n]}

\def\Gn{\Gamma[n]}

\def\gn{G[n]}

\def\hatR{{\hat R}}

\def\isom{{\mathrm{Isom}}}
\def\Isom{\mathfrak{Isom}}

\def\kk{{\mathbf k}}
\def\ktt{\kk[\bt]}
\def\ks{\kk[\![s]\!]}
\def\kw{\kk[w_1,w_2]}
\def\kz{\kk[\![z_1,z_2]\!]}
\def\kt{\kk[t]}

\def\lalpbe{_{\alpha\beta}}

\let\lra=\longrightarrow
\def\lsta{_{\ast}}

\def\lbe{_{\beta}}

\def\lalp{_{\alpha}}

\def\lnpo{_{n+1}}
\def\lxi{_{\xi}}

\def\lbPsi{\underline\Psi}

\let\lbar=\underline

\def\lef{^{-}}
\def\lbb{[\![}

\def\lnode{_{{\rm node}}}

\def\mapright#1{\,\smash{\mathop{\lra}\limits^{#1}}\,}

\def\mapto{\xrightarrow}
\def\mh{\!:\!}

\def\mm{{\mathfrak m}}

\def\Mzdg{\MM(\Zzd,\Gamma)}

\def\mwcg{\MM(\WW,\Gamma)}

\def\Mwg{\MM(\WW,\Gamma)}
\def\Mydo{\MM(\YYo,\Gamma_1)}
\def\Mydt{\MM(\YYt,\Gamma_2)}
\def\Mydi{\MM(\YY_i\urel,\Gamma_i)}
\def\mgwc{\MM(\WW,\Gamma)}

\def\myy{\MM(\YYo\adj\YYt,\eta)}

\let\mwg=\mgwc

\def\mgkst{\MM(\wn,\Gamma)^{st}}

\def\pri{^{\prime}}

\def\pr{{\rm pr}}
\def\Po{{\mathbf P}^1}
\def\pd{_{{\rm pd}}}

\def\rel{_{{\rm rel}}}
\def\rbb{]\!]}
\def\rig{^+}
\def\rhat{\,\hat{}\,}

\def\sub{\subset}
\def\sta{^{\ast}}

\def\unpo{^{n+1}}
\def\upmo{^{-1}}
\def\upvee{^{\vee}}

\def\upsigma{^{\sigma}}

\def\urel{^{\text{rel}}}
\def\uline{\underline}
\def\ucirc{^{\circ}}

\let\ver=\vet

\def\YYo{\YY_1\urel}
\def\YYt{\YY_2\urel}

\def\wn{{W[n]}}

\def\WWc{\WW}

\def\Zzd{{\Z}\urel}

\def\zn{{Z[n]}}

\def\lab{\ }
\def\lab{\label}

\title{Stable morphisms to singular schemes and relative stable morphisms}

\author{Jun Li}
\thanks{Supported partially by NSF grant, Sloan fellowship
and Terman fellowship}
\address{Department of Mathematics\\
Stanford University\\
Stanford, CA 94305 }
\date{}

\begin{document}

\begin{abstract}
Let $W/C$ be a degeneration of smooth varieties so that the special fiber has normal crossing
singularity. In this paper, we first construct the stack of expanded degenerations of $W$. We then
construct the moduli space of stable morphisms to this stack, which provides
a degeneration of the moduli spaces of stable morphisms associated to $W/C$.
Using a similar technique, for a pair $(Z,D)$ of smooth variety and a smooth divisor,
we construct the stack of expanded relative pairs and then the moduli spaces of relative stable morphisms
to $(Z,D)$. This is the algebro-geometric analogue of Donaldson-Floer theory in gauge theory.
The construction of relative Gromov-Witten invariants and the degeneration
formula of Gromov-Witten invariants will be treated in the subsequent paper.

\end{abstract}

\maketitle
\setcounter{section}{-1}

\section{Introduction}

The goal of this paper is to develop a theory of relative
Gromov-Witten invariants and prove a degeneration formula relating
the Gromov-Witten invariants of smooth varieties with the
relative Gromov-Witten invariants of pairs.

Let $C$ be a connected smooth curve, $0\in C$ be a fixed closed
point and let $W\to C$ be a family of projective schemes so that
the fibers $W_t$ of $W$ over $t\ne 0\in C$ are smooth and the
central fiber $W_0$ is the union of two smooth varieties $Y_1$
and $Y_2$ intersecting transversally along a smooth irreducible
divisor. We denote by $D_i\sub Y_i$ the divisor $Y_1\cap
Y_2\sub Y_i$. We will call $(Y_i,D_i)$ the decomposition pairs of
$W_0$. As for $t\ne 0$, since $W_t$ are members of a connected family of
smooth varieties, the Gromov-Witten invariants
of $W_t$ are all equivalent.
The question is how to relate the Gromov-Witten invariants
of $W_t$ with those of the pairs $(Y_i,D_i)$.
In this project, we will construct
the relative Gromov-Witten invariants of the pairs $(Y_i,D_i)$. We
then show that the Gromov-Witten invariants of $W_t$ can be
recovered from the relative Gromov-Witten invariants of the pairs
$(Y_1,D_1)$ and $(Y_2,D_2)$, via a formula of the form
$$GW(W_t)=GW(Y_1,D_1)\ast GW(Y_2,D_2).
$$
Here we use $GW(W_t)$ to denote the full GW-invariants of $W_t$
and use $GW(Y_i,D_i)$ to denote the full relative GW-invariants of
$(Y_i,D_i)$. The operator $\ast$ is an involution type product
linear in both arguments.

Gromov-Witten invariant are essentially a virtual enumeration of
algebraic curves in algebraic varieties (or pseudo-holomorphic
curves in symplectic manifolds). It is a topological theory.
Since their introduction, the Gromov-Witten invariants have become an invaluable tool in studying many
mathematical problems, not to mention their significance to mathematical physics.
A challenging problem of GW-invariants is the understanding of GW-invariants
of singular varieties and the behavior of GW-invariants under the degeneration of
target varieties. This line of research has been pursued by several groups.
In \cite{Tia}, Tian first studied the Gromov-Witten invariants of symplectic sums for
semi-positive symplectic manifolds and derived the decomposition formula of
the Gromov-Witten invariants in this setting.
Later, A.Li-Ruan \cite{LR, Ru3} worked out a version of degeneration formula of Gromov-Witten
invariants in the general setting, along the line of Donaldson-Floer theory.
A parallel theory was developed by
Ionel-Parker around the same time \cite{IP1,IP2,IP3}.
A part of the SFT theory of Eliashberg-Givental-Hofer \cite{EGH} can also be interpreted
as research along this line.
As is now well understood, one approach to degeneration of
moduli spaces (in differential geometry) is to follow Floer's original idea
in his cohomology theory.
This is the case for the
Donaldson-Floer theory, and is the case for Gromov-Witten invariants as demonstrated
by the works mentioned.

However, the algebro-geometric analogue
of this approach to degeneration of moduli spaces has eluded
algebraic geometers up until now. This in part explains why the
Donaldson-Floer theory was never developed fully in algebraic geometry.
This project fills in this gap. We have found the algebro-geometric
approach to Floer's original idea in treating degenerate objects in moduli spaces.
Though this project is mainly about moduli of stable morphisms and Gromov-Witten invariants,
the technique developed can be applied to many other moduli problems, including the moduli
of stable sheaves.

Out approach is geometric. We first construct a
degeneration of the moduli spaces of stable morphisms, associated
to the family $W/C$. Recall that to define the Gromov-Witten
invariants of a projective scheme $Z$, one needs to apply the
machinery of virtual moduli cycles, first developed by G.~Tian and
the author \cite{LT1,LT2}. (An alternative construction was achieved by
K.~Behrend and B.~Fantechi \cite{Beh, BF}. For related work using
analytic techniques, see \cite{FO, LT3, Ru1, Sie}.)
It turns out that the moduli of
stable morphisms to $W$, which can be viewed as a degeneration, is
not suitable for studying Gromov-Witten invariants.
Thus the first task is to construct a new degeneration in line
with the construction of virtual moduli cycles. This is achieved
in the present paper. In this paper, we first construct a stack
$\WW$ of expanded degenerations of $W/C$. We then investigate the
moduli of stable morphisms to this stack $\WW$ of expanded
degenerations. Similarly, to each pair $(Y,D)$ of a smooth divisor
in a smooth variety, we construct the stack of expanded relative
pairs $\YY\rel$ and the moduli space of relative stable morphisms
to the stack of expanded relative pairs $\YY\rel$. The main results
of this paper are the following existence theorems.

\begin{theo}
\label{0.1}
The moduli functor $\mgwc$ of stable morphisms to $\WWc$ of topological type $\Gamma$
is a separated and proper Deligne-Mumford stack over $C$.
\end{theo}

\begin{theo}
\label{0.2}
The moduli functor $\MM(\Zzd,\Gamma)$ of relative stable morphisms to $\Zzd$ of topological type $\Gamma$
is a separated and proper Deligne-Mumford stack.
\end{theo}


In the subsequent paper \cite{Li}, we will develop the obstruction
theory of the moduli stacks $\mwg$ and $\MM(\Zzd,\Gamma)$,
we will show that they admit perfect-obstruction theories
and thus have canonical virtual moduli cycles.
We will construct the relative Gromov-Witten invariants of
the pair $(Y,D)$ and prove the degeneration formula of the
Gromov-Witten invariants mentioned at the beginning of the introduction.
We will address the application of this machinery in our future research.

As indicated and/or shown by the works of \cite{Ion,IP1,LZZ,Ru2}, the degeneration machinery developed in this
project will be useful in enumerating curves in varieties, including algebraic
surfaces and Calabi-Yau manifolds. Some of these problems are crucial to the research in mathematical physics.
The moduli of relative stable morphisms and the relative GW-invariants will help answering
several outstanding conjectures raised by mathematical physicists \cite{GV}. The techniques developed in this
project can be applied to the study of other moduli problems, including stable sheaves over algebraic
varieties. Finally, we point out that based on this project the localization technique can be applied to the
moduli of relative stable morphisms, which should be useful as suggested by \cite{GV}.
We mention that the moduli space of relative stable morphisms and relative GW-invariants
were constructed for a restricted class of varieties for genus 0 curves \cite{Gat}.

We now briefly describe the main idea of this paper.
We let $W\to C$ be a projective family of schemes with general smooth fibers
$W_t$ for $t\ne 0\in C$ and singular fiber $W_0$
with normal crossing singularity, in the situation we mentioned at the
beginning of the introduction. Let $g,n$ be two integers and let $b$ be
a homology class.
We consider the moduli space (stack) $\MM_{g,k}(W_t,b)$ of stable morphisms from $k$-pointed
arithmetic genus $g$ curves to $W_t$ of degree $b$.
The union
$$\MM_{g,k}(W_{C\ucirc},b)=\cup_{t\in C\ucirc}\MM_{g,k}(W_t,b),\quad C\ucirc=C-0,
$$
is a proper family over $C\ucirc$.
Inserting the central fiber $\MM_{g,k}(W_0,b)$ gives one extension (which we
will call a degeneration of moduli spaces).
However, the natural obstruction theory of this new family is no longer perfect
near degenerate stable morphisms.
Here we say a stable morphism $f\mh X\to W_0$ is degenerate if some irreducible components
of $X$ are mapped entirely to the singular locus of $W_0$.
Usually, it is impossible to avoid degenerate stable morphisms if
the extension family is proper.
This makes this choice of degeneration inappropriate to study Gromov-Witten invariants.
In this paper, we will work out a new
degeneration of $\MM_{g,k}(W_{C\ucirc},b)$ that will allow us to
apply the machinery of virtual moduli cycles.

As a warming up to our construction, let us look at a simple case
where degenerate stable morphisms arise. Let $\pi\mh S\to C$ be
a flat morphism from a smooth curve to $C$ with $\pi\upmo(0)$
consists of a single point $s$. We
let $f\mh\cX\to W$ be a flat family of stable morphisms over $S$
with $f$ a $C$-morphism so that $f_s\mh \cX_s\to W_0$ is
degenerate. We now let $W\ucirc=W\times_C C\ucirc$, let
$S\ucirc=S-s$, let $\cX\ucirc=\cX\times_{C\ucirc}S\ucirc$ and let
$f\ucirc=f|_{\cX\ucirc}$ be the restriction of $f$ to over
$\cX\ucirc$. Here we view $W$ as a degeneration (extension) of the
family $W\ucirc$ to $C$ and view $f$ as an extension of $f\ucirc$.
The fact that the family $f$ is degenerate means that the
degeneration $W$ is not a suitable choice to extend $f\ucirc$.
There are other degenerations of $W\ucirc$. For instance, we let
$(\tilde 0,\tilde C)\to (0,C)$ be a base change ramified at
$\tilde 0\in\tilde C$ with ramification index $r$. We let $\tilde
W$ be the desingularization of the fiber product $W\times_C\tilde
C$. Note that the central fiber of $\tilde W$ over $\tilde 0$ has
$(r+1)$ irreducible components: $Y_1$, $Y_2$ and $(r-1)$ copies of
the ruled variety $\bP_{D_1}(\bone_{D_1}\oplus N_{D_1/Y_1})$ over
the divisor $D_1$. We now consider the new family $\tilde
f\ucirc\mh \tilde \cX\ucirc\to\tilde W$ over $\tilde
S\ucirc=S\times_C\tilde C$, where $\tilde\cX\ucirc=
\cX\ucirc\times_C\tilde C$. We can extend it to a family of stable
morphisms $\tilde f\mh \tilde \cX\to\tilde W$ over
$S\times_C\tilde C$. It is conceivable that for some choice of $r$ this new
extension $\tilde f$ will be non-degenerate. It is less obvious,
but will be proven, that there is a minimal choice of $r$ so that
$\tilde f$ is non-degenerate. Thus to avoid degenerate
stable morphisms we will replace the original family $W$ by
some expanded degenerations of $W\ucirc$. The exact choice of
the expanded degeneration, after imposing the minimality
condition, depends uniquely on the individual family.

Using expanded degenerations to study degenerate stable morphisms is certainly not new.
However, to obtain a moduli space that contains stable morphisms to various expanded
degenerations, namely to $\W[r]_0$ for various $r$,
we need to give the space of all expanded degenerations
an algebraic structure.
The new ingredient of this paper is that we can naturally
give such space the structure of an (Artin) stack. This is the stack $\WW$ mentioned before the
statement of Theorem \ref{0.1}.
The moduli space of stable morphisms to all expanded degenerations
of $W$ is then defined to be the moduli space of {\sl stable morphisms to this Artin stack $\WWc$}.
With appropriately defined stability, this moduli space is a proper Deligne-Mumford
stack.

It turns out that this is not sufficient to construct the desired degeneration of moduli spaces.
We need to build the notion of {\sl pre-deformable} morphisms into the definition of the
stable morphisms to $\WWc$.
Let $f\mh X\to W_0$ be a non-degenerate stable morphism so that there is a smooth point $p\in X$
that is mapped to the singular locus $D$ of $W_0$. Then a simple argument shows that
there is no small deformation of $f$ that moves it away from the singular fiber $W_0$.
In other words, there is a local obstruction to deforming
such morphisms away from the singular fiber $W_0$.
As we are interested in the degeneration of the moduli spaces, we
should exclude such morphisms from our total moduli space. To accommodate this we
require all stable morphisms to $\WWc$ to have vanishing
local obstructions to deforming to smooth fibers of $W/C$.
Those morphisms that have vanishing local obstructions will be
called {\sl pre-deformable} morphisms.

In the end, we define stable morphisms to $\WWc$ of given topological
type\footnote{$\Gamma$ contains all
the topological restrictions on the stable morphisms,
including the genus, the number of marked points and the degree, etc.}
$\Gamma$
to be those that are non-degenerate, are minimal among all possible target schemes $\wn_0$ and are pre-deformable.
The main theorem of this paper is that the moduli of such stable morphisms form a
separate and proper Deligne-Mumford stack over $C$, as stated in Theorem \ref{0.1}.

As I mentioned earlier in the introduction, this construction
can be viewed as an algebro-geometric adaptation of
Floer's pioneering work \cite{Flo}.
I should point out that in early 1980's, Gieseker \cite{Gie} (see also Gieseker-Morrison \cite{GM})
constructed a degeneration of the moduli spaces of stable vector bundles over a family of smooth curves
degenerating to a nodal curve. The moduli space of vector bundles over the singular curve $D$
consists of the stable vector bundles over $D$ and the stable
vector bundles over $\tilde D$ that are the result of replacing the node of $D$ by a rational curve.
Bundles over $\tilde D$ are exactly the substitute of those stable non-locally free
sheaves over $D$. The current work can be viewed as a realization of Gieseker's construction to
the general cases.
I should also mention that the notion of admissible cover introduced by Harris-Mumford \cite{HM} in early 80's
utilized similar idea to deal with singular objects.
A more recent work related to this construction is Caporaso's construction of the universal Picard
scheme over $\overline{M}_{g}$ \cite{Cap}.

This technique can be applied to construct the {\sl relative stable morphisms} of pairs.
Let $(D,Z)$ be a pair consisting of a smooth divisor $D$ in a smooth projective variety $Z$.
We consider the moduli of stable morphisms $f\mh X\to Z$. We call a stable morphism $f\mh X\to Z$
non-degenerate (relative to the divisor $D$) if $f\upmo(D)$ is a proper divisor in $X$ and is away
from the nodal and the marked points of $X$.
As before, the issue is to construct a complete moduli space that contains no
degenerate stable morphisms relative the to $D$. As for $W$, we should consider stable
morphisms to expanded relative pairs.
An expanded relative pair of $(Z,D)$ is a reducible scheme that is the union of $Z$ with several copies of the ruled variety
$\bP_D(1_D\oplus N_{D/Z})$.
To define the moduli of relative stable
morphisms we first construct the (Artin) stack $\Zzd$ of all expanded relative pairs.
We then introduce the notion of relative stable morphisms to $\Zzd$, similar to the case of stable morphisms
to $\WWc$. After that we can define the moduli functor of relative stable morphisms to $\Zzd$ of topological type $\Gamma$.
In the end we prove the existence Theorem \ref{0.2}.

The paper is organized as follows. In section one we will construct the stack $\WWc$ of all expanded degenerations
of $W$. The notion of pre-deformable morphisms and related properties will be worked out in section two.
In section three we will define the notion of stable morphisms to $\WWc$ and prove that the moduli functor
of stable morphisms to $\WWc$ is a separated and proper algebraic stack over $C$. Section four is devoted
to study the moduli of relative stable morphisms $\MM(\Zzd,\Gamma)$. The relation of stable morphisms to
$\WW_0=\WWc\times_C 0$ and relative stable morphisms to pairs $(Y_i,D)$ will be investigated there as well.

\subsection{Conventions}

We provide a selected list of conventions used throughout this paper.

\begin{description}
\item[$\wn$] It is the expanded degeneration (of $W/C$) constructed in section 2.

\item[$\cn$] $\cn=C\times_{\Ao}\Anpo$. It has two tautological morphisms $\cn\to\Ao$ and
$\cn\to\Anpo$. It is defined in section 2.

\item[$\gn$] It is $\cong G_m\times\cdots\times G_m$, $n$ copies. For notation related to this group
see (\ref{1.127}) and (\ref{1.126}).

\item[$(Y_i,D_i)$] They are the pairs (of irreducible components with singular locus)
derived from the decomposition of $W_0$.

\item[$X_1\adj X_2$] This is the gluing of two schemes $X_1$ and $X_2$. Let $A_1\sub X_1$ and $A_2\sub X_2$ be
two pairs of closed subschemes in schemes. Assume $A_1$ is isomorphic to $A_2$. Then we can define
a new scheme, denoted by $X_1\adj X_2$, that is the result of gluing $X_1$ and $X_2$
along $A_1$ and $A_2$. $X_1\adj X_2$ is the scheme that satisfies the universal property of the push-out.


\item[$\cn_{[l\dual]}$] $\bH_l\unpo$ is the coordinate hyperplane
given by the vanishing of the $l$-th coordinate axis of $\Anpo$.
We will call it the $l$-th coordinate hyperplane. $\cn_{[l\dual]}=\cn\times_{\Anpo}\bH_l\unpo$.

\item[$\bt$] $\bt=(t_1,\cdots,t\lnpo)$ and hence $\kk[\bt]=\kk[t_1,\cdots,t\lnpo]$.

\item[\hbox{$[n]$}] This is the set of all integers between $1$ and $n$.

\end{description}

\section{The stack of expanded degeneration}

We first fix the notation that will be used throughout this paper.
In this paper, we will work with an algebraically closed field $\kk$ with characteristic $0$.
We let $C$ be a smooth irreducible curve, $0\in C$ a closed point and
$\pi\mh W\to C$ a flat and projective family of schemes satisfying the following
condition: The morphism $\pi$ is smooth away from the central fiber $W_0=W\times_C0$
and the central fiber $W_0$ is
reducible with normal crossing singularity and has
two smooth irreducible components $Y_1$ and $Y_2$ intersecting
along a smooth divisor $D\sub W_0$. When we
view $D$ as a divisor in $Y_i$, we will denote it by $D_i\sub Y_i$.

We now construct the class of expanded degenerations of $W$ mentioned in the introduction of this
paper. We let $\Delta$ be the projective bundle over $D$:
$$\Delta=\bP(\bone_D\oplus N_{D_2/Y_2}),
$$
where $\bone_D$ is the trivial holomorphic line bundle on $D$
and $N_{D_2/Y_2}$ is the normal bundle of $D_2$ in $Y_2$,
viewed as a line bundle on $D$. $\Delta$ has
two distinguished divisors
$$D_-=\bP(\bone_D\oplus 0)\quad {\rm and}\quad
D_+=\bP(0\oplus N_{D_2/Y_2}).
$$
For convenience, we call $D_-$ the left distinguished and $D_+$ the right distinguished divisors.
Of course both $D_-$ and $D_+$ are canonically isomorphic to $D$.
Note that
$$N_{D_-/\Delta}\cong N_{D_2/Y_2}\quad
{\rm and}\quad N_{D_+/\Delta}\cong N_{D_1/Y_1}\cong N_{D_2/Y_2}\upmo.
$$
Using this identification, we can {\sl glue}\footnote{See convention in the Introduction.}
an ordered chain of $n$ $\Delta$'s by
identifying the right distinguished divisor (i.e. $D_+$) in the $k$-th $\Delta$ with the left distinguished
divisor (i.e. $D_-$) of the $(k+1)$-th $\Delta$, for $k=1,\cdots, n-1$. We denote the
resulting scheme by $\Delta[n]$. It is connected, has normal crossing singularity and has
$n$-irreducible component all isomorphic to $\Delta$. We keep the ordering of these $n$ $\Delta$'s.
We then {\sl glue} $Y_1$
to $\Delta[n]$ by identifying $D_1$ in $Y_1$ with the left distinguished divisor $D_-$ in the first
$\Delta$ of $\Delta[n]$, and then {\sl glue} $Y_2$ to this scheme by identifying
the right distinguished divisor $D_+$ of the last $\Delta$ in $\Delta[n]$ with $D_2$ in $Y_2$.
We denote the resulting scheme by $W[n]_0$. Note that $W[n]_0$ has
$(n+2)$-irreducible components.
These $(n+2)$-components form a chain ordered from left to right according to their intersection pattern:
$$W[n]_0=Y_1\adj \Delta\adj\cdots\adj\Delta\adj Y_2.
$$
Often, we will denote these $n+2$ components by
$\Delta_1,\cdots,\Delta_{n+2}$, according to this ordering.
As a result, $W[n]_0$ has $(n+1)$-nodal divisor, ordered so that the $k$-th nodal divisor is
$D_k\defeq \Delta_k\cap\Delta_{k+1}$.
Note that $D_k\sub \Delta_k$ (resp. $D_k\sub\Delta_{k+1}$) is the right (resp. left) distinguished
divisor.
We agree $W[0]_0$ is $W_0$, the fiber of $W$ over $0\in C$.

The purpose of this paper is to construct a stack $\WWc$ representing all degenerations $\tilde W\to\tilde C$,
where $\rho\mh \tilde C\to C$ are base changes, so that the fiber of $\tilde W$ over $t$ is
$W_{\rho(t)}$ in case $\rho(t)\ne 0$ and is one of $\wn_0$ when $\rho(t)=0$.

\subsection{Construction of the standard model}

We begin with the construction of the standard models $\wn$.
By replacing $C$ with an open neighborhood of $0\in C$ we can assume
that there is an \'etale morphism
$C\to\Ao$ so that $0\in C$ is mapped to $0_\Ao\in \Ao$ and
$0\in C$ is the only point in $C$ lies over $0_{\Ao}\in \Ao$.
We fix such a map $C\to\Ao$ once and for all.
We let
$G_m$ be the general linear group GL(1) and let $G[n]$ be the group scheme
\begin{equation}
\label{1.125}
\gn=G_m\times\cdots\times G_m\qquad
n\ {\rm copies}
\end{equation}
which acts on $\Anpo$ via
\begin{equation}
\label{1.127}
\bt^{\sigma}=(\sigma_1t_1,\sigma_1\upmo\sigma_2t_2,\cdots,
\sigma_{n-1}\upmo\sigma_{n}t_{n},\sigma\upmo_{n}t_{n+1}),
\end{equation}
where $\bt=(t_1,\cdots,t\lnpo)$ is the standard coordinate of $\Anpo$.
(Here and in the following we will use superscript to denote the result of the
group action.) The convention throughout this paper is that for $\sigma\in\gn$ we will
denote by $\sigma_i$ its $i$-th component via isomorphism (\ref{1.125}). For convenience, we introduce
\begin{equation}
\label{1.126}
\bar\sigma_i=\sigma_i/\sigma_{i-1}\quad \text{with}\quad \sigma_0=\sigma\lnpo=1.
\end{equation}
In this way, (\ref{1.127}) can be rewritten as
$\bt^{\sigma}=(\bar\sigma_1t_1,\cdots,\bar\sigma\lnpo t\lnpo)$.
Note that if we let
$$\bp: \Anpo\to \Ao, \quad \bp(\bt)= t_1\times\cdots\times t\lnpo,
$$
then $\bp$ is
$\gn$-equivariant with the trivial $\gn$-action on $\Ao$. We define
$$\cn=C\times_{\Ao}\Anpo
$$
with the (unique) distinguished point $\bz=0\times_{\Ao}0_{\Anpo}\in \cn$.
We view $\cn$ as a $C$-scheme via the first projection.
$\cn$ is a $\gn$-scheme as well.

The standard model $\wn$ will be constructed as a desingularization of
$$ W\times_{\Ao}\Anpo=W\times_C\cn.
$$

We begin with some more notation. For any integer $n$, we denote by $[n]$
the set of integers between $1$ and $n$. Let $I\sub[n+1]$ be any subset
with cardinality $|I|=m+1$. The subset $I$ defines a unique increasing map $[m+1]\to[n+1]$ whose image
is $I$. By abuse of notation, we will denote this map by $I$. Hence $I(k)$ is the
$k$-th element in $I$. Given any such $I$ there is a standard embedding
$\Ampo\to \Anpo$ that sends $z\in\Ampo$ to $w\in\Anpo$ via $w_{I(k)}=z_k$ for $k\in I$ and $w_l=1$
for $l\not\in I$. It follows that it induces an immersion
\begin{equation}
\label{1.213}
\underline{\gamma_I}:C[m]\to C[n].
\end{equation}
We call this the {\sl standard embedding} associated to $I\sub [n+1]$.
There is another embedding $\Ampo\to\Anpo$ via the rule $w_{I(k)}=z_k$ for $k\in I$ and $w_l=0$
for $l\not\in I$. This induces an immersion
$\Ampo\to\cn$ whose image lies in a coordinate plane.
We call this the {\sl coordinate plane} associated to $I$ and denoted
it by $\cn_I$. Lastly, when $l$
is an integer in $[n+1]$, we denote by $l$ the set of one element with $l\to [n+1]$ the inclusion.
Hence $\cn_{l}$ is the $l$-th coordinate axis of $\cn$.
Similarly, we let $[l^{\vee}]$ be the complement of $l$ in $[n]$. Hence $C[n]_{[l\upvee]}$
is the coordinate hyperplane spanned by all coordinate axes except the $l$-th axis.

Now we construct $\wn$ by induction on $n$. For $n=0$, $W[0]=W$.
For $n=1$, $W_1=W\times_{C}C[1]$ is smooth except along the
locus $D\times_{C} \bz$, which is a smooth codimension 3 subscheme in $W_1$
and the formal completions (the germs) of its normal slices
in $W_1$ is isomorphic to the formal completion of
\begin{equation}
\label{1.1}
X=\{z_1z_2=t_1t_2\}\sub \bA^4
\end{equation}
along its origin. Here we use $(z_1,z_2,t_1,t_2)$ to denote the coordinate of $\bA^{\! 4}$.
After blowing up $W_1$ along $D\times_C\bzero$ we obtain a smooth scheme
$\tilde W_1$ whose exceptional divisor is a $\Po\times\Po$ bundle over $D$. In
the following, we will show how to contract one of these two
$\Po$-factors from this divisor to
obtain the desired scheme $W[1]$.

We first work out the a desingularization of (\ref{1.1}).

\begin{lemm}
\label{A1}
There is a desingularization $Z$ of (\ref{1.1}) so that its exceptional locus
$\Sigma$ is a $\Po$ and the proper transforms of the $z_1$ and the $t_2$-axis
(resp. the $z_2$ and the $t_1$-axis) in $Z$
intersect in $\Sigma$.
\end{lemm}

\begin{proof}
We first blow up the singular point (the origin)
of the threefold $X$ in (\ref{1.1}).
The resulting threefold, denoted by $\tilde X$, is isomorphic to the total space of the restriction
of the tautological line bundle on $\bP^3$ to the surface
$\{w_1w_2=w_3w_4\}\sub \bP^3$.
The exceptional divisor $B\sub \tilde X$ is isomorphic to the zero section of this line bundle,
which is isomorphic to
$\Po\times\Po$. An isomorphism
$$\Po\times\Po\lra B=\{w_1w_2=w_3w_4\}\sub\bP^3
$$
is given by
\begin{equation}
([a_0,a_1],[b_0,b_1])\mapsto[a_0b_1,a_1b_0,a_0b_0,a_1b_1].
\label{1.2}
\end{equation}
The normal bundle $N_{B/\tilde X}$ is isomorphic to the restriction of the tautological line bundle
to $B$, which is $\cO_{\Po\times\Po}(-1,-1)$, under the above isomorphism. Hence we can contract
either of the $\Po$-factor in $B$ to obtain a smooth threefold.
In the following we will show that we can choose the contraction so that
the proper transforms of the $t_1$ and the $z_2$-axes intersect and the proper transforms
of the $t_2$ and the $z_1$-axes intersect.

Let $\tilde p\mh \tilde X\to X$ be the obvious projection.
Then under the isomorphism (\ref{1.2}) the map $\tilde p$ is of the form
$$\tilde p\mh ([a_0,a_1],[b_0,b_1],\zeta)\mapsto
(z_1,z_2,t_1,t_2)=
(a_0b_1\zeta, a_1b_0\zeta, a_0b_0\zeta, a_1b_1\zeta),
$$
where $\zeta$ is the variable representing points on the total space of the tautological
line bundle. Note that $\zeta$ has homogeneous degree $(-1,-1)$, as $a_i$ and $b_i$ have
homogeneous degree $(1,0)$ and $(0,1)$ respectively\footnote{
By this we mean $([a_0,a_1],[b_0,b_1],\zeta)\sim
([\lambda a_0,\lambda a_1],[\mu b_0,\mu b_1],\lambda\upmo\mu\upmo\zeta)$.}.
Hence the proper transforms of the $z_1$, $z_2$, $t_1$ and the $t_2$-axes intersect $B$ at
$$([1,0],[1,0]),\ ([0,1],[0,1]),\ ([1,0],[0,1])\ {\rm and}\ ([0,1],[1,0])\in\Po\times\Po,
$$
respectively.
Therefore if we contract the first $\Po$-factor the resulting threefold $Z$ satisfies the
required property.
\end{proof}

For later application, we now work out an atlas of $Z$. By out construction,
$Z$ is the total space of
the vector bundle associated to the locally free sheaf $\cO_{\Po}(-1)^{\oplus 2}$.
Since $Z$ is derived from $\tilde X$ by contracting the first $\Po$ factor,
we will use $[b_0,b_1]$ as the homogeneous coordinate of the exceptional locus of
$p\mh Z\to X$. In this way, $p\mh Z\to X$ is of the form
\begin{equation}
\label{1.61}
p:([b_0,b_1],\eta_1,\eta_2)\mapsto (b_1\eta_1,b_0\eta_2, b_0\eta_1, b_1\eta_2)
\end{equation}
and the contraction $\tilde X\to Z$ is given by
\begin{equation}
([a_0,a_1],[b_0,b_1],\zeta)\mapsto ([b_0,b_1],a_0\zeta,a_1\zeta).
\end{equation}
Hence the projection $Z\to\bA^{\!2}$ is given by
\begin{equation}
\label{1.62}
\varphi: ([b_0,b_1],\eta_1,\eta_2)\mapsto (t_1,t_2)=(b_0\eta_1,b_1\eta_2).
\end{equation}
From this we see that $Z\times_{\At}0_{\At}$
is a nodal curve obtained by adjoining an $\Ao$ (with $\eta_1$ variable) to $[0,1]$
in $\Po$ and an $\Ao$ (with $\eta_2$ variable) to $[1,0]$ in $\Po$.
Further if we let $\bA_1$ (resp. $\bA_2$; resp.
$S$) be the $t_1$-axis of $\At$ (resp. the $t_2$-axis; resp.
a curve in $\At$ passing though the origin but not entirely contained in the two coordinate
lines) then $Z\times_{\At}\bA_1$ (resp. $Z\times_{\At}\bA_2$; resp. $A\times_{\At}S$)
is a smoothing of the nodal point $[0,1]$ (resp. the nodal point $[1,0]$;
resp. both nodal points) of $Z\times_{\At}0_{\At}$. We also mention that near $([0,1],0,0)$
(resp. $([1,0],0,0)$) the projection $Z\to\At$ has the form $(b_0,\eta_1,\eta_2)\mapsto
(b_0\eta_1,\eta_2)$ (resp. $(b_1,\eta_1,\eta_2)\mapsto(\eta_1,b_1\eta_2)$).

We continue our construction of $W[1]$. After blowing up the subscheme $D\times_{C}\bz$
in $W\times_C C[1]$
we obtain the scheme $\tilde W_1$ whose exceptional divisor is isomorphic to a $\Po\times\Po$-bundle
over $D$. By Lemma \ref{A1}, we can contract either factor of this $\Po\times\Po$-bundle to obtain a
smooth scheme. We now specify the choice of the $\Po$-factor to be contracted.
Let $W[1]$ be the resulting scheme after contracting one such $\Po$-factor.
Then $W[1]$ has an induced projection
$\pi\mh W[1]\to C[1]$
and the fiber of $\pi$ over the origin $\bz\in C[1]$, denoted by $W[1]_0$, is
reduced with normal crossing singularity whose irreducible components are $Y_1$, $Y_2$ and
a ruled variety $\tilde \Delta$ over $D$. The singular locus of $W[1]_0$ consists of two
disjoint divisors $Y_1\cap\tilde\Delta$ and $Y_2\cap\tilde\Delta$, both isomorphic to $D$.
Following the proof of
Lemma \ref{A1} and the description afterwards, we can contract one of
the $\Po$ factor so that the family
$W[1]\times_{C[1]}C[1]_{1}$ ( resp. $W[1]\times_{C[1]}C[1]_{2}$)
over $C$ is a smoothing
of the nodal divisor $Y_1\cap\tilde\Delta$ (resp. the nodal divisor
$Y_2\cap\tilde\Delta$) of $W[1]\times_{\At}0_{\At}$.
Here $C[1]_i\sub C[1]$ is the $i$-th coordinate line defined in the beginning of this subsection.
It is direct to check that $\tilde\Delta$ is the ruled variety $\bP_D(\bone_D\oplus N_{D_1/Y_1})$
mentioned at the beginning of this section.

We let $\pi_2\mh W[1]\to\Ao$ be the composite of $W[1]\to C[1]\to\At$ with the second
projection $\At\to\Ao$. It is easy to see that the singular locus of $\pi_2$
is the proper transform of $D\times\Ao\sub W\times_{\Ao}\At$ where $D\times\Ao\to
W\times_{\Ao}\At$ is the morphism induced by $D\to W$ and $\Ao\to\At$
as the $t_1$-axis. We denote the singular locus of $\pi_2$ by $\bD_2$.
Further, at each $p\in\bD_2$, there is a formal coordinate chart\footnote{
Namely they generate the maximal ideal of the formal ring $\hat\cO_{W[1],p}$.}
$(z_1,\cdots)$ of $p$ so that $\pi\sta\mh \Gamma(\At)\to \hat \cO_{W[1],p}$
satisfies $t_1\mapsto z_1$ and $t_2\mapsto z_2z_3$.

We now construct $W[2]$. We let $W_2$ be the fiber product $W[1]\times_{\At}\As$ with $\As\to\At$
the morphism $(t_1,t_2,t_3)\mapsto(t_1,t_2t_3)$.
Then $W_2$ has singularity along $\bD_2\times_{\At}\As_{[23\dual]}\sub W[1]\times_{\At}\As$,
where following our convention $\As_{[23\dual]}$ is the codimension two coordinate plane
(line) spanned by all but the second and the third axes in $\As$.
Note that the singularity type of $W_2$ along $\bD_2\times_{\At}\As_{[23\dual]}$
is identical to the one of $W_2$ along $D\times_{\Ao}0_{\Ao}$.
We now blow up $W_2$ along its singular locus and then contract one copy of the $\Po$-bundle
in the resulting exceptional divisor. We will
contract the $\Po$-bundle so that the resulting scheme, denoted by $W[2]\to C[2]$, has the property
that its fibers over the $l$-th axis $C[2]_l\sub C[2]$ is a smoothing of the $l$-th nodal divisor of $W[2]_0\defeq
W[2]\times_{C[2]}\bz$. The proof of this is a word by word repetition
of our earlier construction of $W[1]$ from $W\times_{\Ao}\At$.

The model $\wn/\cn$ is constructed by induction. Assume $W[n-1]/C[n-1]$ was constructed,
we then form $W_n=W[n-1]\times_{\An}\Anpo$, where $\Anpo\to\An$ is defined by $(\cdots,t_{n+2})\mapsto
(\cdots,t_n,t_{n+1}t_{n+2})$. We then resolve the singularity of $W_n$ by a small resolution.
We require that the resulting family $\wn/\cn$ be so that
$\wn\times_{\cn}\bz\cong\wn_0$ and the fibers of $\wn$ over the $l$-th factor $\cn_l\sub\cn$
is a smoothing of the $l$-th nodal divisor of $\wn_0$. Again the detail to work this out is exactly the same as
the construction of $W[1]$.

Throughout this paper, we will use $\wn_0$ to denote the fiber of $\wn$ over $\bz\in\cn$
and denote by $\Delta_1,\cdots,\Delta_{n+2}$ the $n+2$ irreducible components of $\wn_0$,
ordered as specified at the beginning of this section.

\subsection{Construction of the stack of expanded degenerations}

In this part, we will construct the stack $\WWc$ of the expanded degenerations of $W$.

We begin with an investigation of the $\gn$-group action on $\wn$.
Clearly, since $\Anpo\to\Ao$ is $\gn$-equivariant, its action lifts to $\cn\to C$.
To show that it lifts to $\wn$, we need to understand the
corresponding $G_m$-action on the
desingularization $Z$ of $X$, constructed in the subsection 1.1. Following
the notation there, $G_m$ acts on $X$ via
$$(z_1,z_2,t_1,t_2)\upsigma=(z_1,z_2,\sigma t_1,\sigma\upmo t_2),\quad \sigma\in G_m
$$
and its lifting to $Z$ is
\begin{equation}
\label{1.6}
([b_0,1],\eta_1,\eta_2)\upsigma=([\sigma b_0,1],\eta_1,\sigma\upmo\eta_2).
\end{equation}
Based on this, it is easy see that the $\gn$-action on $\cn$ lifts uniquely to $\wn$.

For later application, we now give a precise description of the $\gn$-action
near the nodal divisor of $\wn_0$. This can best be done by looking at the simplest
case where $W=\At$, $C=\Ao$ and $W\to C$ is given by $t=u_1u_2$, where $(u_1,u_2)$
and $t$ are the standard coordinates of $\At$ and $\Ao$.
With this choice of $W\to C$, we can construct the corresponding $\wn\to\cn$.
To distinguish this from the general case, we will denote this special degeneration
by $\Gn\to\Anpo$.

\begin{lemm}
\label{1.63}
Let the notation be as before. Then $\Gn_0$ is a chain of $n+2$ curves of
which the first and the last are isomorphic to $\Ao$ and the remainders are isomorphic to $\Po$.
Further, we can cover $\Gn$ by $\gn$-invariant affine open subsets $U_1,\cdots,U\lnpo$,
each isomorphic to $\Anpt$, so that if we denote by $(u_1^l,\cdots,u_{n+2}^l)$
the standard coordinate of $U_l\cong\Anpt$, then
\newline
(i) the restriction of $\pi\mh\Gn\to\Anpo$ to
$U_l$ is given by
$$\pi^l:(u^l_{1},\cdots,u^l_{n+2})\mapsto(\cdots,u^l_{l-2},u^l_{l-1},u^l_{l}u^l_{l+1},u^l_{l+2},
u^l_{l+3},\cdots);
$$
\newline
(ii) The $\gn$-action on $U_l$ is given by
$$(u^l_{1},\cdots,u^l_{n+2})\upsigma=
(\cdots,\bar\sigma_{l-2}u^l_{l-2},
\bar\sigma_{l-1}u^l_{l-1}, \sigma_{l-1}\upmo u^l_{l},\sigma_l u^l_{l+1},
\bar\sigma_{l+1}u^l_{l+2}, \bar\sigma_{l+2}u^l_{l+3},\cdots),
$$
where $\bar\sigma_i=\sigma_i/\sigma_{i-1}$ with
$\sigma_i=1$ for $i=0, n+1$;
\newline
(iii) The transition function from $(u^l_{\cdot})$ to $(u^{l+1}_{\cdot})$ over $U_l\cap U_{l+1}$
is given by
$$(u^{l+1}_{1},\cdots,u^{l+1}_{n+2})=
(\cdots,u^l_{l-2},u^l_{l-1},u^l_{l}u^l_{l+1},1/u^l_{l+1},u^l_{l+2}u^l_{l+1},
u^l_{l+3},\cdots).
$$
\end{lemm}

\begin{proof}
We prove this Lemma by induction on $n$. For $n=0$, there is nothing to prove. We now assume that the Lemma
holds for $\Gamma[n-1]$.
Namely, we have found coverings $U_1,\cdots,U_n$ of $\Gamma[n-1]$ that satisfy the required property.
Following the construction, $\Gamma[n]$ is a small resolution of
$\Gamma[n-1]\times_{\An}\Anpo$. Using the explicit description of $\pi|_{U_l}\mh U_l\to\An$, we see that
$$V_l\defeq U_l\times_{\An}\Anpo\cong \bA^{\!n+2},\quad\text{where}\ l<n,
$$
are smooth $G[n]$-invariant open subsets of $\Gamma[n]$.
Clearly, we can choose a unique coordinate $(v_{\cdot}^l)$ of
$V_l$ so that the projection $V_l\to U_l$ is given by
$$u^l_{i}=v^l_{i},\quad i\leq n; \quad u^l_{n+1}=v^l_{n+1}v^l_{n+2},
$$
and the projection $\pi^l\mh V_i\to \Anpo$ is given by the formula ($i$)
in the statement of the Lemma.
A routine check shows that for $l<n$ the $G[n]$-action
on $V_l$ and the transition function from $V_{l-1}$ to $V_l$ is exactly as shown in the statement of
the Lemma.

It remains to check the case where $l=n$ and $n+1$. We now look at the chart $U_n$. Let
$\varphi: \Gamma[n]\to\Gamma[n-1]$
be the composite of the desingularization $\Gamma[n]\to\Gamma[n-1]\times_{\An}\Anpo$
with the obvious projection. Then
$$\varphi\upmo(U_n)\cong \bA^{\!n-1}\times \bL^{\oplus 2},
$$
where $L$ is the
tautological line bundle (degree $-1$) on $\Po$ and $\bL^{\oplus 2}$ is the total space
of $L^{\oplus 2}$.
We now cover $\Po$ by $B_0=\{[b_0,1]\}$ and $B_1=\{[1,b_1]\}\sub\Po$. We
choose trivializations $\bL^{\oplus 2}|_{B_0}\cong B_0\times\At$ and
$\bL^{\oplus 2}|_{B_1}\cong B_1\times\At$ so that the transition function from
$B_0\times\At$ to $B_1\times\At$ is given by
$$(b_0,\zeta_1,\zeta_2)\lra (b_1,\zeta\pri_1,\zeta\pri_2)=(1/b_0,\zeta_1b_0,\zeta_2b_0).
$$
We let $V_n$ be $\Anmo\times \bL^{\oplus 2}|_{B_0}$ and let
$V_{n+1}=\Anmo\times \bL^{\oplus 2}|_{B_1}\sub \Anmo\times \bL^{\oplus 2}$.
Then $V_1,\cdots,V_{n+1}$ form a covering of $\Gamma[n]$. Next for $V_n$ we choose
a coordinate
$$(v_1^n,\cdots,v_{n+2}^n)=(z_1,\cdots,z_{n-1}, \zeta_1,b_0,\zeta_2)
$$
and for $V_{n+1}$ we choose
$$(v_1^{n+1},\cdots,v_{n+2}^{n+1})=(z_1,\cdots,z_{n-1}, \zeta_1\pri,b_1,\zeta_2\pri).
$$
We now check that this choice of coordinates satisfy the requirement of the Lemma.

We first check the transition functions. The transition function from $V_n$ to $V_{n+1}$
is determined by the transition function from $\bL^{\oplus 2}|_{B_0}$ to
$\bL^{\oplus 2}|_{B_1}$. Thus we have
$$(v_1^{n+1},\cdots,v_{n+2}^{n+1})=(z_1,\cdots,z_{n-1}, \zeta_1\pri,b_1,\zeta_2\pri)
=(z_1,\cdots,z_{n-1}, b_0\zeta_1,1/b_0,\zeta_2 b_0)
$$
$$
\quad\ =(v_1^n,\cdots,v_{n-1}^n, v_n^nv_{n+1}^n,1/v_{n+1}^n, v_{n+1}^n v_{n+2}^n).
$$
This is the requirement in ($iii$) for $l=n$.

It remains to work out the transition function from $V_{n-1}$ to $V_n$.
We let $\bv^{n-1}=(v_1^{n-1},\cdots,v_{n+2}^{n-1})$ be a general point and let
$\bv^n=(v_1^n,\cdots,v_{n+2}^n)$ be the point defined by the relation
($iii$) with $u$ replaced by $v$ and $l$ replaced by $n-1$. Then
$$v_{n+2}^n=v_{n+2}^{n-1},\ v_{n+1}^n=v_{n+1}^{n-1}v_n^{n-1},\
v_n^n=1/v_n^{n-1},\ v_{n-1}^n=v_{n-1}^{n-1}v_n^{n-1}.
$$
Clearly, to prove ($iii$) we only need to check that the images of $\bv^{n-1}$ and $\bv^n$
in $U_n$ coincide and their images in $\Anpo$ also coincide.
By definition, the image of $\bv^{n-1}$ in $\Anpo$ is
$$(v_1^{n-1},\cdots,v_{n-1}^{n-1}v_n^{n-1},v_{n+1}^{n-1},v_{n+2}^{n-1})
$$
and the image of $\bv^n$ in $\Anpo$ is
$$(v_1^n,\cdots,v_{n-1}^n, v_n^nv_{n+1}^n, v_{n+2}^n).
$$
They coincide using the relations mentioned above.
On the other hand the image of $\bv^{n-1}$ in $U_{n-1}$ is
$(\cdots,v_n^{n-1}, v_{n+1}^{n-1}v_{n+2}^{n-1})$ and its image in $U_n$,
using ($iii$) and the induction hypothesis, is
$$(v^{n-1}_1,\cdots,v_{n-1}^{n-1}v_n^{n-1},1/v_n^{n-1}, v_n^{n-1}v_{n+1}^{n-1}v_{n+2}^{n-1}).
$$
The image of $\bv^n$ in $U_n$ is
$$(v_1^n,\cdots,v_{n-1}^n,v_n^n,v_{n+1}^nv_{n+2}^n).
$$
They coincide as well. This proves the relation ($iii$).
The part ($ii$) concerning the group action can be
checked directly and will be omitted.
\end{proof}

To apply this Lemma to the general case $\wn$, we have the following useful observation.
Let $\hat D$ be the formal completion of $W$ along $D$. For any open $V\sub D$, we
let $\hat V$ be $V\sub\hat D$ endowed with the open subscheme structure.
We choose $V$ so that $\cO_{\hat V}\cong\cO_V\lbb w_1,w_2\rbb$
and the induced morphism $\hat V\to\Ao$ is given by $t\mapsto w_1w_2$.
We let $\tilde V\sub \wn$ be an open subset so that
$\tilde V\cap(\wn\times_WD)=\wn\times_WV$ and let
$\hat\wn$ be the formal completion of $\tilde V$ along $\wn\times_W V$.
Similarly we let $U\sub\Gn$ be the preimage of $0_{\At}\in\At$ under $\Gn\to\Gamma[1]=\At$
and let $\hat\Gn$ be the formal completion of $\Gn$ along $U$.
Then we have a $\gn$-equivariant isomorphism
\begin{equation}
\label{1.67}
\hat\wn\cong V\times \hat\Gn.
\end{equation}
We also denote by
\begin{equation}
\label{1.51}
\lbPsi: \gn\times\cn\lra \cn\quad
{\rm and}\quad
\Psi:\gn\times\wn\lra\wn
\end{equation}
the morphisms of $\gn$-actions.

\begin{coro}
\label{1.111}
Let $\gn_l$ be the subgroup of $\gn$ that is the $l$-th copy of $G_m$ in $\gn$, using the
isomorphism (\ref{1.125}). Then the $\gn_l$
action on $\Delta_{l+1}$ is induced by the linear action of $G_m$ on
$\bone_D$ and $N_{D_2/Y_2}$ of weight $0$ and $1$, respectively.
It acts trivially on all other $\Delta_i$'s.
\end{coro}

\begin{coro}
\label{1.8}
Let $h\mh \wn_0\to\wn_0$ be an isomorphism commuting with $\wn_0\to W_0$.
Then there is a $\sigma\in\gn$ so that $h$ is identical to the automorphism
of $\wn_0$ induced by the action of the $\sigma$. Namely $h(w)=\Psi(\sigma,w)$ for and
$w\in \wn_0$.
\end{coro}

\begin{coro}
\label{1.9}
Let $s\in\cn$ be a closed point and let $\wn_s$ be the fiber of $\wn$ over $s$.
Then the stabilizer of the $\wn_s$ (which is the subgroup of $\gn$ that fixes $\wn_s$)
is the trivial subgroup of $\gn$.
\end{coro}

We now fix the notation for the group action. Let
$G$ be any group scheme acting on an $S$-scheme $X$ via $\Psi\mh G\times X\to X$.
Let $f\mh T\to X$ be any $S$-morphism of schemes and let $\lambda\mh G\pri\to G$ be a
homomorphism of groups. Then we define the induced group action morphism
$$\Psi^{\lambda}_{f}\defeq \Psi\circ(\lambda\times f): G\pri\times T\lra X.
$$
In case $G=G\pri$ and $\lambda=1_G$,
$\Psi^{\lambda}_f$ will be shortened to $\Psi_f$;
In case $T=X$ and $f=1_X$, then $\Psi^{\lambda}_f$ will be shortened to $\Psi^{\lambda}$.
Now let $\rho\mh T\to G$ and $f\mh T\to X$ be two morphisms. We define
\begin{equation}
\label{1.103}
f^{\rho}\defeq\Psi\circ (\rho, f): T\lra X,
\end{equation}
in case the group action $\Psi$ is clear from the context.
Later, we will encounter the situation where $X$ is a scheme over a scheme $Y$ and $\rho$
is a morphism $Y\to G$. By abuse of notation we will use $\rho$ to denote the
composite morphism $X\to Y\to G$ as well and use $f^{\rho}$ to denote the morphism induced by the group action,
as in (\ref{1.103}).

Let $I\sub [n+1]$ be a subset of $l+1$ elements.
It associates an embedding $\underline{\gamma_I}\mh C[l]\to \cn$
(see (\ref{1.213})).
From the construction of $\wn$, it is clear that there is a canonical $W$-isomorphism
$$\gamma_I: W[l]\cong \wn\times_{\cn} C[l].
$$
Now let $J$ be the complement of $I$ in $[n+1]$ with entries $j_1\leq\cdots\leq j_{n-l}$.
We let $\lambda_J\mh G[n-l]\to\gn$
be the homomorphism defined so that the $j_i$-th component of $\lambda_J(\sigma)$ is $\sigma_i$
($\sigma_i$ is the $i$-th component of $\sigma$) and all other components of $\lambda_J(\sigma)$ are $1$.
Here the components of $\gn$ are defined via the isomorphism (\ref{1.125}).
Then by the explicit $\gn$-action on $\Anpo$, we see immediately that
$$\underline{\Psi}_{\uline{\gamma_I}}^{\lambda_J}: G[n-l]\times C[l]\lra \cn
$$
defined by the $\gn$-action $\underline{\Psi}$ is an open immersion.
We denote the image open subscheme by $C[\Sigma I]$ and denote $\wn\times_{\cn} C[\Sigma I]$ by
$W[\Sigma I]$. Clearly, $W[\Sigma I]$ is $G[n-l]$-invariant under the induced action
$$\Psi^{\lambda_J}: G[n-l]\times\wn\lra\wn.
$$

\begin{lemm}
\label{1.92}
The scheme $W[\Sigma I]$ (resp. $C[\Sigma I]$) is a principal fiber bundle over $W[l]$
(resp. $C[l]$) with group $G[n-l]$ via the action $\Psi^{\lambda_J}$ (resp. $\underline{\Psi}^{\lambda_J}$)
and the morphism $\gamma_I$ (resp. $\underline{\gamma}_I$) is a section of this principal
fiber bundle.
\end{lemm}

The proof is based on the parallel result concerning $\bA^{\! l+1}\to \Anpo$,
which is obvious. We shall omit the proof here.

The two principal fiber bundles are compatible in the sense that the following
diagram is commutative
$$
\begin{CD}
G[n-l]\times W[\Sigma I] @>{\Psi^{\lambda_J}}>> W[\Sigma I]\\
@VVV @VVV\\
G[n-l]\times C[\Sigma I] @>{\underline{\Psi}^{\lambda_J}}>> C[\Sigma I].
\end{CD}
$$

Often, we have several morphisms $r\mh S\to \cn$ and like to compare the corresponding
fiber products $\wn\times_{\cn}S$. Adopting the convention in groupoid, we
shall view $\wn\times_{\cn}S$ as the pull-back of $\wn\to\cn$ via $r\mh S\to\cn$.
Accordingly, we will denote such fiber product by $r\sta\wn$.

\begin{coro}
\label{1.93}
Let $(s,S)$ be a  pointed scheme, $\tilde\iota\mh S\to\cn$ be a morphism such that $\tilde\iota(s)\in C[\Sigma I]$
for some subset $I\sub [n+1]$ with $|I|=l+1$. Let $S_0={\tilde\iota}\upmo(C[\Sigma I])$, which is an open neighborhood of
$s$ in $S$. Then there is a morphism $(\rho,\iota)\mh S_0\to G[n-l]\times C[l]$ such that the restriction of
$\tilde\iota$ to $S_0\sub S$ is identical to $\underline{\Psi}^{\lambda_J}_{\underline{\gamma}_I}\circ(\rho,\iota)$.
Further it induces a canonical isomorphism
$$ {\tilde\iota}\sta \wn\times_{S}S_0\cong \iota\sta W[l],
$$
compatible to the projections to $W\times_CS_0$,
using the principal fiber bundle structure $W[\Sigma I]\to W[l]$ and the morphism $\gamma_I$.
\end{coro}

We now introduce the notion of the expanded degenerations of $W$.
Let $S$ be any $C$-scheme. An {\it effective degeneration} over $S$ is a
$C$-morphism $\xi\mh S\to \cn$ for some $n$. Note that such a morphism comes with
an $S$-family
$$\cW=\wn\times_{\cn}S\lra S
$$
and the tautological projection
$\cW\to W\times_CS$.
(For convenience, we will not distinguish the family $\cW$ from the morphism $\xi$,
and vice versa if no confusion arises. As a convention, we will call $\cW$ the associated
family of $\xi$ and call $\xi$ the associated morphism of $\cW$.)
Let $\xi_1\mh S\to C[n_1]$ and $\xi_2\mh S\to C[n_2]$ be two
effective degenerations. An {\it effective arrow} $r\mh \xi_1\to \xi_2$ consists of
a standard embedding $\iota\mh C[n_1]\to C[n_2]$
and a morphism $\rho\mh S\to G[n_2]$ so that
$(\iota\circ\xi_1)^\rho=\xi_2$. By Corollary \ref{1.93}, this identity defines
a canonical $S$-isomorphisms $\xi_1\sta W[n_1]\cong \xi_2\sta W[n_2]$ compatible with their projections
to $W\times_CS$.
We call this the associated isomorphism of the arrow $r\mh\xi_1\to\xi_2$.
By abuse of notation, we will denote this isomorphism by $r$ as well.
Now let $\xi_1$ and $\xi_2$ be any two effective degenerations over $S$. We say $\xi_1$ is equivariant
to $\xi_2$ via an effective arrow if there is either an effective arrow $\xi_1\to \xi_2$ or an effective arrow
$\xi_2\to\xi_1$. We say $\xi_1$ and $\xi_2$ are equivalent via a sequence of effective arrows if there is a
sequence of effective degenerations $\eta_0,\cdots,\eta_m$ so that $\eta_0=\xi_1$, $\eta_m=\xi_2$
and $\eta_i$ is equivariant to $\eta_{i+1}$ via an effective arrow for all $i$.
Note that once such a sequence of effective arrows is given, then there is a unique induced
$S$-isomorphism $\xi_1\sta W[n_1]\cong \xi_2\sta W[n_2]$ compatible with their
projections to $W\times_CS$.

Now let $\cW_1$ and $\cW_2$ be two effective degenerations over $S$.
We say $\cW_1$ is {\sl isomorphic} to $\cW_2$ if $\cW_1$ is $S$-isomorphic to $\cW_2$
and the isomorphism is compatible to their tautological projections to $W\times_C S$.
The following Lemma says that such isomorphisms
are locally generated by a sequence of arrows.

\begin{lemm}
\label{1.167}
Let $\xi_1$ and $\xi_2$ be two effective degenerations over $S$ so that their
associated families $\cW_1$ and $\cW_2$
are isomorphic. Then to each $p\in S$ there is an open neighborhood $S_0$ of $p\in S$
so that the (induced) isomorphism between
$\cW_1\times_SS_0$ and $\cW_2\times_SS_0$ is induced by a sequence of effective
arrows between $\xi_1\times_SS_0$ and $\xi_2\times_SS_0$.
\end{lemm}

\begin{proof}
In case $p$ lies over $C-0$, then there is nothing to prove. Now assume $p$ lies over $0\in C$.
Let $n$ be the integer so that $\cW_1\times_Sp$ has $n+2$ irreducible components. Then
by Corollary \ref{1.93}, there is an open neighborhood $S_0$ of $p\in S$ and two
morphisms $\xi_1\pri,\xi_2\pri \mh S_0\to \cn$ so that $\xi_i\times_S S_0$ is equivalent to
$\xi_i\pri$ via an arrow. Thus to prove the Lemma it suffices to consider the case where
$n_1=n_2=n$. Now we assume, $\xi_i\mh S\to\cn$ and $\xi_i(p)=\bz\in\cn$.
Let $\varphi:\cW_1\lra\cW_2$ be the $W_S$-isomorphism and let
$\phi_i:\cW_i\lra\wn$ be the associated morphisms.
We define a functor $\Isom[S]$ from the category of $S$-schemes to sets
that associates to any $S$-scheme $T$ the set of all morphisms $\rho\mh T\to\gn$
such that with $\iota\mh \cW_1\times_ST\to\cW_1$ the induced morphism (i.e. the first projection) then
$$(\phi_1\circ\iota)^{\rho}=\phi_2\circ\varphi\circ\iota
:\cW_1\times_S T\lra \wn.
$$
(For the notation here see the formula (\ref{1.103}).)
By the existence theorem of Hilbert schemes
developed by Grothendieck coupled with the trivialization of the stabilizers
(see Corollary \ref{1.9}), $\Isom[S]$ is represented by a subscheme $R$ of
$S\times\gn$.

We now show that $R$ is non-empty and dominates a neighborhood of $p\in S$.
Since $\xi_1(p)=\xi_2(p)=0\in\cn$,
the restriction of the isomorphism $\varphi$ to fibers over $p$
induces an isomorphism
$\wn_0\cong\wn_0$ that commutes with $\wn_0\to W_0$.
Then it is induced by an element $g_0\in \gn$ by
Corollary \ref{1.8}. Hence $(p,g_0)\in R$.
Let $\mm$ be the maximal ideal of the local ring $\cO_{S,p}$.
By the previous discussion, we have $\rho_0\mh p\to\gn$ so that
$(i_0,\rho_0)\mh p\to S\times\gn$ factor through $R\sub S\times\gn$,
where $i_0\mh p\to S$ is the inclusion.
Let $B_k=\spec \cO_{S,p}/\mm^{k}$ with $i_k\mh B_k\to S$ the inclusion morphism.
Suppose for some $k\geq 1$ we have
$$(i_{k-1}, \rho_{k-1}): B_{k-1} \lra S\times\gn,
$$
extending $(i_0,\rho_0)$ that factor through $R\sub S\times \gn$. We will
show that $\rho_{k-1}$ extends to $\rho_k$ satisfying similar
condition. Since $\gn$ is smooth, we can extend $\rho_{k-1}$ to
$h_k\mh B_k\to \gn$. Let $\iota_k\mh\cW_1\times_S B_k\to\cW_1$ be the canonical immersion.
We consider the problem of extending
the $W$-morphism
$$(\phi_1\circ\iota_{k-1})^{h_{k-1}}=\phi_2\circ \varphi\circ\iota_{k-1}: \cW_1\times_S B_{k-1}\lra\wn
$$
to a $W$-morphism
$$\cW_1\times_S B_k\lra \wn.
$$
This is a typical extension problem.
Since $(\phi_1\circ\iota_{k})^{h_{k}}$ and $\phi_2\circ\varphi\circ\iota_{k}$ are two such
extensions, by deformation theory they are related by an element in
$\mm^k/\mm^{k+1}\otimes K$, where
$$K=\ker\left\{\Hom_{\wn_0}(\Omega_{\wn},\cO_{\wn_0})
\lra \Hom_{\wn_0}(\Omega_{W}\otimes_{\cO_W}\cO_{\wn_0},\cO_{\wn_0})\right\}.
$$
We claim that $K$ is isomorphic to $\kk^{\oplus n}$ and
is generated by the group action of $\gn$.
We let $i\mh \wn_0\to \wn$ be the immersion and $j\mh \wn_0\to\cn$ be the composite of $i$ with the
projection. Then we have the exact sequence
\begin{equation}
\label{1.99}
0\lra j\sta\Omega_{\cn}\lra i\sta\Omega_{\wn}\lra \Omega_{\wn_0}\lra 0
\end{equation}
and its induced exact sequence
$$
\Hom(i\sta\Omega_{\wn},\cO_{\wn_0})
\mapright{\beta}\Hom(j\sta\Omega_{\cn},\cO_{\wn_0})
\mapright{\delta}\Ext^1(\Omega_{\wn_0},\cO_{\wn_0}) .
$$
At the origin $\bz\in\cn$, we have a canonical isomorphism $T_\bz\cn\cong T_0\Anpo$. We let
$v_l\mh \kk\to\kk^{\oplus n+1}\equiv T_0\Anpo$ be the $l$-th component. The dual
of $v_l$, which is $v_l\dual\mh\Omega_{\cn}\otimes_{\cO_{\cn}}\kk_0\to\kk$,
defines a homomorphism $\bv_l\dual\mh j\sta\Omega_{\cn}\to\cO_{\wn_0}$.
Now let
$$0\ne h\in \image(\beta)\sub \Hom(j\sta\Omega_{\cn},\cO_{\wn_0}),
$$
say given by a homomorphism $j\sta\Omega_{\cn}\to \cO_{\wn_0}$.
Then it must be a linear combination of $\bv_l\dual$, say equal to
$a_1\bv_1\dual+\cdots+a\lnpo\bv\lnpo\dual$.
Now let $c$ be the extension class in $\Ext^1(\Omega_{\wn_0},j\sta\Omega_{\cn})$
of the exact sequence (\ref{1.99}) and let
$$h\sta: \Ext^1(\Omega_{\wn_0},j\sta\Omega_{\cn})\lra
\Ext^1(\Omega_{\wn_0},\cO_{\wn_0})
$$
be the homomorphism
induced by $h$. We claim that $h\sta(c)\ne 0$. Indeed, let $T_2=\spec \kk[t]/(t^2)\to\cn$
be the immersion representing the tangent vector
$\bv=a_1\partial_{z_1}+\cdots+a\lnpo\partial_{z\lnpo}\in T_0\cn$. Then
$$
\wn\times_{\cn}\spec \kk[t]/(t^2)\lra \spec \kk[t]/(t^2)
$$
is an infinitesimal smoothing of the $l$-th nodal divisor of ${\wn_0}$ whenever
$a_l\ne 0$, and consequently the exact sequence
\begin{equation}
\label{1.100}
0\lra \cO_{\wn_0}\lra\Omega_{\wn\times_{\cn}T_2}\otimes_{\wn\times_{\cn}T_2}\cO_{\wn_0}
\lra \Omega_{\wn_0}\lra 0
\end{equation}
does not split since $h\ne 0$. On the other hand, this exact sequence is exactly defined by the extension
class $\delta(h)$.
In case $h\in\image(\beta)$, then $\delta(h)=0$, violating the non-splitness of
(\ref{1.100}). This proves $\beta=0$.
As a consequence
$$\Hom(\Omega_{\wn_0},\cO_{\wn_0})\cong \Hom(i\sta \Omega_{\wn},\cO_{\wn_0}).
$$

Now let $\alpha$ be any element in $K$. Then $\alpha\in\Hom(\Omega_{\wn_0},\cO_{\wn_0})$
by the above identity. Hence the restriction of $\alpha$ to $\Delta_l$,
denoted by $\alpha_l$,
is a vector field of $\Delta_l$. Since $\alpha$ lies in the kernel, we have $\alpha_1=\alpha_{n+2}=0$.
Furthermore, the restriction of $\alpha_l$ to the two (or one if $l=1$ or $n+2$)
distinguished divisors is tangential to these divisors and must satisfy the compatibility
condition $\alpha_l|_{D_-}=\alpha_{l+1}|_{D_+}$.
Now let $V_l$ be the subspace of $H^0(T_{\Delta_l})$ consisting of the vector fields whose
restrictions to the two distinguished divisors are tangential to the two divisors.
Then $V_l$ is canonically isomorphic to the direct sum $H^0(T_D)\oplus \CC$, where
$\CC$ is generated by the vector field generated by the group action $\gn_l$.
Let $e(\alpha_l)$ be the component of $\alpha_l$ in $H^0(T_D)$.
Then the compatibility condition translates to
$e(\alpha_1)=\cdots=e(\alpha_{n+2})$. Adding the fact that $\alpha_1=0$, we conclude that $e(\alpha_l)=0$
for $2\leq l\leq n+1$. Therefore, $\alpha_l$ must be a vector field tangential to the fibers of
$\Delta_l\to D$ and vanishes on the two distinguished divisors. By Corollary \ref{1.8},
such $\alpha$ is generated by the group action $\gn$.

If we apply this to the two extensions $(\phi_1\circ\iota_k)^{h_k}$ and $\phi_2\circ\varphi\circ\iota_k$,
we see that there is a morphism
$h^+_k\mh B_k\to\gn$ with $h_k^+|_{B_{k-1}}\mh B_{k-1}\to\gn$ factor through
$\{e\}\sub\gn$ such that with $\rho_k=h_k^+\cdot \iota_k$,
$$(\phi_1\circ\iota_k)^{{\rho}_k}=\phi_2\circ\varphi\circ\iota_k:
\cW_1\times_S B_k\lra\wn.
$$
Since $k$ is arbitrary, this implies that the projection $R\to S$ dominates
a neighborhood of $p\in S$.

We let $\pi\mh R\to S$ be the projection.
To complete the proof of the Lemma, we need to find a neighborhood
$R_0$ of $(p,\rho_0)\in R$, so that $\pi|_{R_0}\mh R_0\to S$ is an open neighborhood of $p\in S$.
First, by Corollary \ref{1.9}, for any $s\in S$ the set $\pi\upmo(s)$ has
at most one point. Applying Corollary \ref{1.9} again, we see that
the choice of $\rho_k$ in the proof is unique.
Hence the formal completion of $R$ along
$(p,g_0)$ is isomorphic to the formal completion of $S$ along $p$, under $\pi$.
Hence $R\to S$ is \'etale over a neighborhood of $p\in S$. Combining the one-to-one
and \'etale property, we conclude that there is a neighborhood $R_0\sub R$ of
$(p,g_0)$ so that $\pi|_{R_0}\mh R_0\to S$ is an open immersion.
Let $\rho\mh R_0\to\gn$ be the morphism given by the definition of the functor
$\Isom[S]$ and let $\eta\mh R_0\to S$ be the inclusion. Then $\rho$
defines an effective arrow between $\xi_1\times_S R_0$ and $\xi_2\times_S R_0$.
This proves the Lemma.
\end{proof}

We now define the stack $\WWc$ of the expanded degenerations of $W$.

\begin{defi}
Let $S$ be any $C$-scheme. An expanded degeneration of $W$
over $S$ is a pair $(\cW,\rho)$, where $\cW$ is a family over $S$ and $\rho\mh\cW\to W\times_CS$
is an $S$-projection, so that there is an
open covering $S\lalp$ of $S$ such that for each $S\lalp$ the restriction pair
$(\cW\times_S S\lalp,\rho|_{\cW\times_S S\lalp})$
is isomorphic to an effective expanded degeneration of $W$.
Let $(\cW,\rho)$ and $(\cW\pri,\rho\pri)$ be two degenerations over $S$ and $S\pri$ respectively.
An {\it arrow} $\cW\to\cW\pri$ consists of a $C$-morphism $S\to S\pri$ and an $S$-isomorphism
$\cW\to\cW\pri\times_{S\pri}S$ compatible to their projections to $W\times_CS$.
\end{defi}

We define the groupoid of the expanded degenerations of $W$ to be the category $\WWc$ whose objects
are expanded degenerations of $W$ over $S$, where $S$ is a $C$-scheme.
In the following, we will call expanded degenerations of $W$ simply degenerations
if its meaning is clear from the context. The functor $\mathfrak p\mh \WWc
\to({\rm Sch}/C)$ is the one that sends degenerations over $S$ to $S$. Arrows between two
objects are those defined in the above Definition. Obviously, if $\xi\in\WWc$
with ${\mathfrak p}(\xi)=S$ and $i\mh T\to S$ is a morphism, then there is a unique pull-back family
$i\sta\xi\in\WWc$ with an arrow $i\sta\xi\to\xi$. In this way, $(\WWc,\mathfrak p)$ forms a groupoid.
For any $S\in ({\rm Sch}/C)$, we define $\WWc(S)$ to be the category of all degenerations over $S$.

\begin{prop}
The groupoid $\WWc$ is a stack.
\end{prop}

\begin{proof}
This is straightforward and will be left to the readers.
\end{proof}

\section{Pre-deformable morphisms}

Recall that a pre-stable curve is a connected complete
curve with at most nodal singularity.
(Since the marked points will not affect our discussion of pre-deformable morphisms,
we will discuss curves without marked points in this section.)
A morphism $f$ whose domain is a pre-stable curve is called a pre-stable
morphism.
Let $f\mh X\to\wn$ be a pre-stable morphism so that $f(X)\sub\wn_0$.
We say $f$ is non-degenerate
if no irreducible component of $X$ is mapped entirely to the
singular locus of $\wn_0$.
It is well-known that in general $f$ may not be deformable to a pre-stable
morphism to $\wn$ over a general point $t\in\cn$. One set of conditions is purely local
in nature. It reflects the fact that $\wn_0$ has nodal singularity. In this section,
we will investigate such condition in details.

\subsection{Morphisms of pure contact}

We will investigate the following situation in this subsection.
Let $\kt\to\kw$ be the homomorphism defined by $t\mapsto w_1w_2$ and let $\phi\mh\ks\to\kz$
be the homomorphism defined by $\phi(s)=z_1z_2$ or $\phi(s)=z_1$.
Note that the former corresponds to a smoothing of the node and the later corresponds to a
family of smooth curves. We let $A$ be a $\ks$-algebra and let $R=\kz\otimes_{\ks}A$.
We assume $A$ is $(s)$-adically complete and let $\hat R$ be the $(z_1,z_2)$-adic
completion of $R$.
In this part, we will investigate the situation where there is a homomorphism
$\psi\mh \kt\to A$ with induced homomorphism $\tilde\psi\mh\kt\to \hatR$ via
$\tilde\psi(t)=1\otimes\psi(t)$ and a $\kt$-homomorphism
\begin{equation}
\label{12.1}
\varphi: \kw\lra \hatR
\end{equation}
satisfying the following non-degeneracy condition:
{\sl The ideals $(\varphi(w_1),\varphi(w_2))$
and $(z_1,z_2)$ in $\hatR$ satisfy}
$$(z_1,z_2)^m\sub (\varphi(w_1),\varphi(w_2))\sub (z_1,z_2)
$$
for some integer $m$.

Geometrically, this states that there are
no irreducible components of the fibers of $\spec\hatR\to \spec A$ that
are mapped entirely to the node in $\spec\kw/(w_1w_2)$.

We now state and prove our first technical result which explains the notion of curves of pure contact.
We begin with the notion of normal forms of elements in $\hatR$.

\begin{lemm}
\label{2.5}
Let $\hatR$ be as before. Then any $\alpha\in \hatR$
has a unique expression, called its normal form, as
\begin{equation}
\label{12.3}
\alpha=a_0+\sum_{i>0}a_i z_1^i+\sum_{i>0} b_iz_2^i,\quad a_i,b_i\in A.
\end{equation}
\end{lemm}

We call such expansion the normal form of $\alpha$.

\begin{proof}
Clearly, any $\alpha$ can be expressed as above. We now show that
such expression is unique.
Let $I_n\sub \hatR$ be the $A$-sub-module spanned by all $z_1^iz_2^j$ satisfying $|i-j|>n$.
Clearly, the quotient $A$-module
$\hatR/I_n$ is a free $A$-module with basis
$z_1^n,\cdots,z_1,1,z_2,\cdots,z_2^n$.
Let $\pi_n: \hatR\lra \hatR/I_n$
be the quotient homomorphism of $A$-modules.
Assume $\alpha=0$ has an expression as in (\ref{12.3}), then
$0=\pi_n(\alpha)=a_0+\sum_{j\geq 1}^n a_j z_1^j+b_j z_2^j$ in $\hatR/I_n$.
Since $\hatR/I_n$ is a free $A$-module and $z_1^n\cdots,z_2^n$ is a basis,
all $a_j$ and $b_j$ are zero for $j\leq n$. Since $n$ is arbitrary, this implies that
all $a_j$ and $b_j$ are zero. This proves the Lemma.
\end{proof}

\begin{prop}
\label{2.4}
Let the notation be as before.
We assume further that $A$ is a local ring and is flat over $\kt$.
Then of the two choices of $\phi\mh\ks\to\kz$ given before,
only $\phi(s)=z_1z_2$ is possible. Further, possibly after exchanging $z_1$ and $z_2$,
the homomorphisms $\varphi$  and $\psi$ must be of the
forms
$$\varphi(w_1)=z_1^n\beta_1,\quad \varphi(w_2)=z_2^n\beta_2\quad{\rm and}
\quad \psi(t)=s^n\epsilon
$$
for some integer $n$, units $\beta_1$
and $\beta_2\in \hatR$ and a unit $\epsilon\in A$ such that
$\beta_1\beta_2=\epsilon$ in $\hatR$.
\end{prop}

\begin{proof}
Let $\mm\sub A$ be the maximal idea and let
\begin{equation}
\label{2.6}
\varphi(w_i)=\sum_{j\geq 0}a_{i,j}z_1^j + \sum_{j\geq 0} b_{i,j} z_2^j,\qquad a_{i,0}=b_{i,0}
\end{equation}
be the normal form of $\varphi(w_i)$.
Since $\varphi$ is non-degenerate, the homomorphism
$$\varphi_0:\kw\lra \kz/(z_1z_2)
$$
induced by $\varphi$ and $\kz/(z_1,z_2)\cong \hatR/\mm\hatR$
must be of the form (possibly after exchanging $z_1$ and $z_2$)
$$
\varphi_0(w_1)= c_1 z_1^{n_1}\!\!\mod z_1^{n_1+1}\quad {\rm and}\quad
\varphi_0(w_2)= c_2 z_2^{n_2}\!\!\mod z_2^{n_2+1},\quad c_1,c_2\ne 0,
$$
for some positive integers $n_1$ and $n_2$.
Hence in (\ref{2.6}), all $a_{1,j}$ for $j<n_1$ and
$b_{2,j}$ for $j<n_2$ are in the maximal ideal $\mm\sub A$.
We let
$$I=\bl\{a_{1,j}s^j,b_{1,j},a_{2,j},b_{2,j}s^j|j\geq 0\}\br\sub A
$$
be the ideal generated by the listed terms for all $j\geq 0$.

We first show that $I\sub \mm$.
After expressing $\varphi(w_1)\varphi(w_2)$ in its normal form
$\varphi(w_1)\varphi(w_2)=c_0+\sum_{j>0}c_jz_1^j+d_jz_2^j$ we have
$$
c_k=\sum_{i+j=k}a_{1,i}a_{2,j}+\sum_{i\geq0} a_{1,k+i}b_{2,i}s^i
+\sum_{i\geq 0}a_{2,k+i}b_{1,i} s^i,\quad k>0
$$
and
$$
d_k=\sum_{i+j=k}b_{1,i}b_{2,j}+\sum_{i\geq0} b_{1,k+i}a_{2,i}s^i
+\sum_{i\geq 0}b_{2,k+i}a_{1,i} s^i, \quad k>0.
$$
Because $\varphi(w_1)\varphi(w_2)=\varphi(w_1w_2)$ is the image of $t$
in $A$, all $c_k$ and $d_k$ for $k>0$ are zero.
Now let $l$ be the largest integer so that $b_{1,j}\in\mm$ for $j<l$.
Then by considering the relation $d_{l+n_2}=0$, we obtain
$$\sum_{i+j=l+n_2} b_{1,i}b_{2,j}+b_{1,n_2+l}a_{2,0}+b_{2,n_2+l}a_{1,0}\in\mm.
$$
By assumption, all terms in the above expression except $b_{1,l}b_{2,n_2}$ are in
$\mm$. Hence $b_{i,l}b_{2,n_2}\in\mm$. Since $b_{2,n_2}\not\in\mm$, $b_{1,l}\in\mm$.
This shows that there is no such $l$ and hence all $b_{1,j}$ are in $\mm$.
Similarly, by considering the relation $c_k=0$ we
can show that all $a_{2,j}\in\mm$.

Now let $n=\min(n_1,n_2)$ and let $\bar I$ be the associated ideal of $I$ in the
quotient ring $A/(s^n)$. Let $\bar \mm$ be the maximal ideal of $A/(s^n)$. We
now show that $\bar I\sub\bar \mm\cdot\bar I$. Note that
in the expression of $c_k$ and $d_k$, terms $a_{2,k+i}b_{1,i}s^i$ and
$b_{1,k+i}a_{2,i}s^i$ are in $\mm\cdot I$ for all $k,i\geq 0$.
We let $J=\mm\cdot I+(s^n)$. Note that $b_{2,k}s^k\in J$ for $k\geq n$.
Now let $k$ be the smallest integer so that $b_{2,k+i}s^{k+i}\in J$ for
all $i>0$. Hence $k<n$. Using $c_{n_1-k}=0$, we obtain
$$\sum_{n_1-k+i\geq 0} a_{1,n_1-k+i}b_{2,i}s^i+
\sum_{i+j=n_1-k}a_{1,i}a_{2,j} \, \in\,\mm\cdot I,
$$
which implies
$$a_{1,n_1}b_{2,k}s^k+\sum_{k-n_1<i<k}a_{1,n_1-k+i}b_{2,i}s^i
+\sum_{k<i}a_{1,n_1-k+i}b_{2,i}s^i
\equiv a_{1,n_1}b_{2,k}t^k \! \mod J.
$$
Hence $b_{2,k}s^k\in J$. By minimality of $k$ we have
$b_{2,j}s^j\in J$ for all $j\geq 0$.
For similar reason, we can prove $a_{1,k}s^k\in J$ for all $k\geq 0$.
To show $a_{2,k}\in J$, we use the relation $c_{n_1+k}=0$, which
combined with $a_{2,k+i}b_{1,i}\in J$ and
$b_{2,j}s^j\in J$ implies
$\sum_{i+j=n_1+k} a_{1,i}a_{2,j}\in J$.
Now assume $a_{2,j}\in J$ for $j<k$, then the above inclusion
implies
$$a_{1,n_1}a_{2,k}+\sum_{i>0}a_{1,n_1-i}a_{2,k+i}+\sum_{i>0}a_{1,n_1+i}a_{2,k-i}
\in J.
$$
Since $a_{1,n_1-i}\in\mm$ and $a_{2,k+i}\in I$, we have $a_{1,n_1}a_{2,k}\in J$
and hence $a_{2,k}\in J$. For similar reason, $b_{1,k}\in J$
for all $k>0$. Combined, this proves that $\bar I\sub\bar\mm\cdot\bar I$.
Since $A$ is Noetherian, this is possible only if $\bar I=0$. This
proves $I\sub (s^n)$.

We next show $n_1=n_2$. Assume $n=n_1<n_2$. Then $d_{n_2-n_1}=0$ implies that
$$
\sum_{i+j=n_2-n_1}b_{1,i}b_{2,j}+\sum_{i\geq 0} b_{1,n_2-n_1+i}a_{2,i}s^i
+\sum_{i\geq 0} b_{2,n_2-n_1+i}a_{1,i}s^i=0.
$$
Note that all terms in this expression except $b_{2,n_2}a_{1,n_1}s^{n_1}$ belong to $\mm\cdot(s^n)$. This
is impossible. Hence $n_1=n_2$.

We now look at the normal form of $\varphi(w_1)$.
Since $b_{1,j}\in (s^n)$ for all $j\geq 0$, there are
$b_{1,j}\pri\in A$ so that $b_{1,j}=b_{1,j}\pri s^n$.
On the other hand, since $a_{1,j}s^j\in (s^n)$,
there are $a_{1,j}\pri$ so that $a_{1,j}=a_{1,j}\pri s^{n-j}$ for $j<n$.
Using $s=z_1z_2$, we see that there are $\beta_1$ and $\beta_2$ in $\hatR$ so that
$\varphi(w_i)=z_i^n \beta_i$. Further by our choice of $n$ both $\beta_1$ and $\beta_2$
are units in $\hatR$.
It remains to show that there is a unit $\epsilon\in A$ so that $\psi(t)=s^n\epsilon$ and
$\beta_1\beta_2=\epsilon$. First, since $\psi(t)=\varphi(w_1w_2)$, which is in the ideal $(s^n)$,
there is an $\epsilon\in A$ so that $\psi(t)=s^n\epsilon$. Hence $s^n(\beta_1\beta_2-\epsilon)=0$
in $\hatR$. Because $\hatR$ is flat over $\kt$, this implies that $\beta_1\beta_2-\epsilon=0$.
Hence $\epsilon$ is a unit in $A$.
This completes the proof of the Lemma.
\end{proof}

This Lemma inspires the notion of homomorphisms of
pure contacts as follows.

\begin{defi}
\label{2.10}
Let the notation be as stated in the beginning of this subsection with $A$ only assumed to be $(s)$-adically
complete. We say $\varphi$ has pure contact if,
possibly after exchanging $z_1$ and $z_2$, there is a unit $\beta$
in $\hatR$, a unit $\epsilon$ in $A$ and an integer $n$ such that $\varphi(w_1)=\beta z_1^n$
and $\varphi(w_2)=\epsilon\beta\upmo z_2^n$.
\end{defi}

The integer $n$ is called the order of contact. Sometimes to emphasize this order we
will say $\varphi$ is of pure $n$-contact.
The following is the main existence Lemma of this section.

\begin{lemm}
\label{2.11}
Let $n$ be any positive integer. Suppose $a_{1,n}$ and $b_{2,n}$ in the
normal forms of $\varphi(w_1)$ and $\varphi(w_2)$ in (\ref{2.10}) respectively are units in $A$.
Then there is an ideal $I_A\sub A$ so that for any $A\to T$ the induced homomorphism
$$\varphi_T:\kw\lra \hatR\otimes_A T
$$
has pure $n$-contact if and only if $A\to T$ factor through $A/I_A\to T$.
The ideal $I_A$ satisfies the following base change property: Let $A\to A\pri$ be a
homomorphism of rings and let $\varphi\pri\mh\kw\to(\kz\otimes_{\ks}A\pri)\rhat$,
where $(\cdot)\rhat$ is the $(z_1,z_2)$-adic completion of the respective ring,
be the induced $\kt$-homomorphism.
Then $I_{A\pri}=I_A\cdot A\pri$.
\end{lemm}

\begin{proof}
We let $\zeta,\epsilon, \xi_j,\eta_j$, where $j$ runs through all positive integers, be indeterminants.
We consider the ring
$B=A[\![\zeta,\zeta\upmo,\epsilon,\epsilon\upmo,\bxi,\bleta]\!]$,
where $\bxi$ and $\bleta$ mean $(\xi_1,\xi_2\cdots)$ and $(\eta_1,\eta_2,\cdots)$ respectively.
We let $C$ be the $(z_1,z_2)$-adic completion of $\kz\otimes_{\ks}B$.
We let $\gamma\in C$ and its inverse be
$$\gamma=\zeta(1+\sum_{j>0}\xi_i z_1^i+\eta_j z_2^j)
\quad{\rm and}\quad
\gamma\upmo=\zeta\upmo(1+\sum_{j>0}\tilde\xi_i z_1^i+\tilde\eta_j z_2^j),
$$
where $\tilde\xi_j$ and $\tilde\eta_j$ are elements in $B$.
We consider the ideal $J\sub C$ that is generated by all
coefficients of $z_i^j$ in the normal forms of
\begin{equation}
\label{2.60}
\Phi_1=\varphi(w_1)-z_1^n\zeta(1+\sum_{j>0} \xi_j z_1^j+\sum_{j>0} \eta_j z_2^j)
\end{equation}
and
\begin{equation}
\label{2.61}
\Phi_2=\varphi(w_2)-z_2^n\zeta\upmo(1+\sum_{j>0} \tilde\xi_j z_1^j+\sum_{j>0} \tilde\eta_jz_2^j)\epsilon.
\end{equation}
We now investigate these generators of $J$, or equivalently the relations in $C/J$.
First note that the coefficient of $z_1^{n+k}$ in $\Phi_1$ is $a_{1,k+n}-\xi_k$. Hence
we obtain elements $a_{1,k+n}-\xi_k$ in $J$ for all $k\geq 1$.
Similarly, the coefficient of $z_1^n$ in $\Phi_1$ is $a_{1,n}-\zeta$, which by definition is
in $J$. In the following we will call them the canonical relations of $\xi_k$ and $\eta$,
respectively.
We let $J_1$ be the ideal generated by $a_{1,n}-\zeta$ and all $a_{1,k+n}-\zeta_k$.
We next find relations of $\eta_j$.
First, by expanding $1/\gamma$ in power series and using $s=z_1z_2$, we can express
$\tilde\eta_k$ canonically as
$\tilde\eta_k=-(\eta_k+\sum_{j\geq 1}c_{k,j}s^j)$, where
$c_{k,j}\in \kk[\![\bxi,\bleta]\!]$.
Then the coefficient of $z_2^{k+n}$ in $\Phi_2$ become
$$b_{2,n+k}+\zeta\upmo\eta_k\epsilon+\sum_{j\geq 1}\zeta\upmo c_{k,j}s^j\epsilon.
$$
We let $J_2$ be the ideal generated by $J_1$ and this set of elements for all $k>0$.
Since $c_{k,j}\in \kk[\![\bxi,\bleta]\!]$, modulo $J_2$
we can substitute $\zeta$ and all $\xi_l$ by elements in $A$ using relations in $J_1$
and then substitute all $\eta_l$ in $c_{k,j}$ by $-\zeta\epsilon\upmo b_{2,n+k}-\sum c_{k,j}s^j$.
This way we obtain
$$\eta_k-\Bigl( d_{1,k}+\sum_{j\geq 2}\tilde c_{k,j}s^j\Bigr)\in J_2,
$$
where $d_{k,1}\in A$ and $\tilde c_{k,1}\in A[\![\bleta,\epsilon\upmo]\!]$. Note that
here we have used the fact that $A$ is complete with respect to the ideal $(s)$. Repeating this procedure,
we eventually obtain an element $\eta_k-f_k(\epsilon\upmo)\in J_2$, where $f_k(\epsilon\upmo)\in A[\![\epsilon\upmo]\!]$.
It remains to find relation in $\epsilon$. By using the coefficient of $z_2^n$ in $\Phi_2$, we
obtain an element
$$b_{2,n}-\epsilon\zeta\upmo(1+\sum_{j\geq 1}h_js^j)\in J,
$$
where $h_j\in \kk[\![\bxi,\bleta]\!]$. We then let $J_0\sub C$ be the ideal generated by
$J_2$ and this element. Then again modulo $J_0$, which amounts to
repeatedly replace $\zeta$, $\xi_k$ and $\eta_k$
by relations in $J_0$, we obtain
an element $-1+\zeta\epsilon\upmo b_{2,n}+s\phi(\epsilon\upmo)\in J_0$,
where $\phi(\epsilon\upmo)\in A[\![\epsilon\upmo]\!]$.
By iteration and taking the limit, we obtain a canonical element $\epsilon-h\in J$ for
some $h\in A$. It is clear from this construction that $h$ is a unit since $b_{2,n}$ is a unit.
This way we have found elements $c_{1,k}$ and $c_{2,k}\in A$ for $k\geq 1$ and a unit $h\in A$ so that
\begin{equation}
\label{12.10}
J_0=\Bigl(\{\xi_k-c_{1,k},\ \eta_k-c_{2,k}\mid k>0\};\ \zeta-a_{1,n},\ \epsilon-h\Bigr)\sub C.
\end{equation}
From this set of generators, we see that the canonical homomorphism
$A\to C/J_0$ is an isomorphism. Since $J_0$ is a sub-ideal of $J$, we have quotient homomorphism $C/J_0\to C/J$.
We let $I_A\sub A$ be the kernel of $A\to C/J_0\to C/J$. We will show that
the ideal $I_A$ satisfies the property specified in the statement of the Lemma.
First, from the construction it is direct to check that $I_A$ satisfies the base change property.
Now let $A\to T$ and $\varphi_T$ be as in the Lemma so that $\varphi_T$ has pure $n$-contact.
Then it follows from our construction that the ideal
$I_{T}=\{0\}$.
Therefore $A\to T$ factor through $A/I_A\to T$ since $I_T=I_A\cdot T$.
This proves one direction of the Lemma.
Now assume $A\to T$ is a homomorphism so that it factor through $A/I_A\to T$.
Then by the base change property we have $I_T=\{0\}$.
It follows from the construction of $I_T$ that $\varphi_T$ does have
pure $n$-contact. This completes the proof of the Lemma.
\end{proof}

\subsection{Families of pre-deformable morphisms}

In this subsection, we will introduce the notion of pre-deformable morphisms. We will
then prove that for any family of non-degenerate stable morphisms $f\mh\cX\to\wn$
over $S$, the locus $S_{pd}\sub S$ of all pre-deformable morphisms admits a natural
closed subscheme structure.

We begin with some notation that will be used throughout this section.
We let $l$ be any integer in $[n+1]$.
Following our convention in section 1, we let $l\dual$ be the complement of $l$ in $[n+1]$ and
$\cn_{[l\dual]}$ the smooth divisor $C\times_{\Ao} \bH_l\unpo\sub C[n]$,
where $\bH_l\unpo\sub \Anpo$ is the $l$-th coordinate hyperplane.
We define
$$C_0[n]\defeq C[n]\times_C 0=\cup_{l=1}^{n+1} \cn_{[l\dual]}.
$$
We denote by $\bD_l$ the nodal divisor of $\wn\times_{\cn}\cn_{[l\dual]}$
and by $\Delta_l^-$ and $\Delta_l^+$ its two smooth
irreducible components.
According to our ordering of the irreducible
components of $\wn_0$, $\Delta_l^-$ is on the left of $\Delta_l^+$. Then $\wn\times_{\cn}\cn_{[l\dual]}$ is the
union of $\Delta_l^-$ and $\Delta_l^+$ intersecting transversely along
$\bD_l$.
Let $\pi_l\mh\cn\to \Ao$ be induced by the $l$-th projection $\Anpo\to\Ao$.
For any point $q\in \bD_l$ there
is a neighborhood $\cW$ of $q\in\wn$ with two regular
functions $w_1$ and $w_2$ so that $w_1=0$ (resp. $w_2=0$) defines the divisor $\cW\cap \Delta_l^+$
(resp. $\cW\cap \Delta_l^-$). By adjusting $w_2$
if necessary, we
can assume $w_1w_2$ is the image of $t_l$ via $\Gamma(\Anpo)\to\Gamma(\cW)$.
Hence by shrinking $\cW$ if necessary, we obtain a smooth morphism
\begin{equation}
\label{4.109}
\psi:\cW\lra\Theta_l\defeq \spec\kw\otimes_{\kk[t_l]}\ktt.
\end{equation}
Here $\kk[t_l]\to\kw$ is defined by $t_l\mapsto w_1w_2$ and $\kk[\bt]=\kk[t_1,\cdots,t\lnpo]$.
Later we will call $\cW$ with $\psi$ implicitly understood an admissible chart of $\wn$ at $q\in\bD_l$.

Let $f\mh X\to\wn$ be a pre-stable morphisms.
We say $f$ is non-degenerate if no irreducible component of $X$ is mapped entirely to any of $\bD_1,
\cdots,\bD\lnpo$. Now assume $f$ is non-degenerate and
let $p\in f\upmo(\bD_l)$. We define the notion of
upper and lower contact order of $f$ along $p$.
Let $\hat X$ be the formal completion of $X$ along $p$ and let $\hat \Delta_l^-$ (resp.
$\hat\Delta_{l}^+$; resp. $\hat\bD_l$) be the formal completion of $\Delta_l^-$
(resp. $\Delta_l^+$; resp. $\bD_l$)
along $f(p)$. For any irreducible component $\hat V$ of $\hat X$, we can define
the contact order of $f$ along $\hat V$ as follows. Assume $f(\hat V)\sub
\hat \Delta_l^-$. We let $z$ be a generator of the maximal ideal of $\cO_{\hat V}$
and let $w$ be a generator of the ideal
$\cI_{\hat \bD_l\sub\hat\Delta_l^-}\sub \cO_{\hat \Delta_l^-}$.  Let $f\sta\mh\cO_{\hat\Delta_l^-}\to
\cO_{\hat V}$ be the induced homomorphism. Then $f\sta(w)=z^{\delta}\beta$ for a
positive integer $\delta$ and a unit $\beta\in\cO_{\hat V}$.
The integer $\delta$ is the contact order of $f$ along $\hat V$. The contact
order of $f$ along $\hat V$ when
$f(\hat V)\sub \hat \Delta_l^+$ is defined similarly.
We define the left contact order (resp. right contact order)
of $f$ at $p$, denoted
$\delta(f,p)\lef$ (resp. $\delta(f,p)\rig$),
to be the sum of the contact orders of $f$ along all irreducible components
of $\hat X$ that are mapped to $\hat \Delta_l^-$ (resp. $\hat \Delta_l^+$).

\begin{defi}
Let $f\mh X\to\wn$ be a non-degenerate pre-stable morphism.
We say $f$ is pre-deformable along $\bD_l$ if $\delta(f,p)\lef=\delta(f,p)\rig$
for all $p\in f\upmo(\bD_l)$.
We say $f$ is pre-deformable if $f$ is pre-deformable along all $\bD_l$.
\end{defi}

Note that if $f\upmo(\bD_l)=\emptyset$, then $f$ is automatically pre-deformable along $\bD_l$.

In the remainder part of this section, we fix a scheme $S$ over $C$ and a family
\begin{equation}
\label{2.30}
f\mh\cX\to\wn
\end{equation}
of non-degenerate pre-stable morphisms over $S$.
Our goal is to construct a maximal closed subscheme $S\pd\sub S$ so that
it parameterizes all pre-deformable morphisms in $S$.

We begin our investigation with the decomposition of $f$ along $\bD_l$.
We let $S_l$ be $S\times_{\cn}\cn_{[l\dual]}$
endowed with the reduced scheme structure. We have the following decomposition Lemma
whose proof follows from the standard deformation theory of nodal singularity of curves.

\begin{lemm}
\label{3.111}
Let the notation be as before. We let $\cX_l=\cX\times_S S_l$ and let $f_l$ be the restriction of
$f$ to $\cX_l$. Then there are two flat families of curves (not necessarily connected) $\cY_-$
and $\cY_+$ over $S_l$,
a pair of $S_l$-morphisms $g_{\pm}\mh \cY_{\pm}\to \Delta_l^{\pm}$,
a scheme $\Sigma$ finite and \'etale over $S_l$ and a pair of $S_l$-immersions
$\iota_{\pm}\mh \Sigma\to\cY_{\pm}$ of which the followings hold:
The original family $\cX_l$ is the gluing of $\cY_-$ and
$\cY_+$ along $\iota_-(\Sigma)\sub \cY_-$ and $\iota_+(\Sigma)\sub\cY_+$,
and the restriction of $f_l$ to $\cY_{\pm}\sub\cX_l$ is $g_{\pm}\mh\cY_{\pm}\to \Delta_l^{\pm}\sub \wn$.
\end{lemm}

Later, we will write $\cX=\cY_-\adj\cY_+$ and $f=g_-\adj g_+$ to indicate that $(\cY_{\pm},g_{\pm})$
are the decomposition of $(f,\cX)$.

\begin{lemm}
\label{2.63}
Let $S_{pd,l}\sub S$ be the set of all closed points $\xi\in S$ so that $f\lxi$ is pre-deformable
along $\bD_l$. Then $S_{pd,l}$ is closed in $S$.
\end{lemm}

\begin{proof}
It suffices to prove the case where $S$ is reduced and irreducible, which we assume now.
Let $S\to \Ao$ be induced by $S\to\cn\mapto{\pi_l}\bA^{\!n+1}_l\cong\Ao$.
In case $S$ is dominant over $\Ao$, then Lemma \ref{2.4} implies that
$S_{pd,l}=S$ as topological spaces and hence the Lemma holds. Now assume $S\to \Ao$ factor through
$0\in \Ao$. Since $S$ is reduced, $S\to\cn$ factor though $\cn_{[l\dual]}\sub \cn$
and then $S=S_l$ as schemes.
Let $\iota_\pm\mh \Sigma\to\cY_\pm$ and $g_\pm\mh\cY_\pm\to \Delta_l^\pm$ be
the decomposition of $f$ provided by Lemma \ref{3.111}.
Let $\pi\mh\Sigma\to S$ be the projection.
We let $\iota\mh\Sigma\to\cX$ be the composite $\iota_\pm\mh\Sigma\to\cY_\pm\to\cX$.
As before, for any $\xi\in S$ we denote
by $f\lxi$ the restriction of $f$ to the fiber over $\xi$ and denote by $\Sigma\lxi$
the fiber of $\Sigma$ over $\xi$.
Then the left and the right contact order
of $f$ along $\iota(\Sigma)$ define two functions
\begin{equation}
\label{2.51}
\delta(f_{\pi(\cdot)},\iota(\cdot))^{\mp}:\Sigma\lra\ZZ.
\end{equation}
It follows from the definition that
$\xi\in S_{pd,l}$ if the following two conditions hold:
\newline
(i) {\sl For any $x\in \Sigma\lxi$, $\delta(f\lxi,\iota(x))^+=\delta(f\lxi,\iota(x))^-$};
\newline
(ii) {\sl The preimage $f\lxi\upmo(\bD_l)$ is identical to $\iota(\Sigma\lxi)$}.
\newline
Applying the usual semi-continuity theorem, we see that both $\delta(f_{\pi(\cdot)},\iota(\cdot))^{\pm}$
are upper semi-continuous. Hence (i) and (ii) combined implies that
$S_{pd,l}$ is a constructible set.
Therefore, to show that $S_{pd,l}$ is closed it suffices to show that if
$R$ is a discrete valuation domain with residue field $\kk$ and quotient field $K$,
and if $\spec K\to S_{pd,l}$ extends to $\spec R\to S$,
then the extension $\spec R\to S$ factor through $S_{pd,l}\sub S$.

Now we let $\spec R\to S$ be a morphism so that $\spec K$ factor through $S_{pd,l}\sub S$.
We let $f_R\mh\cX_R\to \wn$ be the pull back family over $\spec R$ via $\spec R\to S$.
By abuse of notation, we still denote by
$\cY_\pm$, $g_\pm$ and $(\iota_\pm,\Sigma)$
the decomposition data of the $R$-curve $(f_R,\cX_R)$.
Let $\eta=\spec K$ be the generic point of $\spec R$. Then
by assumption the conditions (i) and (ii) hold for $\xi=\eta$.
Now let $\xi\in\spec R$ be any element.
As before, for $\alpha=-$ or $+$ we let $g_{\alpha,\xi}$ and $\cY_{\alpha,\xi}$ be the restrictions of the respective
families to the fiber over $\xi$.
 We let $g_{\alpha,\xi}\sta(\bD_l)$ be the pull-back divisor where $\bD_l$ is viewed as a divisor
in $\Delta_l^\alpha$.
Since $f_R$ is non-degenerate,
\begin{equation}
\label{2.25}
\deg(g_{\alpha,\xi}\sta(\bD_l),\cY_{\alpha,\xi})\geq\sum_{x\in \iota(\Sigma)\cap \cX_\xi}
\delta(f\lxi,\iota(x))^\alpha,\quad \alpha=-\ {\rm and}\ +.
\end{equation}
Clearly, the equality in (\ref{2.25}) holds for $\alpha$ if $g\sta_{\alpha,\xi}(\bD_l)=
\iota\lalp(\Sigma)\cap\cX\lxi$. Hence
condition (ii) is equivalently to that the equalities in the above two
inequalities hold.
Because the function of contact order is upper semi-continuous,
for the closed point $\eta_0\in\spec R$,
$$\sum_{x\in \iota(\Sigma)\cap \cX_{\eta_0}}\delta(f_{\eta_0},\iota(x))^{-}
\geq
\sum_{x\in \iota(\Sigma)\cap \cX_{\eta}}\delta(f_{\eta},\iota(x))^{-}.
$$
On the other hand, since $f$ is a flat family of non-degenerate morphisms, we have
$$\deg(g_{-,\eta_0}\sta(\bD_l),\cY_{-,\eta_0})=\deg(g_{-,\eta}\sta(\bD_l),\cY_{-,\eta}).
$$
Combined with the inequality in (\ref{2.25}) for $\alpha=-$ and $\xi=\eta_0$ and the equality in
(\ref{2.25}) for $\alpha=-$ and $\xi=\eta$,
we see that the equality in (\ref{2.25}) for $\alpha=-$ and $\xi=\eta_0$ must also hold. For similarly reason,
the equality in (\ref{2.25}) also hold for $\alpha=+$ and $\xi=\eta_0$.
Hence the two contact order functions in (\ref{2.51})
must be constant functions. Therefore, (i) and (ii) hold at $\eta_0$ and hence
$\spec R\to S$ factor through $S_{pd,l}\sub S$.
This shows that $S_{pd,l}$ is closed in $S$.
\end{proof}

Let $S\pd$ be the subset of $S$ consisting of $\xi\in S$ so that $f\lxi$ is pre-deformable.
Clearly, $S\pd=\cap_{l=1}^{n+1} S_{pd,l}$. Hence $S\pd$ is closed in $S$.
In the remainder of this section, we will give $S\pd$ a canonical closed scheme structure.

We begin with the notion of parameterizations of the formal neighborhood of nodes of
pre-stable curves. Let $\pi\mh\cX\to S$ be a family of pre-stable curves and let $\cX\lnode\sub\cX$
be the locus of all nodes of the fibers of $\cX$ over $S$, endowed with the reduced scheme structure.
Let $p\in\cX\lnode$ be any closed point and let $V\sub\cX\lnode$ be an affine neighborhood of
$p\in\cX\lnode$. Then by shrinking $V$ if necessary
we can find a scheme $U$ containing $V$ as its closed subscheme
so that there is an \'etale $r\mh U\to S$ extending $\pi|_V\mh V\to S$.\footnote{Such $(U,r)$
can be constructed by elementary method. For instance, this follows from the proof of Lemma 1.5
in \cite{Li}.}
We let $\hat V$ be the formal completion of $U$ along $V$.
We now argue that $\hat V$ is independent of the choice of $U$.
Indeed, in case $r\pri\mh U\pri\to S$ is another
such $U$ and $\hat V\pri$ is the formal completion of $U\pri$ along $V$. Then since
both $\hat U\to S$ and $\hat U\pri\to S$ are \'etale, by the topological invariance of \'etale
morphisms \cite[p30]{Mil}\footnote{In \cite{Mil} this theorem was stated for the case
$V\sub\hat V$ is defined by a
nilpotent ideal. However, the proof can be adopted without change to cover our case.},
the scheme $\hat V$ is canonically isomorphic to $\hat V\pri$. Thus $\hat V$ is
independent of the choice of $U$. Now let $\cX_{\hat V}$ be $\cX\times_S\hat V$.
Clearly, $V\sub\cX\lnode$ is a closed subscheme of $\cX_{\hat V}$.
We let $\hat\cX_{\hat V}$ be the formal completion of $\cX_{\hat V}$ along $V\sub\cX_{\hat V}$.

\begin{lemm}
\label{2.231}
Let $\pi\mh\cX\to S$ be a family of pre-stable morphisms as before and let $p\in\cX\lnode$
be any closed point. Then there is an affine neighborhood $V\sub\cX\lnode$ of $p$ so that
with $\hat V$ and $\hat\cX_{\hat V}$ constructed as before, we have a morphism
$\lbar{\phi}\mh \hat V\to\spec\ks$ and a $\hat V$-isomorphism
$$
\phi: \hat\cX_{\hat V}\lra\spec \bl\kz\otimes_{\ks}\Gamma(\cO_{\hat V})\br\rhat,
$$
where $(\cdot)\rhat$ is the $(z_1,z_2)$-adic completion of the respective ring.
We will call $\phi$ a parameterization of the formal neighborhood of the nodes of
$\cX$ along $V$. Furthermore, all such parameterizations satisfy the following uniqueness property.
Let $\phi\pri$ be another parameterization of $\cX$ along $V$.
Then there is a unit $u$
in $\bl\kz\otimes_{\ks}\Gamma(\cO_{\hat V})\br\rhat$ and a unit $\epsilon$ in $\Gamma(\cO_{\hat V})$ such that
$\uline{\phi\pri}\sta(s)=\epsilon \uline{\phi}\sta(s)$ and that the $\Gamma(\cO_{\hat V})$-automorphism
$(\phi\pri\circ\phi\upmo)\sta$ of $\hat R$ is define by
$$ (\phi\pri\circ\phi\upmo)\sta(z_1)=z_1 u\quad{\rm and}\quad
(\phi\pri\circ\phi\upmo)\sta(z_2)=z_2 u\upmo\epsilon.
$$
Finally, the choice of such $(u,\epsilon)$ is unique.
\end{lemm}

\begin{proof}
The existence of such a parameterization follows immediately from the deformation of nodes \cite{DM}.
This also follows from the local parameterization constructed in \cite[Section 1]{Li}.
Now let $\phi$ and $\phi\pri$ be two such parameterizations of $\cX$ along $V$.
Because $\lbar{\phi\pri}\sta(s)$ and $\lbar{\phi}\sta(s)$ define the same divisor in $\hat U$ corresponding to the
subscheme where the nodes are not smoothed, there is a unit $\epsilon\in\Gamma(\cO_{\hat V})$
such that $\lbar{\phi\pri}\sta(s)=\epsilon\lbar{\phi}\sta(s)$. Hence
$\phi^{\prime\ast}(z_1z_2)=\epsilon \phi\sta(z_1z_2)$ and hence
\begin{equation}
\label{3.220}
(\phi\pri\circ\phi\upmo)\sta(z_1z_2)=(\phi^{\ast-1}\circ \phi^{\prime\ast})(z_1z_2)=\epsilon z_1z_2.
\end{equation}
It follows from the deformation of nodes that there are units $u_1$ and $u_2$ in $\hat R$
such that
\begin{equation}
\label{2.230}
(\phi\pri\circ\phi\upmo)\sta(z_1)=z_1u_1
\quad{\rm and}\quad
(\phi\pri\circ\phi\upmo)\sta(z_2)=z_2u_2.
\end{equation}
Using the above relation, we obtain $su_1u_2=s\epsilon$. Hence
$u_2=\epsilon(u_1\upmo+e)$ for some $e\in \hat R$ satisfying $se=0$.
We now demonstrate that if we choose $u_1$, $u_2$ and $\epsilon$ appropriately, then
$e$ will be zero. We first write $e$ in its normal form
$$e=\alpha_0+\sum_{i\geq 1} z_1^i\alpha_i+z_2^i\beta_i,\quad \alpha_i,\beta_i\in \Gamma(\hat U).
$$
Since $se=0$, by the uniqueness of the normal form we have $z_1\beta_i=0$ and $z_2\alpha_i=0$ for $i\geq 1$ while
$s\alpha_0=0$. We let $u_1=\gamma_0+\sum z_1^i\gamma_i+z_2^i\delta_i$ be its normal form. Since $u_1$ is
a unit, $\gamma_0$ must be a unit. We now let $\epsilon\pri=\epsilon+\gamma_0\alpha_0$. Then since $s\alpha_0=0$,
we can replace $\epsilon$  in (\ref{3.220}) by $\epsilon\pri$. Then the new error term
$e\pri=u_2-u_1\upmo\epsilon\pri$, which obeys $u_2=\epsilon\pri(u_1\upmo+e\pri)$,
belongs to the ideal $(z_1,z_2)$. Hence without loss of generality we can
assume the $\epsilon$, $u_1$ and $u_2$ are chosen appriori that the $\alpha_0$ in the normal form of $e$ is zero.
On the other hand, since $z_2\alpha_i=0$ for $i>0$, we can replace $u_2$ by $u_2\pri=u_2-\sum_{i>0}z_1^i\alpha_i$
so that (\ref{3.220}) and (\ref{2.230}) still hold. Hence we can assume appriori that all $\alpha_i=0$.
Lastly, we let $u_1\pri=u_1(1+eu_1\epsilon\upmo)\upmo$. Clearly, $z_1u_1=z_1u_1\pri$ and $u_2=\epsilon u_1^{\prime-1}$
since $u_2=u_1\upmo\epsilon+e$. Thus $(\epsilon,u)$ with $u=u_1\pri$ is what we want.

The uniqueness of $(u,\epsilon)$ can be proved following same strategy.
Since a more general version will be proved in \cite[Section 1]{Li},
we will omit the proof here.
\end{proof}

We keep the integer $l\in[n]$ and continue to use the notation developed
so far. Let $\pi\mh\Sigma\to S_l$ be the projection and let $p\in\Sigma$ be any closed
point so that $\pi(p)\in S_{pd,l}$.
We pick a connected neighborhood $V$ of $p\in\Sigma$ and a
local parameterization $\phi$ of the formal neighborhood of the nodes of $\cX$ along $V$.
Let $A=\Gamma(\cO_{\hat V})$ and let $\hat V=\spec A$.
Let $\hat R$ be the $(z_1,z_2)$-adic completion of $\kz\otimes_{\ks}A$.
Without loss of generality, we can assume that the image $f(V)$ is contained in an
admissible chart $(\cW,\psi)$ of $\wn$ constructed in the beginning of this subsection.
Let
\begin{equation}
\label{2.31}
\varphi: \kw\lra\hat R
\end{equation}
be the morphism induced by $f$, the parameterization $\phi$ and $\psi$, and the
projection $\Theta_l\to\spec\kk[w_1,w_2]$.
Let $\mm$ be the ideal of the point $\xi\in V$, where $\xi=\pi(p)$.
Since $f\lxi$ is a pre-deformable morphism,
the left and the right contact order of $f\lxi$ at $p$ coincide, which we denote by $n$. Hence
possibly after exchanging $z_1$ and $z_2$, near $p$ we have $f\lxi\sta(w_i)=c_iz_i^n$
for some units $c_i$ in $\kz$.
Hence in the normal form $\varphi(w_i)=\sum a_{i,j}z_1^j+b_{i,j}z_2^j$,
as in (\ref{2.6}), both $a_{1,n}$ and $b_{2,n}$ are not in $\mm$.
Hence for sufficiently small neighborhood $V_0$ of $p\in V$ both
$a_{1,n}$ and $b_{2,n}$ are units in $B=\Gamma(\cO_{V_0})$. Therefore we
can apply Lemma \ref{2.11} to obtain an ideal $I_B\sub B$ for such $B$.

\begin{defi}
\label{2.35}
Let the notation be as before.
Let $V\sub\Sigma$ be an affine neighborhood of $p\in\Sigma$ so that $f(V)$ is contained in an
admissible neighborhood $(\cW,\psi)$ of $f(p)\in\wn$.
We say $f$ is pre-deformable near $p$ if for some neighborhood $V_0$ of $p\in V$
both $a_{1,n}$ and $b_{2,n}$ are units in $B=\Gamma(\cO_{V_0})$ and further
the ideal $I_B\sub B$ just constructed is the zero ideal in $B$.
The family $f$ is said to be pre-deformable along $\bD_l$ if it is pre-deformable
at every $p\in f\upmo(\bD_l)$. We say $f$ is pre-deformable if it is pre-deformable along
all divisors $\bD_l$.
\end{defi}

A remark about this definition is in order. First, in defining the notion of pre-deformable at $p$
we need to choose a parameterization of the formal
neighborhood of the nodal singularity of $\cX$ along $V$.
Because of the second part of Lemma \ref{2.231}, this
definition does not depend on this choice.
For the same reason, this definition is also independent of the choice of $(\cW,\psi)$.

We still denote by $f$ a family of non-degenerate morphisms to $\wn$ over $S$.
We will conclude this section by showing that there is a canonical closed subscheme structure on
$S\pd\sub S$ making it the closed subscheme parameterizing all pre-deformable morphisms over $S$.

\begin{defi}
Let $f$ be a family of non-degenerate morphisms to $\wn$ over $S$. We define a functor
$\F_{pd,l}[f,S]$ (resp. $\F_{pd}[f,S]$) that associates to any scheme $T$ over $C$ the set of morphisms
$T\to S$ so that the pull-back family $f_T\mh\cX_T\to\wn$ is pre-deformable along the divisor $\bD_l$
(resp. all divisors $\bD_l$).
\end{defi}

When the choice of $f$ is apparent from the context, we will omit the reference $f$
and abbreviate $\F_{pd,l}[f,S]$ to $\F_{pd,l}[S]$.

\begin{theo}
Let $f\mh\cX\to\wn$ be a family of non-degenerate stable morphisms. Then
the functor $\F_{pd,l}[S]$ and $\F_{pd}[S]$ are represented by closed subschemes of $S$.
\end{theo}

\begin{proof}
Clearly, should $\F_{pd,l}[S]$ be represented by a subscheme in $S$, then set-theoretically it
must be $S_{pd,l}$. Hence it suffices to give
$S_{pd,l}$ a closed scheme structure so that it represents the
functor $\F_{pd,l}[S]$.
As before, we let $S_l$ be
$S\times_{\Anpo}\bH_l\unpo$ endowed with the reduced scheme structure and $\hat S_l$ be the formal completion
of $S$ along $S_l$. We let $\hat V$ be any connected component of $\hat S_l$ endowed with the open
scheme structure and let $V$ be $\hat V$ endowed with the reduced scheme structure. Our first objective
is to show that $S_{pd,l}\cap \hat V$ admits a canonical scheme structure (as closed scheme of $\hat V$)
so that it represents the functor $\F_{pd,l}[\hat V]$. We continue to use the notation developed
so far. For instance
$(\cY_\pm,g_\pm)$ and $(\iota_\pm,\Sigma)$
are the decomposition data of the restriction of $f$ to $V\sub S$ along the divisor $\bD_l$.
We let $\hat \Sigma\to\hat V$ be the \'etale morphism so that
$\hat \Sigma\times_{\hat V} V\cong\Sigma$. We let $\pi\mh\Sigma\to V$ and $\hat \pi\mh\hat \Sigma\to\hat V$
be the projections.

Let $\xi\in S_{pd,l}\cap V$ be any closed point and let $p_1,\cdots,p_k$ be the set of points in $\pi\upmo(\xi)$.
Note that $p_1,\cdots,p_k$ are exactly the preimage set $f\lxi\upmo(\bD_l)$. Hence each $p\lalp$ is
associated with an integer $n\lalp$ that is the contact order of $f\lxi$ at $p\lalp$.
For each $\alpha $ between $1$ and $k$, we pick an affine open neighborhood
$\hat U\lalp$ of $p\lalp\in\hat\Sigma$.
We let $\cX\lalp=\cX\times_S \hat U\lalp$ and let
$f\lalp\mh\cX\lalp\to\wn$ be the pull back of $f$ to $\cX\lalp$.
We let $U\lalp$ be $\hat U\lalp\cap\Sigma$. Then the multi-section
$\iota\mh\Sigma\to\cX$ defines a section $\iota\lalp\mh U\lalp\to\cX\times_{\hat U\lalp} U\lalp$.
Without loss of generality, we can assume $\cX\lalp$ admits a parameterization of the formal
neighborhood of its nodes along $\iota\lalp(U\lalp)$. By shrinking $U\lalp$ if necessary, we can assume
that there is an admissible chart $(\cW\lalp,\psi\lalp)$ of $\wn$ near $f\lxi(p\lalp)$
so that $f\lalp(\iota\lalp(U\lalp))$ is entirely contained in $\cW\lalp$.
Therefore, following the construction associated to Definition \ref{2.35},
by shrinking $U\lalp$ if necessary we have
an ideal sheaf $\cI\lalp\sub\cO_{\hat U\lalp}$ so that it defines the subscheme of $\hat U\lalp$
that parameterizes all morphisms in the family  $f\lalp\mh\cX\lalp\to \wn$ that are pre-deformable and have
contact order $n\lalp$ along the nodes $\iota\lalp(U\lalp)\sub\cX\lalp$.

To obtain the ideal sheaf of $\cO_{\hat V}$ that defines the desired closed subscheme structure of
$S_{pd,l}\cap\hat V$ we need to descend the ``intersection'' of $\cI\lalp$ since $S_{pd,l}$ is defined by the
requirement that all nodes are pre-deformable.
As before, we denote by $k$ the number of elements in $\pi\upmo(\xi)$ for $\xi\in V$. Since $V$ is connected
and $\Sigma\to V$ is \'etale and finite, $k$ is independent of the choice of $\xi\in V$.
We consider the fiber product, denoted $\hat \Sigma^{\times k}$,
of $k$-copies of $\hat \Sigma$ over $\hat V$. We denote by $\Sigma^{\times k}$
accordingly the fiber product of $k$ copies of $\Sigma$ over $V$.
Note that points in $\Sigma^{\times k}$ are $(a_1,\cdots,a_k)$ with not necessarily distinct
$a_1,\cdots,a_k\in\pi\upmo(\xi)$
for some $\xi\in V$. We let $\Gamma\sub \Sigma^{\times k}$ be the subset of those points
$(a_1,\cdots,a_k)\in \Sigma^{\times k}$
so that all $a_1,\cdots,a_k$ are distinct. Since $\Sigma\to V$ is finite and \'etale, and that
each fiber of $\pi\mh\Sigma\to V$ has cardinal $k$, $\Gamma$ is an open subset of $\Sigma^{\times k}$.
Then since $\Sigma^{\times k}$ is homeomorphic to $\hat \Sigma^{\times k}$, $\Gamma$ is an open
subset of $\hat \Sigma^{\times k}$. We let $\hat \Gamma$ be $\Gamma\sub \hat \Sigma^{\times k}$ endowed with
the open scheme structure. Then the obvious projection $\hat \Gamma\to \hat V$
is \'etale and finite. Now let $\hat \pr\lalp\mh\hat \Gamma\to \hat \Sigma$ be the morphism induced by the $\alpha$-th
projection of $\hat \Sigma^{\times k}$, which is \'etale.

Now consider any point $\eta=(p_1,\cdots,p_k)\in\hat\Gamma$. To each $\alpha$ between $1$ and $k$ we
let $\hat U\lalp\sub\hat \Sigma$ be the neighborhood of $p\lalp\in\hat\Sigma$ and $\cI\lalp\sub\cO_{\hat U\lalp}$
the ideal sheaf just constructed. We let $\hat\cU\sub\hat\Gamma$ be an open neighborhood of $\eta$
so that $\hat \pr\lalp(\hat\cU)\sub\hat U\lalp$ for each $\alpha$. We then let $\tilde\cI\lalp$ be
the ideal sheaf of $\cO_{\hat\cU}$ that is the pull back of $\cI\lalp\sub\cO_{\hat U\lalp}$ under
$\hat\pr\lalp|_{\hat \cU}\mh \hat\cU\to\hat U\lalp$,
and let $\cJ\sub\cO_{\hat\cU}$ be the intersection $\cJ=\tilde\cI_1\cap\cdots\cap\tilde\cI_l$.
This ideal sheaf defines a close subscheme of $\hat\cU$, denoted by $Z(\cJ)$.
Now let $\Phi\mh\hat \Gamma\to\hat S$ be the composite $\hat\Gamma\to\hat V\to\hat S$.
It may happens that $Z(\cJ)$ and $\Phi\upmo(S_{pd,l})$ are different as subsets of $\hat\Gamma$.
However, by the proof of the closedness of $S_{pd,l}$, the connected component of $Z(\cJ)$ containing
$\eta$ coincides with the connected component of $\Phi\upmo(S_{pd,l})$ containing $\eta$. Hence by shrinking $\hat \cU$
if necessary, we can assume $Z(\cJ)$ is connected.
Now let $f_{\hat\cU}$ be the family over $\hat\cU$ that is the pull-back of the family $f\mh\cX\to\wn$
via $\hat\Gamma\mapright{\Phi}\hat S\to S$. Since for $\xi\in\hat \cU$ the preimage $f\lxi\upmo(\bD_l)$ has
exactly $k$ points, it is direct to check that the subscheme $Z(\cJ)$
represents the functor $\F_{pd,l}[\hat \cU]$.

Now by repeating this procedure, we can cover $\hat \Gamma$ by open subsets $\hat\cU\lbe$
with ideal sheaves $\cJ\lbe\sub\cO_{\hat\cU\lbe}$, for $\beta=1,\cdots,L$,
such that the closed subscheme $Z(\cJ\lbe)$ defined
by the ideal sheaf $\cJ\lbe$ represents the functor $\F_{pd,l}[\hat\cU\lbe]$.
By the universal property of these schemes,
whenever $\hat\cU\lalp\cap\hat\cU\lbe\ne\emptyset$ the ideal sheaves $\cJ\lalp$ and $\cJ\lbe$
coincide over $\hat\cU\lalp\cap\hat\cU\lbe$.
Thus they patch together to form an ideal sheaf $\cJ\sub\cO_{\hat\Gamma}$. Again, by the universal property, this
ideal sheaf forms a descent data for the faithfully flat \'etale morphism $\hat\Gamma\to\hat V$. Hence it
descends to an ideal sheaf $\cI\sub\cO_{\hat V}$
which defines a closed subscheme $Z(\cI)$.
It follows from the previous construction that $Z(\cI)$ represents the functor $\F_{pd,l}[\hat V]$.

We now show that the functor $\F_{pd,l}[S]$ is represented by a closed subscheme.
Since each connected component $\hat V$ of $\hat S_l$ has a closed scheme structure
that represents the functor $\hat V$, the functor $\F_{pd,l}[\hat S_l]$ is
represented by a closed subscheme of $\hat S_l$. We denote this
subscheme by $\hat S_{pd,l}$.
It remains to show that there is a closed scheme structure on $S_{pd,l}$ so that it
represents the functor $\F_{pd,l}[S]$.
As before, we let $t_l$ be the $l$-th standard parameter of $\An$.
We let $\cK\sub\cO_S$ be the subsheaf of elements
annihilated by some power of $t_l$. Since $S$ is noetherian, there is an integer $m_0$ so that
$t^{m_0}_l\cK\sub \cO_S$ is zero. Similarly, we let $\hat\cK$ be the subsheaf of $\cO_{\hat S}$
consisting of elements annihilated by some power of $t_l$. Then $\hat\cK$ is canonically isomorphic
to $\cK$. As an ideal sheaf, $\hat \cK$ defines a closed subscheme $\hat S_{fl}\sub\hat S_l$.
Clearly, $\hat S_{fl}$ is flat over $\spec\kk[t_l]$ near $0$. Hence by Lemma \ref{2.4},
the restriction of the family $f$ to $\hat S_{fl}$ is pre-deformable along $\bD_l$.
Consequently, the ideal sheaf of $\hat S_{pd,l}\sub \hat S_l$ is contained in $\hat \cK$.
Using the natural isomorphism $\hat \cK\cong\cK$, this ideal sheaf defines an ideal sheaf $J$ of
$\cO_S$. By our construction the support of the subscheme defined by $J$ coincides with $S_{pd,l}$.
We define the closed scheme structure of $S_{pd,l}$ to be the one defined by $J$.
It is straightforward to check that with this closed scheme structure, $S_{pd,l}$ represents the
functor $\F_{pd,l}[S]$.

As to $\F_{pd}[S]$, it is clear that it is represented by the intersection of the closed subschemes
$\cap_{l=1}^n S_{pd,l}$. This completes the proof of the Theorem.
\end{proof}

\section{Moduli stack of stable morphisms}

In this section, we will introduce the notion of stable morphisms
to the stack $\WWc$. We will then show that the moduli of such stable
morphisms form an algebraic stack (i.e. a Deligne-Mumford stack).

\subsection{Stable morphisms}

We first define the notion of stable morphisms in our setting.
Let $S$ be any $C$-scheme. A pre-stable morphism to $\WWc(S)$ consists of a
triple $(f,\cX,\cW)$ where $\cX$ is a family\footnote{All families are flat families unless
otherwise is mentioned.}
of marked connected, pre-stable curves
over $S$, $\cW$ is a member in $\WWc(S)$ and $f$ is an $S$-morphism $f\mh\cX\to\cW$ that
is stable as a family of morphisms to $\cW$.
When there is no confusion, we will abbreviate the triple $(f,\cX,\cW)$ to $f$.
Let $(f\pri,\cX\pri,\cW\pri)$ be another pre-stable morphism over $S$.
An isomorphism from $f$ to $f\pri$ consists of an $S$-isomorphism $r_1\mh\cX\to\cX\pri$
(as pointed curves) and
an arrow $r_2\mh\cW\to\cW\pri$ so that $r_2\circ f=f\pri\circ r_1$.
When $f\equiv f\pri$, then such isomorphisms are called the automorphisms
of $f$. To define the group scheme of the automorphisms of $f$, we follow
the usual procedure. Let $\AUT_{{\mathfrak W}}(f)$ be the functor that associates any $S$-scheme $T$
the set of all automorphisms of $f\times_ST$. $\AUT_{{\mathfrak W}}(f)$ is represented by a group scheme
over $S$, denoted $\Autw(f)$ and called the automorphism group of $f$.

Let $f$ be a map $f\mh X\to \wn$ with $X$ proper and reduced.
We fix an ample line bundle $H$ on $W$. Using $\wn\to W$, $H$ pulls back to a line bundle over $\wn$,
which we denote by $H$ if there is no confusion.
We define the degree of $f$ to be the degree of $f\sta H$ over $X$.
We fix a triple of integers $\Gamma=(b,g,k)$ where $k$ will be the number of marked points of the domain curve,
$g$ will be the arithmetic genus of the domain curve and $d$ will be
the degree of the map.
We say $f\mh X\to \wn$ is a pre-stable map of topological type $\Gamma$ if
$f\mh X\to\wn$ is an ordinary stable morphism whose domain is
a connected, $k$-pointed arithmetic genus $g$ nodal curve and $f$ has degree $b$.
We fix such an $H$ and $\Gamma$ once and for all in this section.

\begin{defi}
\label{31.1}
Let $S$ be any $C$-scheme and $f\mh\cX\to\cW$ be a family of pre-stable morphisms
to $\WWc(S)$. Suppose $\cW$ is an effective degeneration, say given by a
$C$-morphism $r\mh S\to\cn$. Then we say $f$ is a family of pre-deformable morphisms
to $\WWc(S)$ if the associated morphism $\tilde f\mh\cX\to\wn$ is a family of pre-deformable
morphisms; We say $f$ is a family of stable morphisms to $\WWc(S)$ if in addition to
$f$ being a family of pre-deformable morphisms we have that for every
point $\xi\in S$ the automorphism group $\Autw(f_{\xi})$ of $f\lxi\mh\cX\lxi\to \wn$
is finite. We say $f$ is a family of topological type $\Gamma$ if for every $\xi\in S$
the restriction of $f$ to the fiber of $\cX$ over $\xi$ is of topological type
$\Gamma$. In general, given a family of pre-stable morphisms $f$ to $\WWc(S)$,
we say that $f$ is a family of pre-deformable (resp. stable, resp. of topological
type $\Gamma$)
morphisms to $\WWc(S)$ if there is an open covering $S\lalp$ of $S$ such that to each $\alpha$ the family
$\cW\times_S S\lalp$ is an effective family over $S\lalp$ and the restriction of $f$ to $\cX\times_S S\lalp$
is a family of pre-deformable (resp. stable, resp. of topological type $\Gamma$) morphisms to $\WWc(S\lalp)$.
\end{defi}

Clearly, the definition of stability does not depend on the choice of the local representatives of
the family $\cW$.
The goal of this section is to construct the moduli stack of all marked stable morphisms to
$\WWc$ of given topological type $\Gamma$. In the sequel of this paper, we will
use this moduli stack to study the degeneration of GW-invariants.

We begin with the definition of the groupoid of stable morphisms to $\WWc$.
We fix the relative ample line bundle $H$ on $W$ and the triple $\Gamma$.
We define $\mgwc$ to be the category whose objects
are all families of stable morphisms
to $\WWc$ of topological type $\Gamma$.
We define ${\mathfrak p}\mh\mgwc\to ({\rm Sch}/C)$ to be the morphism that sends families over $S$ to $S$.
Let $\xi=\{f\mh\cX\to\cW\}\in\mgwc(S)$, $\xi\pri=\{f\pri\mh\cX\pri\to\cW\pri\}\in\mgwc(S\pri)$
and let $h\mh S\pri\to S$ be an arrow (a morphism). An arrow $r\mh\xi\pri\to\xi$
with ${\mathfrak p}(r)=h$ consists of an arrow $\cW\pri\to\cW$ covering $h$
(i.e. a $W_{S\pri}$-isomorphism $\cW\times_{S\pri} S\pri\cong\cW\pri$) and
an $S\pri$-isomorphism $\cX\times_S S\pri\cong\cX\pri$ (as marked
curves) so that
$$
\begin{CD}
\cX\pri @>{f\pri}>> \cW\pri\\
@VV{\cong}V @VV{\cong}V \\
\cX\times_SS\pri @>{f\times_S S\pri}>> \cW\times_S S\pri
\end{CD}
$$
is commutative. In such case, we call $\xi\pri$ the pull-back of $\xi$ via $h\mh S\pri\to S$.
Obviously, this makes $(\mgwc,\mathfrak p)$ a groupoid.

Before proceeding, we need to study the structure of stable morphisms to $\WWc$.
Let $f\mh X\to \wn$ be a pre-deformable morphism with $B\sub X$ the divisor of the
marked points. We assume that $f$ factor through $\wn_0\sub\wn$.
By the decomposition Lemma, $f$ splits $X$, according to the nodal divisors of $\wn_0$,
into components $X_1,\cdots,X_{n+2}$
(some of them may be empty), marked points $\iota_i^+\mh\Sigma_i\to X_{i+1}$
and $\iota_i^-\mh \Sigma_i\to X_{i}$, where $\Sigma_i$ is a finite set,
so that $f(X_i)\sub\Delta_i\sub\wn_0$ and that
$X$ is the result of successively gluing
$X_i$ and $X_{i+1}$ along $\iota_i^+\mh\Sigma_i\to X_{i+1}$ and $\iota_i^-\mh\Sigma_i\to X_i$.
Namely,
$$X=X_1\adj\cdots\adj X_{n+2}\quad {\rm and}\quad
f=f_1\adj\cdots\adj f_{n+2}.
$$
Now let $A$ be any purely one-dimensional closed subset of $X$. We define the weight of $f$ along $A$ to be
$$\omega(f,A)=\deg(A,f\sta H)+2g(A)-2+\#(B\cap A)+\#\tau(A),
$$
where $\tau(A)$ is the set of smooth points of $A$ that are nodes of $X$.
(We agree $\omega(f,A)=0$ if $A=\emptyset$.)
Note that if $A_1$ and $A_2$ are two such purely one-dimensional closed subsets with no
common irreducible components, then
$$\omega(f,A_1\cup A_2)=\omega(f,A_1)+\omega(f,A_2).
$$

\begin{lemm}
\label{3.90}
The pre-deformable morphism $f\mh X\to\wn_0$ is a stable morphism in $\mgwc(\kk)$\footnote{
For simplicity we agree $\mgwc(\kk)=\mgwc(\kk)$.}
if and only if $\omega(f,X_i)>0$ for all $i=2,\cdots,n+1$.
\end{lemm}

\begin{proof}
Let $i$ be any integer between $2$ and $n+1$. We let
$$\Psi_i\mh G_m\times \Delta_i\lra \Delta_i
$$
be the group action given in Corollary \ref{1.111}. For any closed
purely one-dimensional subset $A$ of $X_i$ we define the
relative automorphism group $\Autw(f,A)\rel$ to be the set of pairs $(a,b)$ where $a\mh A\to A$ is an
automorphism leaving $B\cap A$ and $\tau(A)$ fixed and $b$ is an
element in $G_m$ such that
$\Psi_i(b,f_i(w))=f_i(a(w))$ for all $w\in A$.
Let $\sigma\in\Autw(f)$ be any automorphism. Then $\sigma$ induces a permutation
of $\Sigma_i$ for each $i$. As a result, it defines a homomorphism
$$\Autw(f)\lra S_{\Sigma_1}\times\cdots\times S_{\Sigma\lnpo},
$$
where $S_{\Sigma_i}$ is the permutation group of $\Sigma_i$. Clearly, its kernel is
$$\Autw(f,X_1)\rel\times\cdots\times\Autw(f,X_{n+2})\rel.
$$
Since $f$ is an ordinary stable morphism, $\Autw(f,X_1)\rel$ and $\Autw(f,X_{n+2})\rel$
are finite.
Hence $\Autw(f)$ is finite if and only if $\Autw(f,X_i)\rel$ are finite for $i=2,\cdots,n+1$.

Now assume $A$ is a connected component of $X_i$ so that $f(A)$ is not a single point set and
that $\Autw(f,A)\rel$ is infinite. Then the following must hold:
$f(A)$ must be contained in a fiber of $\Delta_i\to D$;
$A$ must be a smooth rational curve; $A$ contains no marked points and $\tau(A)$ consists of
exactly two nodal
points of $X$, one lies in $\iota_{i-1}^+(\Sigma_{i-1})$ and the other lies in $\iota^-_i(\Sigma_i)$. Further if we
let $\Po\sub\Delta_i$ be the image of $A$ then $f|_A\mh A\to\Po$ is a branched covering
ramified at $\tau(A)$. We call such connected components
{\sl trivial components}. Clearly, if $A$ is not a trivial component, then $\Autw(f,A)\rel$ is
finite. On the other hand, if $A$ is a trivial component then $\Autw(f,A)\rel\to G_m$ is
surjective. Hence $\Autw(f,X_i)\rel$ is finite if and only if not all connected components of
$X_i$ are trivial.

To complete the proof of the Lemma, it suffices to show that $\omega(f,X_i)=0$ if and only if
all connected components of $X_i$ are trivial.
Let $A$ be any connected component of $X_i$. In case $\deg(A,f\sta H)>0$,
then since $H$ is very ample we can assume without loss of generality that $\deg(A,f\sta H)>3$,
and hence $\omega(f,A)>0$. In case $\deg(A,f\sta H)=0$, then $f(A)$ must be contained in a fiber
of $\Delta_i\to D$. Since $\deg(X,f\sta H)>0$ and $X$ is connected, $f(A)$ can not be a single point.
Hence $f(A)$ is a fiber of $\Delta_i\to D$. Because $f$ is pre-deformable,
$A$ contains at least two nodal points of $X$, one lies in $f\upmo(D_{i-1})$ and the other lies in
$f\upmo(D_i)$. Hence $\omega(f,A)\geq 0$.
Moreover, the equality holds if and only if $g(A)=0$, $A$ contains no marked points and
contains exactly two nodal points of $X$.
Again since $f$ is pre-deformable, it is easy to see that this is possible only if $A$
is a trivial component.
This shows that any connected component $A$ of $X_i$ has
$\omega(f,A)\geq 0$ and the equality holds if and only if $A$ is a trivial component.

Now we assume $f$ is stable. Then to each $i$ between $2$ and $n+1$ at least one connected component of
$X_i$ is non-trivial. Then combined with the fact that
the weight of any connected components of $X_i$ are non-negative, we have $\omega(f,X_i)>0$.
Conversely, if for some $i$ between $2$ and $n+1$ the weight $\omega(f,X_i)=0$, then all connected components
of $X_i$ are trivial. Since to each component $A$ of $X_i$ the automorphism group
$\Autw(f,A)\rel$ surjects onto $G_m$, $\Autw(f,X_i)\rel$ will surjects onto $G_m$ as well.
Thus $\Autw(f,X)\rel$, and hence $\Autw(f)$ is infinite. This proves the Lemma.
\end{proof}

For convenience we define the norm of the triple $\Gamma=(b,g,k)$ to be
$\vert \Gamma\vert =b+2g-2+k$.

\begin{lemm}
\label{3.20}
Let $f\mh X\to\wn_0$ be a stable morphism in $\mgwc(\kk)$ and let $X_1,\cdots,
X_{n+2}$ be the splitting components of $X$. Then
$$\vert \Gamma\vert =\omega(f,X_1)+\cdots+\omega(f,X_{n+2}).
$$
\end{lemm}

We need a more general form of this Lemma. Let $S=\spec R$ be a $C$-scheme where
$R$ is a discrete valuation domain. Let
$f\mh\cX\to\cW$ be a family of stable morphisms in $\mwcg(S)$ and $\iota\mh S\to\cn$ be the morphism
associated to $\cW$.
Let $\eta_0\in S$ be the closed point and
$\eta\in S$ be the generic point.
Let $\Delta_1^0,\cdots,\Delta_{n_0+2}^0$ be the $(n_0+2)$-irreducible components of $\wn\times_{\cn}\iota(\eta_0)$
and $\Delta_1,\cdots,\Delta_{n_1+2}$ be the $(n_1+2)$-irreducible components of $\wn\times_{\cn}\iota(\eta)$.
We say $\Delta_i^0\preceq \Delta_j$ if $\Delta_i^0$ is contained in the closure of $\Delta_j$.

\begin{lemm}
\label{3.21}
Let the notation be as before and let $X_i^0$ (resp. $X_j$) be the splitting components of
$\cX\times_X\eta_0$ (resp. $\cX\times_S\eta$). Then for each $j\in[1,n_1+2]$,
$$\omega(f_{\eta},X_j)=\sum_{\Delta_i^0\preceq \Delta_j} \omega(f_{\eta_0}, X_i^0).
$$
\end{lemm}

\begin{proof}
The proof of these two Lemmas are standard and will be omitted.
\end{proof}

We next introduce the standard choice of ample line bundle on the domain and the target of
a stable morphism in $\mwcg$. We begin with line bundles on $\wn$. As before, we fix the relative
ample line bundle $H$ on $W/C$ and its pull back on $\wn$, still denoted by $H$.
For $l\in[n]$, the divisor $\wn\times_{\cn}\cn_{[l\dual]}\sub\wn$ is the union of two smooth irreducible components,
one {\sl left} and one {\sl right} according to our convention of ordering. We denote
the {\sl right} component by $\Xi_l$.
(Note that $\Xi_l\times_{\cn}\bz$ contains the last irreducible component of $\wn_0$.)
Then associating to any sequence of integers
$[a]=(a_1,\cdots,a_n)$ we define a line bundle
$$H_{[a]}=H(a_1\Xi_1+\cdots+a\lnpo\Xi\lnpo)
$$
on $\wn$. Since $\wn\times_C(C-0)\to W\times_C(C-0)$ is $\gn$-equivariant and since restricting
to $\wn\times_C(C-0)$ the line bundle $H_{[a]}$
is the pull-back of $H$ from $W\times_C(C-0)$, the trivial $\gn$-linearization of
$H$ over $W$ pull back to a unique $\gn$-linearization
of $H_{[a]}$ over $\wn\times_C(C-0)$. Further, this linearization extends to a $\gn$-linearization
of $H_{[a]}$ over $\wn$. We call this the canonical $\gn$-linearization of $H_{[a]}$.

Let $\Delta\to D$ be the ruled variety that is one of the middle components of $\wn_0$, where $n>1$.
We denote by $H_{\Delta}$ the pull back of the ample line bundle $H$ via $\Delta\to D\to W$.
Without lose of generality, we can assume $H$ on $W$ is sufficiently ample so that for any integer
$0\leq a<b\leq\vert \Gamma\vert$ the divisor $H_{\Delta}(aD_+ +bD_-)$ is ample on $\Delta$,
where $D_{\pm}$ are the two distinguished divisors of $\Delta$.
We now choose a standard choice of ample line bundle on the
target of $f\mh X\to\wn_0$.
Let $X_1,\cdots,X_{n+2}$ be the splitting components of $X$. We define
$$a_i(f)=\sum_{j=1}^i\omega(f, X_j),
$$
and define the standard line bundle on $\wn_0$ (associated to $f$) to be
$$H_{[f]}=H( a_1(f)\Xi_1+\cdots+a\lnpo(f)\Xi\lnpo).
$$
Because $\{a_i(f)\}$ is a strictly increasing sequence, $H_{[f]}$ is ample on $\wn_0$.

Let $S$ be any scheme and $f\mh\cX\to \cW$ be a family in $\mgwc(S)$.
We now extend this construction to the family $\cW/S$.
Let $\xi\in S$ be any closed point and let $U$ be an open neighborhood of $\xi\in S$
so that $\cW\times_S U$ is isomorphic to the pull-back $\rho\sta\wn$ for a
$\rho\mh U\to\cn$. Without lose of generality, we can assume $\rho(\xi)=0\in\cn$.
Let $\Xi_1,\cdots,\Xi_n$ be the divisors of $\wn$ mentioned before. Let
$f\lxi\mh\cX\lxi\to\wn_0$ be the restriction of $f$ to the fiber over $\xi$. We
then form the line bundle
$$H_{[f_{\xi}]}=H(a_1(f\lxi)\Xi_1+\cdots+a\lnpo(f\lxi)\Xi\lnpo)
$$
on $\wn$. Because of Lemma \ref{3.21}, after shrinking $U$ if necessary,
for any $\eta\in U$ the restriction of $H_{[f\lxi]}$ to $\cW_{\eta}=\wn\times_{\cn}\rho(\eta)$
is the standard choice of ample line bundle associated to $f_{\eta}$.
We define the standard choice of ample line bundle on $\cW\times_S U$ to be $\rho\sta H_{[f\lxi]}$
and denote it by $H_{[f,U]}$.
Repeating this procedure, we can find an open covering $U\lalp$ of $S$ so that
each $\cW\times_S U\lalp$ is isomorphic to $\rho\lalp\sta W[n\lalp]$ for some
$\rho\lalp\mh U\lalp\to C[n\lalp]$ and over $\cW\times_S U\lalp$ we have the standard
choice of ample line bundle $H_{[f,U\lalp]}$. By Lemma \ref{1.167} we can assume that
whenever $U\lalpbe=U\lalp\cap U\lbe\ne\emptyset$ then there is an arrow
$\varphi\lalpbe$ making the diagram
$$\begin{CD}
\rho\lalp\sta W[n\lalp]\times_{U\lalp} U\lalpbe
@>{\varphi_{\beta\alpha}}>> \rho\lbe\sta W[n\lbe]\times_{U\lbe} U\lalpbe\\
@VVV@VVV\\
\cX\times_S U\lalpbe @= \cX\times_S U\lalpbe\\
\end{CD}
$$
is commutative.
Using the $\gn$-linearization, such an arrow defines a canonical isomorphism
$$\psi_{\beta\alpha}: \varphi_{\beta\alpha}\sta H_{[f,U\lbe]}\cong H_{[f,U\lalp]}.
$$
Since $\psi_{\beta\alpha}$ is canonical, the collection $\{H_{[f,U\lalp]}\}$ forms a line bundle on $\cW$.
In the following, we will call the
collection $H_{[f,U\lalp]}$ the standard choice of relative ample line bundle on $\cW$ and denote it
by $H_{[f,\cW]}$.

The standard choice of the ample line bundle on the domain of $f$ is
$$H_{[f,\cX]}=f\sta H_{[f,\cW]}(5B)\otimes\omega_{\cX/S}^{\otimes 5},
$$
where $\omega_{\cX/S}$ is the relative dualizing sheaf and
$B\sub\cX$ is the divisor of the marked points.

\subsection{Basic Properties of $\mgwc$}

In this subsection, among other things, we will show that $\mwcg$ is a separated and proper
algebraic stack over $C$.

\begin{prop}
The groupoid $\mwcg$ over $C$ is a stack.
\end{prop}

\begin{proof}
Following the definition of stacks (cf. \cite{Vis}), it suffices to show the following:
\newline
(i) For any scheme $S$ over $C$ and two stable morphisms $\xi_1$ and $\xi_2\in\mwcg(S)$, the functor
$$\Isom_S(\xi_1,\xi_2): ({\rm Sch}/S)\lra ({\rm Sets})
$$
which associates to any morphism $\varphi\mh T\to S$ the set of isomorphisms in $\mwcg(T)$ between
$\varphi\sta\xi_1$ and $\varphi\sta\xi_2$ is a sheaf in the \'etale topology;
\newline
(ii) Let $\{S_i\to S\}$ be a covering of $S$ (over $C$) in the \'etale topology. Let
$\xi_i\in\mwcg(S_i)$ and let $\varphi_{ij}\mh\xi_i|_{S_i\times_S S_j}\to \xi_j|_{S_i\times_S S_j}$
be isomorphisms in $\mwcg(S_i\times_S S_j)$ satisfying the cocycle condition. Then there is a
$\xi\in\mwcg(S)$ with isomorphism $\psi_i\mh\xi|_{S_i}\to \xi_i$ so that
$$\varphi_{ij}=(\psi_i|_{S_i\times_S S_j})\circ(\psi_j|_{S_i\times_S S_j})\upmo.
$$

We first prove (i). Let $\xi_i$, $i=1$ and $2$, be represented by $f_i\mh \cX_i\to\cW_i$.
To each morphism $\rho\mh U\to S$ we define $\Isom_S(\xi_1,\xi_2)(U)$ to be the set of all arrows
from $\rho\sta\xi_1$ to $\rho\sta\xi_2$. This defines a functor
$$\Isom_S(\xi_1,\xi_2): (\hbox{\'etale opens in S})\lra ({\rm Sets}).
$$
Clearly it is a pre-sheaf.
To show that this is a sheaf it suffices to show that if $U=\cup_{j\in J} U_j$
is an \'etale open covering of $U$, then
$$
\Isom_S(\xi_1,\xi_2)(U) \mapright{e}
\prod_{j\in J} \Isom_S(\xi_1,\xi_2)(U_j)
\mapright{p_1,p_2}  \prod_{i,j\in J} \Isom_S(\xi_1,\xi_2)(U_i\times_S U_j)
$$
is an equalizer diagram. This means that if $\bar r=\{r_j|j\in J\}$ is an element in
the middle term above, then
\begin{equation}
\label{3.1}
r_i|_{U_i\times_S U_j}\equiv r_j|_{U_i\times_S U_j}
\end{equation}
for all $i,j\in J$ if and only if there is an $r\in\Isom_S(\xi_1,\xi_2)(U)$ so that
its restriction to $U_i$ is $r_i$.
By definition, each $r_i$ consists of two isomorphisms
$r_{i,1}$ and $r_{i,2}$ and a
commutative diagram
$$
\begin{CD}
\cX_1|_{U_i} @>{f_1|_{U_i}}>>  \cW_1|_{U_i}\\
@V{r_{i,1}}VV @V{r_{i,2}}VV\\
\cX_2|_{U_i} @>{f_2|_{U_i}}>> \cW_2|_{U_i}.
\end{CD}
$$
The condition (\ref{3.1}) makes $\{r_{i,1}|i\in J\}$ a descent data for an isomorphism
between $\cX_1$ and $\cX_2$ and makes $\{r_{i,2}|i\in J\}$ a descent data for an
isomorphism between $\cW_1$ and $\cW_2$. Since $\coprod U_i\to U$ is an \'etale covering,
by Descent Lemma there is an isomorphism $r_1\mh \cX_1\to\cX_2$ and $r_2\mh\cW_1\to\cW_2$
restricting to $r_{i,1}$ and $r_{i,2}$ respectively. Hence $r=(r_1,r_2)\in
\Isom_S(\xi_1,\xi_2)(U)$ is the element so that $e(r)=\bar r$. This proves (i).

The proof of (ii) is similar. It amounts to construct a family
$\cX$ over $S$, a family $\cW$ in $\W(S)$ and an $S$-morphism $f\mh\cX\to\cW$,
which form an element $\xi\in\mwcg(S)$, so that $\xi_i=\xi|_{S_i}$ and the identity map
is $\varphi_{ij}$. This is again a standard application of Descent Lemma and will be omitted.
\end{proof}

\begin{lemm}
Let $S$ be any $C$-scheme and let $\xi_1,\xi_2\in\mwcg(S)$ be any two families over $S$.
Then the functor $\Isom_S(\xi_1,\xi_2)$ is represented by a scheme quasi-projective over
$S$.
\end{lemm}

\begin{proof}
We let $\xi_i$ be represented by $f_i\mh \cX_i\to\cW_i$. Let $H_{[f_i,\cW_i]}$ and $H_{[f_i,\cX_i]}$
be the standard choices of relative ample line bundles on $\cW_i$ and on $\cX_i$.
It follows from the uniqueness of these line bundles that whenever
$T\to S$ is a morphism and $(r_1,r_2)$ is a pair of isomorphisms shown below that makes the
diagram
$$
\begin{CD}
\cX_1\times_S T @>{r_1}>> \cX_2\times_S T\\
@V{f_1}VV @V{f_2}VV\\
\cW_1\times_ST @>{r_2}>>  \cW_2\times_S T.
\end{CD}
$$
commutative then
$$r_1\sta H_{[f_2,\cX_2]}\cong H_{[f_1,\cX_1]}\quad
{\rm and}\quad
r_2\sta H_{[f_2,\cW_2]}\cong H_{[f_1,\cW_1]}.
$$
Hence $\Isom_S(\xi_1,\xi_2)$ is isomorphic to the functor associating to $T\to S$ the set of pairs of
isomorphisms of polarized projective schemes over $T$ satisfying the obvious compatibility
condition. It is a routine application of Grothendieck's results on the representability of the
Hilbert scheme and related functor \cite{Gro}, see also \cite[p84]{DM}, that the functor
$\Isom_S(\xi_1,\xi_2)$ is represented by a scheme quasi-projective over $S$.
This completes the proof of the Lemma.
\end{proof}

In the following, we will denote the scheme representing the functor $\Isom_S(\xi_1,\xi_2)$
by $\isom_S(\xi_1,\xi_2)$.

\begin{lemm}
The scheme $\isom_S(\xi_1,\xi_2)$ is quasi-finite and unramified over $S$.
\end{lemm}

\begin{proof}
To prove that $\isom_S(\xi_1,\xi_2)$ is unramified over $S$, we only need to check the
case where $S$ is itself a closed point. Namely $S=\spec\kk$, in which case $\isom_S(\xi_1,\xi_2)$
is either empty or is isomorphic to $\Autw(\xi_1)$. Since ${\rm char}\,\kk=0$, it is well-known that
$\Autw(\xi_1)$ is unramified over $\kk$ if there are no vector fields $v_1$ of $\cX_1$ and
vector fields $v_2$ of $\cW_1$ (whose push-forward to $W$ is trivial) so that they are compatible
via $f_1\mh\cX_1\to\cW_1$. Such vector field does not exist because $\xi_1$ is stable. Hence
$\isom_S(\xi_1,\xi_2)$ is unramified over $S$. Further, since $\isom_S(\xi_1,\xi_2)$ is of finite type over
$S$, it is quasi-finite over $S$. This proves the Lemma.
\end{proof}

\begin{lemm}
\label{3.45}
Let $S=\spec R$ be a $C$-scheme where $R$ is a discrete valuation domain with residue field $\kk$ and
quotient field $K$. We denote by $\eta_0$ and $\eta$ the closed and the generic point of $S$.
Let $\xi_1,\xi_2\in\mwcg(S)$ be two families over $S$ so that there is an arrow $r_K$ from $\xi_1\times_S\eta$
to $\xi_2\times_S\eta$. Then possibly after a base change
$r_K$ extends to an arrow $r_R$ from $\xi_1$ to $\xi_2$.
\end{lemm}

\begin{proof}
Without lose of generality we can assume that $\xi_i$ is represented by $f_i\mh\cX_i\to W[n_i]$
with $\iota_i\mh S\to C[n_i]$ the associated morphism so that $\iota_i(\eta_0)=\bz\in C[n_i]$.
We will prove the case where $S\to C$ does not factor through $0\in C$.
The other case is similar and will be left to readers.

We let $\eta$ and $\eta_0$ be the generic and the closed point of $S$, as stated in the Lemma.
Let $\cX_K=\cX_1\times_S\eta$, which is isomorphic to $\cX_2\times_S\eta$ by assumption. Let
$\bar f_K\mh \cX_K\to W$ be the composite $\cX_1\times_S\eta\to W[n_1]\to W$.
Since $S\to C$ does not factor through $0\in C$, $\bar f_K$ is an ordinary stable morphism.
By the property of stable morphisms, possibly after a base change $\tilde S\to S$
the morphism $\bar f_K\times_S\tilde S$ extends to a family of stable morphisms $\bar f\mh\cX\to W$
over $\tilde S$. For simplicity by replacing $S$ with $\tilde S$ we can assume $\tilde S=S$.
We consider the induced morphism
$$p_i\circ f_i: \cX_i\lra W,
$$
where $p_i\mh W[n_i]\to W$ is the canonical projection. By our assumption, the restriction of
$p_i\circ f_i$ to $\cX_i\times_S\eta$ is isomorphic to $\bar f_K$. Therefore by the property of stable
morphisms $\bar f\mh\cX\to W$ is the stabilization of $p_i\circ f_i\mh\cX_i\to W$ for $i=1,2$.
In particular, there is a unique contraction morphism $q_i\mh \cX_i\to\cX$
so that
$$p_i\circ f_i=\bar f\circ q_i:\cX_i\lra W.
$$

Now let $X_i=\cX_i\times_S\eta_0$ and let $A\sub X_1$ be any irreducible component contracted under $q_1$.
Following the same argument in the proof of Lemma \ref{3.90},
one can show that $\deg(f_1\sta H,A)=0$, $f_1(A)$ is contained in
$\Delta_l\sub W[n_1]_0$ for some $2\leq l\leq n_1$ and $f_1(A)$ is contained in a fiber of $\Delta_l\to D$.
There are two cases we need to distinguish. One is when $f_1(A)$ is a point.
Then that $A$ is contracted by $q_1$ implies $A$ contains at least two nodal points of $X_1$.
On the other hand, when $f_1(A)$ is not a point, then $f_1(A)$ is a fiber of
$\Delta\to D$ and hence
because $f_1\times_S\eta_0$ is stable and hence pre-deformable, the component
$A$ contains at least two nodal points of $X_1$ as well. Because this is
true for every irreducible components in $X_1$ contracted under $q_1$, by the property of stable
contraction, an irreducible component $A\pri$ of $X_1$ is contracted under $q_1$ if and only if
it is isomorphic to $\Po$, it contains exactly two nodal points but no marked points and
$\deg(f_1\sta H, A\pri)=0$. Then since $f_1\times_S\eta_0$ is stable, such $A\pri$ must be a trivial
component of $X_1$. Conversely, any trivial component of $X_1$ obviously will be
contracted under $q_1$. This shows that the morphism $q_i\mh \cX_i\to\cX$ contracts exactly
all trivial components of $X_i$. Consequently, if we let $E_i\sub X_i$ be the union of all
trivial components in $X_i$, then the contractions $q_1$ and $q_2$ define an isomorphism
\begin{equation}
\label{3.25}
\begin{CD}
\varphi: \cX_1-E_1 @>{\cong}>> \cX_2-E_2
\end{CD}
\end{equation}
extending the isomorphism $\cX_1\times_S\eta\cong\cX_2\times_S\eta$.

To proceed, we assume $n_1\geq n_2$ without lose of generality. Note that $W[n_1]$ is
a resolution of
$W[n_2]\times_{\bA^{\!n_2+1}}\bA^{\!n_1+1}$
\footnote{$\bA^{\!n_1+1}\to\bA^{\!n_2+1}$ is defined by
$(t_1,\cdots,t_{n_1+1})\mapsto (t_1,\cdots,t_{n_2},\Pi_{j=n_2+1}^{n_1+1}t_j)$.},
which induces a morphism $q\mh W[n_1]\to W[n_2]$. We consider
$$q\circ f_1: \cX_1\lra W[n_2].
$$
Clearly no irreducible components of $X_1$ are mapped entirely to the nodal divisor of
$W[n_2]_0$ unless $n_1>n_2$ and in which case some non-trivial components of $X_1$
will be mapped to and only to the $n_2$-th nodal divisor of $W[n_2]_0$.
Let $\tilde\iota_1\mh S\to C[n_2]$ be the morphism induced by $q\circ f_1$.
Because $S\to C$ does not factor
through $0\in C$, there are isomorphisms $r_{1,\eta}$ and $r_{2,\eta}$ fitting into the commutative diagram
\begin{equation}\label{3.70}
\begin{CD}
\cX_1\times_S\eta @>{q\circ f_1}>> W[n_2]\times_{C[n_2]}\tilde\iota_1(\eta)\\
@VV{r_{1,\eta}}V @VV{r_{2,\eta}}V\\
\cX_2\times_S\eta @>{f_2}>> W[n_2]\times_{C[n_2]}\iota_2(\eta)
\end{CD}
\end{equation}
By Lemma \ref{1.167} and Corollary \ref{1.93}, $r_{2,\eta}$ is induced by the
$G[n_2]$-action on $W[n_2]$ via a morphism $\rho\mh\spec K\to G[n_2]$.
Let $v$ be the uniformizing parameter of $R$ and let $\rho$ be defined by
$$\rho\sta(\sigma_1)=c_1v^{a_1},\cdots,\rho\sta(\sigma_{n_2})=c_{n_2}v^{a_{n_2}},
$$
where all $c_i$ are units in $R$ and $a_i$ are integers. Here we follow the convention
$$G[n_2]=(G_m)^{\times n_2}
\quad {\rm and}\quad
\Gamma(G[n_2])=\kk[\sigma_1,\sigma_1\upmo,\cdots,\sigma_{n_2},\sigma_{n_2}\upmo].
$$
Now let $A$ be an irreducible component of $X_1$ whose image under $f_1$ is contained in the first
irreducible component of $W[n_1]_0$. Then $p_1\circ f_1(A)$ is not contained in $D\sub W$,
$A$ is a non-trivial component and $f_2\circ\varphi(A)$ is also contained in the
first component of $W[n_2]_0$.
Now let $A$ be a non-trivial irreducible component of $X_1$ whose image under $f_1$ is contained
in the second irreducible component of $W[n_1]_0$. (Here we assume $n_1\geq 1$ since otherwise there
is nothing to prove.) We let
$h\mh(\tilde\eta_0, T)\to \cX_1-E_1$, where $E_1$ is defined in (\ref{3.25}),
be a curve so that $h(\tilde\eta_0)$ is a closed point in $A$ in general position and
the composite $T\mapto{h}\cX_1\to S$ is a branched covering ramified at $\tilde\eta_0$.
We pick an affine open $V\sub W$
so that the image of $h(\tilde\eta_0)$ under $p_1\circ f_1$ is contained in $V \cap D$.
Let $\hat T$ be the formal completion of $T$ along $\tilde\eta_0$, let $\hat V$ be the formal
completion of $V$ along $D\cap V$ and let $\hat h\mh\hat T\to
\cX_1$ be the induced morphism.
Then the morphism
$$q\circ f_1\circ \hat h: \hat T\to W[n_2]
$$
factor through $W[n_2]\times_W\hat V\to W[n_2]$, which we denote by
$$\hat h_1: \hat T\lra W[n_2]\times_W\hat V\cong (V\cap D)\times\hat \Gamma[n_2].
$$
Here $\hat\Gamma[n_2]$ is the scheme defined before (\ref{1.67}). We let
$$\hat h_2:\hat T\lra W[n_2]\times_W\hat V\cong (V\cap D)\times\hat \Gamma[n_2]
$$
be the morphism induced by $f_2\circ\varphi\circ h \mh T\to W[n_2]$.
Because of (\ref{3.70}), if we let $\pr_2$ be the second projection of $(V\cap D)\times\hat\Gamma[n_2]$ and
let
$$\Phi_{\Gamma}\mh G[n_2]\times\hat\Gamma[n_2]\to\hat\Gamma[n_2]
$$
be the group action morphism, then
$$(\pr_2\circ\hat h_1)^{\rho}=
\pr_2\circ\hat h_2:\hat T\lra\hat\Gamma[n_2].
$$
To utilize this identity, we use the open subset $U_1$ of $\Gamma[n_2]$ constructed in Lemma
\ref{1.63}. Since $\pr_2\circ \hat h_1(\tilde\eta_0)$
is contained in the second irreducible component of $\hat\Gamma[n_2]$,
if we let $(u_i^1)$ be the coordinate chart of $U_1$
introduced in Lemma \ref{1.63}, then
$$(\pr_2\circ \hat h_1)\sta(u_j^1)=\alpha_j\in\cO_{\hat T},\quad
j=1,\cdots,n_2.
$$
Because $\pr_2\circ \hat h_1(\tilde\eta_0)$ is in the smooth
locus of $\Gamma[n_2]_0$, all $\alpha_j$ are units in $\cO_{\hat T}$
except $\alpha_1$, which is in the maximal ideal $\mm_{\hat T}$ of $\cO_{\hat T}$.
On the other hand, because $T\to S$ is flat, there is a $0\ne c\in\mm_{\hat T}$ so that
$R\to \cO_{\hat T}$ is given by $v\mapsto c$.
Therefore using the formula in Lemma \ref{1.63},
\begin{equation}
\label{3.71}
\bl (\pr_2\circ\hat h_1)^{\rho}\br\sta(u_i^1)=
(\alpha_1,c_1c^{a_1}\alpha_2, w_3,\cdots,w_{n_2+2})
\end{equation}
for some $w_3,\cdots,w_{n_2+2}$ in the quotient field of $\cO_{\hat T}$. Hence
if we let $\pi\mh \Gamma(n_2)\to\bA^{\!n_2+2}$ be the projection, then
$((\pi\circ(\pr_2\circ\hat h_1)^\rho)\sta$ maps
$$t_1\mapsto c_1c^{a_1}\alpha_1\alpha_2, \ t_2\mapsto w_3,\cdots, t_{n_2+1}\mapsto w_{n_2+2}.
$$
Because this composition specializes to $0_{\bA^{\!n_2+1}}\in \bA^{\!n_2+1}$ by assumption, all $w_i$ are regular and
are in the maximal ideal $\mm_{\hat T}$.
Furthermore, we know $\pr_2\circ\hat h_2$ does not specialize to the first nodal
point of $\Gamma[n_2]_0$, therefore $c_1c^{a_1}\alpha_2$ is not in $\mm_{\hat T}$.
Since $\alpha_2$ is a unit in $\cO_{\hat T}$, $c_1$ is a unit in $R$ and $c$ is in the maximal ideal $\mm_{\hat T}\sub\cO_{\hat T}$,
this implies $a_1\leq 0$.
Similarly, we can pick a non-trivial irreducible component $A$ of $X_2$ so that its image under
$f_2$ is contained in the second component of $W[n_2]_0$. Then because $f_1\circ \varphi\upmo(A)$ is
not contained in the first nodal divisor of $W[n_2]_0$, a parallel argument shows that $a_1\geq 0$.
Therefore, $a_1=0$.

Repeating this argument, using the fact that no irreducible components of $X_1$ are mapped entirely
to the first $n_2$ nodal divisors of $W[n_2]_0$ under $q\circ f_1$, we can show that all
$a_1,\cdots,a_{n_2}$ are
zero. Therefore, $\rho\mh\spec K\to G[n_2]$ extends to $\rho_R\mh S\to G[n_2]$.

Now we prove that the arrow $r_K$ extends to an arrow $r_R$
from $\xi_1$ to $\xi_2$.
Let $\tilde f_1\mh\tilde\cX_1\to W[n_2]$ be the stable contraction of $q\circ f_1\mh\cX_1\to W[n_2]$
as a morphism to $W[n_2]$. Then $\tilde\cX_1\times_S\eta$ is $\cX_1\times_S\eta$, which is isomorphic to
$\cX_2\times_S\eta$ via the arrow $r_K$. We consider the morphism
$$(\tilde f_1)^{\rho_R}:\tilde\cX_1\lra W[n_2].
$$
By assumption, its restriction to $\tilde\cX_1\times_S\eta$ is identical to the restriction of $f_2$ to
$\cX_2\times_S\eta$ via $\tilde\cX_1\times_S\eta\cong\cX_2\times_S\eta$. On the other hand, both
$(\tilde f_1)^{\rho_R}$ and $f_2$ are stable extensions of a stable morphism from $\cX_1\times_S\eta
\cong \cX_2\times_S\eta$ to $W[n_2]$. Hence by the uniqueness of the extensions of ordinary stable morphisms,
the isomorphism $\cX_1\times_S\eta\cong\cX_2\times_S\eta$ extends to an isomorphism $\tilde\cX_1\cong\cX_2$
and the morphism $(\tilde f_1)^{\rho_R}$ is identical to $f_2$. This shows that
$(\tilde f_1)^{\rho_R}$ maps no irreducible components of $\tilde\cX_1\times_S\eta_0$ to the nodal
divisors of $W[n_2]_0$.  We claim that $n_1$ must be equal to $n_2$.
Otherwise, all irreducible components of $X_1$ that are mapped to $\Delta_{n_2+1}\sub W[n_1]_0$ (under $f_1$)
are trivial components, contradicting to Lemma \ref{3.90}.
Once we have $n_1=n_2$, then
 $\tilde f_1=f$ and hence $\tilde\cX_1=\cX_1$ and $f_2=
f_1^{\rho_R}$ follows immediately. This defines an arrow $r_R$ that is an extension of $r_K$. This proves the Lemma.
\end{proof}

\begin{lemm}
\label{3.95}
Let $S=\spec R$ be a $C$-scheme where $R$ is a discrete valuation domain as before and let
$\eta$ be the generic point of $S$. Let
$\xi_K\in\mwcg(\eta)$ be a family of stable morphisms over $\spec K$. Then
possibly after a base change $\tilde S\to S$, $\xi_K\times_S\tilde S$ extends
to a family $\xi\in\mwcg(\tilde S)$.
\end{lemm}

\begin{proof}
Let $\xi_K$ be represented by $f_K\mh \cX_K\to\wn$. We will prove
the case where $\spec K\to C$ does not factor though $0\in C$
and leave the other case to the readers.

Since $\spec K\to C$ does not factor through $0\in C$, we can assume without loss of generality
that $f_K$ is represented by $f_K\mh\cX_K\to W$. Then $f_K$ is an ordinary stable morphism over $\spec K$.
By the property of stable morphism, possibly after a base change $\tilde S\to S$,
the morphism $f_K\times_S\tilde S$ extends to an $\tilde S$-family of ordinary stable morphism
$f_1\mh \cX_1\to W$. Again by replacing $S$ with $\tilde S$ we can assume $S=\tilde S$.
Let $\eta_0$ and $\eta$ be the closed and the generic points of $S$. In case $S\to C$
sends $\eta_0$ to $C-0$, then $f_1$ is already a family of stable morphisms in $\mwcg(S)$. Now
assume $\eta_0$ is mapped to $0\in C$.
Let $\iota\mh S\to\cn$ be any $C$-morphism. Then since $\wn\times_{\cn}\eta$ is canonically isomorphic to
$W\times_C\eta$, $f_K\mh\cX_K\to W$ induces a morphism $\tilde f_K\mh\cX_K\to\wn$.
We say $f_K$ admits a partial extension to $\wn$ if after a base change $\tilde S\to S$
the morphism $\tilde f_K\times_S\tilde S$ extends to
a family of quasi-stable morphisms over $\tilde S$
\begin{equation}
\label{3.80}
f_n: \cX_n\lra \wn
\end{equation}
so that the associated $\iota\mh \tilde S\to\cn$ maps the closed point of $\tilde S$ to $\bzero$.
Here we say $f_n$ is quasi-stable if $\Autw(f_{n,\eta_0})$ is finite, where the automorphism group
is the set of pairs $(a,b)$ such that $a$ is an automorphism of the domain of $f_{n,\eta_0}$ and
$b\in\gn$ such that $f_n\circ a=(f_{n,\eta_0})^b$.
(As before, we let $f_{n,\eta_0}$ be the restriction of $f_n$ to the fiber $\cX_{n,\eta_0}
=\cX_n\times_S\eta_0$.)
It is clear that the extension $f_1$ is a
quasi-stable extension to $W=W[1]$. On the other hand, by
the proof of Lemma \ref{3.90} any quasi-stable extension $f_n$ satisfies
$n\leq \vert \delta\vert$.

Now we let $n$ be the largest possible integer so that there is a quasi-stable extension $f_n$. We show that
$f_n$ defines an extension of $f_K$.
Let $f_n$ be a quasi-stable extension as in (\ref{3.80}) via $\iota\mh S\to\cn$. If
$$f_{n,\eta_0}: \cX_{n,\eta_0}\lra \wn_0
$$
is non-degenerate, then Lemma \ref{2.4} implies that $f_n$
is pre-deformable. Hence $f_n\in\mwcg(S)$ as desired. Now assume $f_{n,\eta_0}$ is
degenerate. Let $A$ be an irreducible component of $\cX_{n,\eta_0}$ so that $f_{n,\eta_0}(A)$ is
contained, say, in the $l$-th nodal divisor of $\wn_0$. We distinguish the case where $l\leq n$ from the
case where $l=n+1$. We first consider the case $l\leq n$. We let $v$ be a uniformizing parameter of $R$ and let
$\rho\mh\spec K\to G_m$ be defined by $\rho\sta(\sigma)=v^{\alpha}$ for an integer $\alpha$.
Let $\lambda_l\mh G_m\cong\gn_l\to\gn$ be the $l$-th one-parameter subgroup of $\gn$. We consider
\begin{equation}
\label{3.88}
(f_{n,\eta})^{{\lambda_l\circ\rho}}: \cX_{n,\eta}\lra\wn.
\end{equation}
As in the proof of Lemma \ref{3.45}, we can pick a curve
$h\mh(\tilde\eta_0,T)\to\cX_n$ covering $\underline{h}\mh (\tilde\eta_0,T)\to (0,C)$ so that
$\underline{h}$ is flat,
$h(\tilde\eta_0)\in A$ and is in general position of $A$. Further we can find an integer
$\alpha$ so that
$$(f_n\circ h)^{{\lambda_l\circ\rho\circ\uline{h}}}:T \lra\wn
$$
specializes to a point in the $(l+1)$-th irreducible component of $\wn_0$ and is away from the nodal
divisors of $\wn_0$. A simple analysis using the covering $U_l$ of $\Gamma[n]$ shows that
if we let
$$\tilde f_n:\tilde\cX_n\to \wn
$$
be the extension of (\ref{3.88}) to a family of stable morphisms, then
$\tilde f_n$ is still quasi-stable. Further, some irreducible component of
$\tilde\cX_{n,\eta_0}$ is mapped entirely to the $(l+1)$-th nodal divisor of
$\wn_0$. This shows that we can assume without loss of generality that $f_n$ maps
an irreducible component $A$ of $\cX_{n,\eta_0}$ to $D\lnpo\sub\wn_0$.

We now show that there is a quasi-stable extension of $f_K$ to $W[n+1]$.
Let $f_n$ and $A\sub\cX_{n,\eta_0}$ be as before so that $f_n(A)\sub D_n$.
We let $\cn\sub C[n+1]$ and $\wn\to W[n+1]$ be the standard embedding associated
to $[n]\sub [n+1]$. Namely, they are induced by the embedding $\Anpo\to\Anpt$ that keep
the last coordinate 1. Let $\tilde f_n\mh \cX_n\to W[n+1]$ be the
composite of $f_n$ with $\wn\to W[n+1]$. Then the same technique as before shows that
we can find an integer $\alpha$ so that if we let $\rho\mh\spec K\to G_m$ be defined
by $\rho\sta(\sigma)=v^{\alpha}$, then the extension (possibly after a base change $\tilde S\to S$) of
$$(\tilde f_n)^{{\lambda_n\circ\rho}}:\cX_{n,\eta}\lra W[n+1]
$$
to a family of ordinary stable morphisms
$$ f_{n+1}:\cX_{n+1}\lra W[n+1]
$$
is quasi-stable. This contradicts to the maximal assumption of $n$. Hence $f_n$ must
already be a family of non-degenerate morphism. Since $\tilde S\to C$ is flat, it must be pre-deformable and thus
in $\mwcg(\tilde S)$, extending $\xi_K\times_S\tilde S\in\mwcg(\eta)$.
This proves the Lemma.
\end{proof}

\begin{theo}
\label{32.22}
The moduli stack $\mgwc$ of stable morphisms to $\WWc$ of topological type $\Gamma$
is separated and proper over $C$. Furthermore, it is a Deligne-Mumford stack.
\end{theo}

\begin{proof}
The fact that the moduli stack $\mgwc$ is separate over $C$ follows from Lemma \ref{3.45} and that it
is proper over $C$ follows from Lemma \ref{3.95}. It remains to show that it is algebraic. Namely,
$\mgwc$ admits an \'etale cover by a scheme of finite type. We now show that it
admits an \'etale covering by a quasi-projective scheme.

Let $\MM(\wn,\Gamma)$ be the moduli stack of stable morphisms from pointed curves
to $\wn$ of type $\Gamma$, and let $\mgkst$ be the substack of $\MM(\wn,\Gamma)$ consisting of
all pre-deformable morphisms that are stable in the sense of Definition \ref{31.1}. Clearly,
$\mgkst$ is a locally closed substack of $\MM(\wn,\Gamma)$. Since members of $\mgkst$ are stable
morphisms to $\WWc$, there is a natural morphism $\mgkst\to\mwg$. Because of Lemma \ref{3.20},
there is an integer $\bar n$ so that
\begin{equation}
\label{32.1}
\bigcup_{n=1}^{\bar n} \mgkst\lra\mwg
\end{equation}
is surjective.
Hence $\mwg$ is surjected onto by a quasi-projective scheme since each $\mgkst$ does.

We now show that we can find a quasi-projective scheme $Y$ and a surjective \'etale morphism
$Y\to\mwg$.
Let $p$ be any closed point in $\mwg$.
Since (\ref{32.1}) is surjective, $p$ is contained in the image of $\mgkst$ for
some $n\leq \bar n$. Now let $\rho\mh\bar S\to\mgkst$ be a chart\footnote{All charts of
stacks are \'etale charts unless otherwise is specified.}
so that $p$ is contained in the image of $\bar S$ in $\mwg$. To obtain a chart of
$\mwg$ that contains $p$, we need to investigate the $\gn$-action on $\mgkst$.
Since $\wn$ is a $\gn$-variety, $\mgkst$ admits a natural $\gn$-action.
By the definition of stability (cf. Definition \ref{31.1}), the $\gn$ action on $\mgkst$
has only finite stabilizer. Let
$$\Phi: \bar S\times \gn\lra \mgkst
$$
be the morphism induced by the group action. Namely, $\Phi(s,\sigma)=\rho(s)^{\sigma}$.
Then there is an open subset $U\sub \bar S\times\gn$ that contains $\bar S\times\{e\}$
so that $\Phi|_U$ lifts to a morphism $\tilde\Phi\mh U\to\bar S$
so that $\Phi|_U=\rho\circ\tilde\Phi$. Let $s_0\in\bar S$ be a point so that $\rho(s_0)=p$. Since the
stabilizer of $p\in\mgkst$ (of the $\gn$-action) is finite, $\tilde\Phi\mh U\to\bar S$ is
a smooth morphism near $(s_0,e)$. Hence we can find a locally closed subscheme $S\sub\bar S$ containing
$s_0$ so that $\tilde\Phi|_{U\cap(S\times\gn)}$ is \'etale at $(s_0,e)$. By shrinking $S$
if necessary, we can assume $\tilde\Phi|_{U\cap(S\times\gn)}$ is \'etale near $S\times\{e\}$.
Therefore, the composite
$$S\lra \mgkst\lra\mwg
$$
is \'etale. Since $\mwg$ is bounded, we can find a finite number of \'etale charts of $\mwg$,
each quasi-projective, so that their union covers $\mwg$. This proves that $\mwg$ is an algebraic
(i.e. Deligne-Mumford) stack.
\end{proof}


\section{Stacks of relative stable morphisms}

Let $D\sub Z$ be a smooth connected divisor (called the distinguished divisor)
in a smooth projective variety.
In this section, we will construct the moduli stack of relative stable morphisms
with prescribed contact with $D$. We will prove some basic property of this
moduli stack and describe its relation with the
stack of stable morphisms constructed in the previous section.
This construction can easily be generalized to the case where $D$ is smooth
but not necessarily connected.

\subsection{Standard models}

We first construct the standard models $(\zn,\dn)$ which are the building blocks of the
stack of expanded relative pairs of $(Z,D)$. For $n=0$, we let $(Z[0],D[0])=(Z,D)$.
For $n=1$, $Z[1]$ is the blowing-up of $Z\times\Ao$ along $D\times\{0\}\sub Z\times\Ao$.
Its distinguished divisor $D[1]$ is the proper transform of $D\times\Ao\sub Z\times\Ao$
in $Z[1]$. We now let $W\to\Ao$ be $Z[1]\to\Ao$. We shall view the ruled variety
$\Delta=\bP(\bone_D\oplus N_{D/Z})$ and $Z$ in $W_0$($=W\times_{\Ao}0$)
as its top (left) and its bottom (right) components.
With this understanding, we can construct a sequence of varieties
$\wn$ over $\Anpo$ for all $n\geq 1$ as in Section 1 with the
projection $\pi\mh \wn\to Z\times\Anpo$.
We define $\zn=W[n-1]$ and define its distinguished divisor $D[n]$ to be
the proper transform of $D\times \An\sub Z\times\An$ in $Z[n]$.
We call $(\zn,\dn)$ with the tautological projection $\zn\to Z\times \An$
an expanded relative pair of $(Z,D)$.
Note that under this projection $\dn$ is isomorphic to $D\times\An
\sub Z\times\An$.
We denote the projection $\zn\to\An$ by $\pi$.

Clearly, the fiber of $\zn$ over $0\in\An$, denoted $\zn_0$,
is reduced with normal crossing singularity.
It consists of a chain of smooth varieties, one $Z$ and $n$ ruled varieties $\Delta$.
We give these $n$ components an ordering so that
the first component of $Z[n]_0$ contains the distinguished divisor $D[n]_0$ and
the remainder components are ordered according to their intersection chain structure. Namely, the
$k$-th component will intersect the $(k+1)$-th component.
We call this the standard ordering of $Z[n]_0$.
We continue to denote by $\gn$ the product of $n$ copies of $G_m$. Since the $G_m$-action
on $Z\times \Ao$ defined via $(z,t)^{\sigma}=(z,\sigma\upmo t)$
leaves $D\times\{0\}\sub Z\times\Ao$ fixed, it lifts to an action on $Z[1]\to \Ao$.
In general, we let $\gn$ acts on $\An$ via
\begin{equation}
\label{4.13}
(t_1,\cdots,t_{n})^{\sigma}=(\bar\sigma_2t_1,\bar\sigma_3t_2,\cdots,\bar\sigma_{n+1}t_{n}).
\end{equation}
(Recall $\bar\sigma_i=\sigma_i/\sigma_{i-1}$, see (\ref{1.126}).)
Note that if we embed $\An$ in $\Anpo$ in the standard way as the $n+1$-th coordinate hyperplane
$\bH_{n+1}\unpo$ and let $\gn$ acts on $\Anpo$ as defined in Section 1,
then the embedding $\An\sub\Anpo$ is $\gn$-equivariant.
It is direct to check that there is a unique $\gn$-action on $\zn$ that makes $\zn\to Z\times\Anmo$
$\gn$-equivariant, where the $\gn$ action on $Z\times\Anmo$ is $(z,t)^{\sigma}=(z,t^\sigma)$
defined in (\ref{4.13}).

For later application, we need to consider $\zn$ with the reversed ordering of its components.
In the following for any $n$ we define $\bone_{\An}^o\mh\An\to\An$ be the morphism
$(t_1,\cdots,t_n)\mapsto(t_n,\cdots,t_1)$.
We let $\zn^o$ be $\zn$ with the tautological projection
$\pi^o\equiv\bone_{\An}^o\circ\pi\mh \zn\to\An$.
At the same time we reverse the order of the components of
$\zn_0$ so that the first irreducible component of $\zn^o_0$ is $Z$ while its
last component contains the distinguished divisor. Note that the restriction of $\zn_0^o$
to the $l$-th coordinate line is still a smoothing of the $l$-th nodal divisor of $\zn_0^o$.
Accordingly, we let the $\gn$-action on $\An$ be
\begin{equation}
\lab{4.1.2}
(t_1,\cdots,t_{n})^{\sigma^o}=(\bar\sigma_1 t_1,\cdots,\bar\sigma_{n}t_{n}).
\end{equation}

We next describe the connection between the standard model constructed in Section 1 with the relative
standard models so constructed. Let $W\to C$ be the pair in Section 1 with $D_1\sub Y_1$ and $D_2\sub Y_2$
the two irreducible components of $W_0$ with their distinguished divisors. Using the pair
$(Y_i,D_i)$ we can construct the relative pair $(Y_i[n],D_i[n])$ with the associated projection
$Y_i[n]\to Y_i\times\An$. Now let $l$ be any integer in $[n]$. We consider the pairs
$$
D_1[l]\times\Anl\sub Y_1[l]^o\times\Anl\quad{\rm and}\quad
\Al\times D_2[n-l]\sub\Al\times Y_2[n-l].
$$
We have the canonical isomorphisms
$$D_1[l]\times\Anl\mapright{\cong}(D_1\times\bA^{\!l})\times\Anl\mapright{\cong}
D_1\times\bA^{\!l}\times\Anl
$$
and canonical isomorphisms
$$\Al\times D_2[n-l]\mapright{\cong} \Al\times(D_2\times \Anl)\mapright{\cong}
D_2\times\Al\times\Anl,
$$
where the second arrow is induced by exchanging the factor $\Almo$ and $D_2$.
Combined with the canonical isomorphism $D_1\cong D_2$, we have
\begin{equation}
\label{4.2}
D_1[l]\times\Anl\cong \Al\times D_2[n-l].
\end{equation}
Hence we can glue $Y_1[l]^o\times \bA^{\!n-l}$ with $\bA^{\!l}\times Y_2[n-l]$ along
$$ D_1[l]\times\Anl\sub Y_1[l]^o\times \Anl
\quad{\rm and}\quad
\Al\times D_2[n-l]\sub \Al\times Y_2[n-l]
$$
according to the isomorphism (\ref{4.2}) to obtain a new scheme, denoted by
\begin{equation}
\label{4.3}
Y_1[l]^o\times\Anl\adj \Al\times Y_2[n-l].
\end{equation}

The group $\gn$ acts on (\ref{4.3}) as follows:
We let $\Al\times\Anl\to\Anpo$ be the embedding as coordinate plane defined by
$(r,s)\mapsto (r,0,s)$.
Then the $\gn$-action on $\Anpo$ (see (\ref{1.127})) lifts to an action on
$\Al\times\Anl$. Namely for $\sigma=(\sigma_1,\cdots,\sigma_{n-1})\in\gn$,
$$(r_1,\cdots,r_{l})^{\sigma}=(\bar\sigma_1 r_1,\cdots,\bar\sigma_{l}r_{l})\
{\rm and}\
(s_1,\cdots,s_{n-l})^{\sigma}=(\bar\sigma_{l+2}s_1,\cdots,\bar\sigma_{n+1} s_{n-l}).
$$
Let $\rho_-\mh\gn\to G[l]$ and $\rho_+\mh\gn\to G[n-l]$ be
the homomorphism
$$\rho_-(\sigma)=(\sigma_1,\cdots,\sigma_{l})
\quad {\rm and}\quad
\rho_+(\sigma)=(\sigma_{l+1},\cdots,\sigma_{n}).
$$
Thus $r^{\sigma}=r^{\rho_-(\sigma)^o}$ and $s^{\sigma}=s^{\rho_+(\sigma)}$,
where $r^\sigma$ and $s^\sigma$ are the actions
defined above and $r^{\rho_-(\sigma)^o}$ and $s^{\rho_+(\sigma)}$ are the actions
defined in (\ref{4.13}) and (\ref{4.1.2}).
We define the $\gn$-action on $Y_1[l]^o\times\Anl$ and on $\Al\times Y_2[n-l]$ be
$$(x,s)^{\sigma}=(x^{\rho_-(\sigma)^o},s^{\rho_+(\sigma)})
\quad{\rm and}\quad
(r,y)^{\sigma}=(r^{\rho_-(\sigma)^o},y^{\rho_+(\sigma)}).
$$
It is direct to check that the isomorphism (\ref{4.2}) is $\gn$-equivariant,
and hence (\ref{4.3}) is a $\gn$-scheme.

\begin{prop}
\label{4.8}
Let $\pi\mh \wn\to\Anpo$ be the associated morphism and $\pi_l\mh\wn\to\Ao$
be the composite of $\pi$ with the $l$-th projection of $\Anpo$. Then we have a canonical
isomorphism
$$ Y_1[l]^o\times\Anl\adj \Al\times Y_2[n-l]\cong \wn\times_{\pi_l} 0_{\Ao}\sub\wn.
\footnote{By $\wn\times_{\pi_l}0_{\Ao}$ we mean $\wn\times_{\Ao}0_\Ao$ with $\wn\to\Ao$
given by $\pi_l$.}
$$
Further the above inclusion and isomorphism are $\gn$-equivariant.
\end{prop}

\begin{proof}
The proof is straightforward and will be omitted.
\end{proof}

Now let $(\cY_i,\cD_i)$, $i=1$ and $2$, be effective relative pairs over $S$ associated to
$\tau_i\mh S\to\bA^{\!l_i}$ defined by
$$(\cY_i,\cD_i)=(Y_i[l_i],D_i[l_i])\times_{\bA^{\!l_i}}S.
$$
As before, we define $\tau_1^o=\bone_{\bA^{\!l_1}}^o\circ \tau_1\mh S\to \bA^{\!l_1}$ and
let $\cY_1^o$ be $Y_1[l_1]^o\times_{\tau_1^o}S$.
Obviously $\cY_1^o=\cY_1$ with only the ordering of the components of its fibers reversed.
Since $\cD_1$ (resp $\cD_2$) is canonically isomorphic to
$D_1\times S$ (resp. $D_2\times S$), we can construct a new scheme over $S$ by
gluing $\cY_1^o$ and $\cY_2$ along $\cD_1$ and $\cD_2$ using $\cD_1\cong \cD_2$.
Since both $\cD_1$ and $\cD_2$ are smooth over $S$, the resulting scheme
$\cY_1^o\adj\cY_2$ is a flat family of schemes with
normal crossing singularity relative to $S$.
Let $\An\cong\bA^{\!l_1}\times\Ao\times\bA^{\!l_2}$, $n=l_1+l_2+1$,
be the canonical isomorphism keeping the order of each
copies of $\Ao$ in $\Anpo$. Then we have an embedding $\bA^{\!l_1}\times\bA^{\!l_2}
\to \Anpo$ sending $(r,s)\in \bA^{\!l_1}\times\bA^{\!l_2}$ to $(r,0,s)\in\Anpo$. Let
$\tau\mh S\to\Anpo$ be the composite
\begin{equation}
\label{51.69}
\begin{CD}
\tau: S @>{(\tau_1^o,\tau_2)}>>\bA^{\!l_1}\times\bA^{\!l_2} @>>> \Anpo.
\end{CD}
\end{equation}

\begin{coro}
\label{7.1}
Let the notation be as above and let $\tau\mh S\to\An$ be defined in (\ref{51.69}).
Then the associated family $\cW=\tau\sta\wn$ is canonically isomorphic to $\cY_1^o\adj\cY_2$.
\end{coro}

\begin{proof}
Clearly, the composite of
$\tau\mh S\to\Anpo$ with $\pr_{l_1}\mh\Anpo\to\Ao$ is trivial. Hence
by Proposition \ref{4.8}
\begin{eqnarray*}
\wn\times_{\Anpo}S&\cong &\bl Y_1[l_1]^o\times \bA^{\!l_2}\adj \bA^{\!l_1}\times Y_2[l_2]\br
\times_{\bA^{\!l_1}\times\bA^{\!l_2}}S\\
&\cong & Y_1[l_1]^o\times_{\bA^{\!l_1}} S
\mathop{\adj}_{D_1[l_1]\times_{\bA^{\!l_1}}S=D_2[l_2]\times_{\bA^{\!l_2}}S}
Y_2[l_2]\times_{\bA^{\!l_2}}S, \\
\end{eqnarray*}
which is isomorphic to $\cY_1^o\adj\cY_2$. This proves the Corollary.
\end{proof}

We now define the stack of expanded relative pairs of $(Z,D)$. Let $S$ be any scheme. An {\sl effective}
family of expanded {\sl relative pair} of $(Z,D)$
(in short relative pair) is the associated family $(\cZ,\cD)$ of a morphism $\tau\mh S\to\An$
for some $n$ defined by
$$\cZ=\tau\sta\zn=\zn\times_{\An}S\quad{\rm and}\quad
\cD=\tau\sta \dn=\dn\times_{\An}S.
$$
We call $\cD$ the
distinguished divisor of $\cZ$. Now let
$$\underline{\Phi}:\gn\times\An\lra \An\quad{\rm and}\quad
\Phi:\gn\times\zn\lra\zn
$$
be the $\gn$-action on $\An$ and on $\zn$ defined earlier in this section. Then for any $\tau\mh S\to\An$
and $\rho\mh S\to\gn$, the group action $\uline{\Phi}$ defines a new
morphism $\tau^{\rho}\mh S\to\An$ and thus a new family
of relative pair $({\tau^{\rho}}\sta\zn,{\tau^{\rho}}\sta\dn)$.
Using the group action of $\gn$ on $\zn$, there is a canonical isomorphism between
the pair $(\tau\sta\zn, \tau\sta\dn)$ and $({\tau^{\rho}}\sta\zn,{\tau^{\rho}}\sta\dn)$.
As in the case for $\wn$, given two effective relative pairs $\xi_1=(\cZ_1,\cD_1)$ and
$\xi_2=(\cZ_2,\cD_2)$ over $S$ associated to morphisms $\tau_1\mh S\to\bA^{\!n_1}$
and $\tau_2\mh S\to\bA^{\!n_2}$ respectively, an effective arrow from $\xi_1$ to $\xi_2$ consists of
a standard embedding $\iota\mh \bA^{\!n_1}\to\bA^{\!n_2}$ (associated to an increasing
map $[n_1]\to [n_2]$, see discussion at the beginning of subsection 1.1) and a morphism
$\rho\mh S\to G[n_2]$ such that $(\iota\circ\tau_1)^{\rho}=\tau_2$.
Similar to the discussion before Lemma \ref{1.167}, this identity defines a canonical $S$-isomorphism
of pairs
$$(\tau_1\sta Z[n_1],\tau_1\sta D[n_1])\cong (\tau_2\sta Z[n_2], \tau_2\sta D[n_2]),
$$
compatible to their projections to $Z\times S$.
Along the same line, we can define an effective arrow between $\xi_1$ and $\xi_2$
to be either an effective arrow from $\xi_1$ to $\xi_2$ or an effective arrow from $\xi_2$ to $\xi_1$.
We say $\xi_1$ and $\xi_2$ are equivalent via a sequence of effective arrows if there is a
sequence of effective families $\eta_0,\cdots,\eta_m$ so that $\xi_1=\eta_0$, $\xi_2=\eta_m$ and
that $\eta_i$ is equivalent to $\eta_{i+1}$ via an effective arrow.
Note that this sequence of effective arrows induces a canonical isomorphism of
$(\cZ_1,\cD_1)$ and $(\cZ_2,\cD_2)$.
We have a partial inverse to this, similar to Lemma \ref{1.167}.

\begin{lemm}
Let $\xi_i=(\cZ_i,\cD_i)$, $i=1$ and $2$, be two effective relative pairs over $S$
associated to $\tau_i\mh S\to\bA^{\!n_i}$. Suppose there is an isomorphism
$(\cZ_1,\cD_1)\cong(\cZ_2,\cD_2)$ compatible to their projections to $Z\times S$.
Then to each $p\in S$ there is an open neighborhood $S_0$ of $p\in S$
such that over $S_0$ the induced isomorphism $\cZ_1\times_SS_0\cong\cZ_2\times_SS_0$
is induced by a sequence of effective arrows between $\xi_1\times_SS_0$ and $\xi_2\times_SS_0$.
\end{lemm}

\begin{proof}
The proof is similar to that of Lemma \ref{1.167} and will be omitted.
\end{proof}

We define a family of relative pair over $S$ to be a pair $(\cZ,\cD)$, where
$\cZ$ is an $S$-family with an $S$-projection $\cZ\to Z\times S$ and $\cD$ is a Cartier divisor of
$\cZ$, such that there is an open covering $S\lalp$ of $S$ so that each $(\cZ\times_SS\lalp,\cD\times_SS\lalp)$
is isomorphic to an effective relative pair associated to some $\tau\lalp\mh S\lalp\to \bA^{\!n\lalp}$
and the isomorphism $\cZ\times_S S\lalp\cong\tau\lalp\sta Z[n\lalp]$ is compatible to their
projections to $Z\times S\lalp$.
Note that then $\cD$ is isomorphic to $D\times S\sub Z\times S$ under
the projection $\cZ\to Z\times S$.
We call $\cZ\times_S S\lalp\cong Z[n\lalp]\times_{\bA^{\!n\lalp}}S\lalp$
with its associated data a local representative of $(\cZ,\cD)$.
We say two such families $(\cZ_1,\cD_1)$ and $(\cZ_2,\cD_2)$ are isomorphic if there is an $S$-isomorphism
$(\cZ_1,\cD_1)\to(\cZ_2,\cD_2)$ compatible to the projections $\cZ_1\to Z\times S$ and $\cZ_2\to Z\times S$.
Clearly, if $(\cZ,\cD)$ is a relative pair over $S$ and $\rho\mh S\pri\to S$ is a morphism, then
$(\cZ\times_SS\pri,\cD\times_SS\pri)$ coupled with the obvious projection $\cZ\times_S S\pri\to Z\times S\pri$
is a relative pair over $S\pri$, called the pull back relative pair.

\begin{defi}
We define $\Zzd$ to be the category whose objects are families of expanded relative pairs of $(Z,D)$ over $S$ for
some scheme $S$. Let $S_i$, $i=1$ and $2$, be any two schemes and $\xi_i\in\Zzd(S_i)$
be two objects in $\Zzd$. An arrow from $\xi_1$ to $\xi_2$ consists of an arrow (a morphism)
$\rho\mh S_1\to S_2$ and an isomorphism of relative pairs $\xi_1\cong\rho\sta\xi_2$.
We define ${\mathfrak p}\mh\Zzd\to({\rm Sch})$ to
be the functor that send families over $S$ to $S$. The pair $(\Zzd,{\mathfrak p})$ forms a groupoid.
\end{defi}

\begin{prop}
The groupoid $(\Zzd,{\mathfrak p})$ is a stack.
\end{prop}

\begin{proof}
The proof is straightforward and will be omitted.
\end{proof}

In the following, we will call $\Z\urel$ the stack of expanded relative pairs of $(Z,D)$.

Now let $(Y_1,D_1)$ and $(Y_2,D_2)$ be the two relative pairs associated to the family $W/C$
and let $\YY_1\urel$ and $\YY_2\urel$ be the associated stacks of expanded relative pairs
of $(Y_1, D_1)$ and $(Y_2,D_2)$, respectively.
Then the correspondence $(\cY_1,\cY_2)\mapsto\cY_1^o\adj\cY_2$ defined in Corollary \ref{7.1}
induces a canonical morphism of stacks
\begin{equation}
\label{41.15}
\YY_1\urel\times \YY_2\urel\lra \WW\times_C 0\sub\WW.
\end{equation}

\subsection{Relative stable morphisms}

Let $(Z,D)$ be the pair and $\Zzd$ the stack of expanded relative pairs of $(Z,D)$.
In this part we define the notion of
relative stable morphisms to $\Zzd$. We fix an ample line bundle $H$ on $Z$,
called a polarization of $(Z,D)$.
We first introduce the notion of admissible weighted graphs (in short admissible graphs).

In this paper, by a graph $\Gamma$ we mean a finite collection of vertices, edges, legs and roots. Here
an edge is as usual a line segment with both ends attached to vertices of $\Gamma$.
A leg or a root is a line segment with only one end attached to a vertex of $\Gamma$.
We will denote by $V(\Gamma)$ the set of vertices of $\Gamma$.

\begin{defi}
\label{14.8}
An admissible weighted graph $\Gamma$ is a graph without edges coupled
with the following additional data:
\newline
1. An ordered collection of legs, an ordered collection of weighted roots
and two weight functions $g\mh\ver(\Gamma)\to\ZZ_{{}^{\geq 0}}$
and $b\mh V(\Gamma)\to A_2 Z/\sim_{{\rm alg}}$.
\newline
3. The graph is relatively connected in the sense that either $\ver(\Gamma)$ consists of
a single element or each vertex in $\ver(\Gamma)$ has at least one root attached to it.
\end{defi}

In this section, we will reserve the integer $k$ to denote the number of
legs and the integer $r$ to denote the number of roots of $\Gamma$. We will
denote by $\mu_1,\cdots,\mu_r$ the weights of the $r$ roots.

Given two admissible graphs $\Gamma_1$ and $\Gamma_2$, an isomorphism $\Gamma_1\cong\Gamma_2$
is an isomorphism between their respective vertices, legs and roots
(order preserving in the later two cases) so that it preserves
the weights of the roots and the two sets of weights of the vertices.

We now define the relative stable morphisms to $(\zn,\dn)$ of type $\Gamma$.

\begin{defi}
\label{14.6}
Let $S$ be an $\An$-scheme and $\Gamma$ be an admissible graph
with $r$ roots, $k$ legs and $l$ vertices $v_1,\cdots,v_l$.
An $S$-family of relative stable morphisms to $(\zn,\dn)$ of type $\Gamma$
is a quadruple $(f,\cX,q_i,p_j)$ as follows:
\newline
1. $\cX$ is a disjoint union of $\cX_1,\cdots,\cX_l$ such that each $\cX_i$ is
a flat family of pre-stable
curves\footnote{Recall by pre-stable we mean that $\cX/S$ is a family of
connected, complete and nodal curves.} over $S$ of arithmetic genus $g(v_i)$.
\newline
2. $q_i\mh S\to \cX$, $i=1,\cdots,r$ and $p_j\mh S\to\cX$, $j=1,\cdots,k$, are disjoint sections
away from the singular locus of the fibers of $\cX/S$ so that
$q_i(S)\sub \cX_j$ (resp. $p_i(S)\sub \cX_j$) if the $i$-th root
(resp. $i$-th leg) is attached to the $j$-th vertex of $\Gamma$.
\newline
3. $f\mh\cX\to\zn$ is an $\An$-morphism (i.e. an $S$-morphism $\cX\to\zn\times_{\An}S$)
so that as Cartier divisors $f\upmo(\dn)=\sum_{i=1}^r\mu(x) q_i(S)$,
and that to each closed $s\in S$ the class $\varphi\lsta(f(\cX_s))$ is algebraically equivalent
to $b(v)$.
Here $\varphi\mh \zn\to Z$ be the tautological projection.
\newline
4. Finally, each morphism $f|_{\cX_i}$,
considered as a morphism whose domain is $\cX_i$ with all marked sections in $\cX_i$,
is a family of stable morphisms to $\zn$.
\end{defi}

Note that a necessary condition for the existence of relative stable morphisms of type $\Gamma$
is that to each $v\in V(\Gamma)$,
$\sum_{i\prec v} \mu_i=b(v)\cdot D$,
where $i\prec v$ means that the $i$-th root is attached to the vertex $v$.
We will call $q_i$ the distinguished marked sections and call $p_j$ the ordinary marked sections.

We now fix an admissible graph $\Gamma$.
Let $f\mh\cX\to\zn$ be a relative stable morphism to $(\zn,\dn)$ of type $\Gamma$.
Let $l$ be any integer in $[n]$.
Then $\zn\times_{\An}\bH_l^{n}$
has normal crossing singularity whose singular locus is isomorphic to $D\times\bA^{\!n-1}$.
We denote this nodal divisor by $\bB_l$. As in section 2, we call $f$ non-degenerate if to each
$s\in S$ there are no irreducible components of $\cX_s$ mapped entirely to the union of
$\bB_1,\cdots,\bB_{n}$ under $f$.
We call $f$ pre-deformable if $f$ is pre-deformable along $\bB_1,\cdots,\bB_{n}$ as defined
in section 2.

We now change our view point of $f$. Since $S$ is an $\An$-scheme,
$$(\zn\times_{\An}S,\dn\times_{\An}S)\in Ob(\Zzd(S)),
$$
and hence $f$ can be viewed as an $S$-family
of morphisms to $\Zzd$. We say $f$, considered as an $S$-family of morphisms to $\Zzd$, is
relative pre-stable (resp. pre-deformable; resp. of type $\Gamma$) if $f$,
considered as an $S$-family of morphisms to
$\zn$, is relative stable (resp. pre-deformable; resp. or type $\Gamma$).
As before, we define $\autz(f)$ to be the functor associating to each $S$-scheme $T$ the set of all $(a,b)$,
where $a$ is an automorphism of $\cX\times_S T$ as pointed nodal curve over $T$, namely $a$ leaves
all marked (distinguished and ordinary) sections of $\cX\times_S T$ fixed, and $b$ is
a morphism $T\to\gn$ such that if we let
$f_T$ be the morphism $\cX\times_S T\to \zn$ induced by $f$ then
$$(f_T)^b=f_T\circ a: \cX\times_S T\to \zn.
$$
Here $(f_T)^b$
is the morphism resulting from the $\gn$-action
on $\zn$ associated to $b\mh T\to\gn$ (see (\ref{1.103})).
Obviously, this functor is represented by a group scheme over $S$, called the automorphism
group of $f$ and denoted by $\autz(f)$.

Let $S$ be any scheme and $(\cZ,\cD)\in Ob(\Zzd(S))$ be any object. Then $S$ has an open covering
by $\{S\lalp\}$ so that over each $S\lalp$ we have isomorphism
$$(\cZ,\cD)\times_S S\lalp\cong(Z[n\lalp]\times_{\bA^{\!n\lalp}}S\lalp,D[n\lalp]
\times_{\bA^{\!n\lalp}}S\lalp)
$$
for some $S\lalp\to \bA^{\!n\lalp}$.
For any $S$-family of morphisms $f\mh\cX\to\cZ$, we will call the induced morphisms
$f\lalp:\cX\times_S S\lalp\lra Z[n\lalp]$ local representatives of $f$.

\begin{defi}
\label{14.7}
An $S$-family of relative pre-stable morphisms to $\Zzd$ of type $\Gamma$ is a family $f\mh\cX\to\cZ$,
where $(\cZ,\cD)\in Ob(\Zzd(S))$, so that all its local representative $f\lalp$ are
relative stable (as morphisms to $\cZ\lalp$) of type $\Gamma$.
We say $f$ is stable if in addition to it being pre-stable,
all its local representatives are pre-deformable and that for all closed $s\in S$ the automorphism
groups $\autz(f_s)$ are finite.
\end{defi}

It is routine to check that the definition of $\autz(f_s)$ and the definition
of relative stable morphisms to $\Zzd$ is independent of the
choices of local representatives.

\begin{defi}
\label{14.9}
Let $(Z,D)$ be as before and let $\Gamma$ be an admissible
graph. We define $\Mzdg$ to be the
category whose objects are families of relative stable morphisms to $\Zzd$ of
type $\Gamma$. As usual, $\Mzdg(S)$ is the sub-category of families over $S$.
Let $\xi_1$ and $\xi_2$ be objects in $\Mzdg(S)$ represented by relative stable
morphisms $f_1\mh\cX_1\to\cZ_1$ and $f_2\mh\cX_2\to\cZ_2$.
An arrow $\xi_1\to\xi_2$ in $\Mzdg(S)$
covering $1_S\mh S\to S$ consists of a pair $(\varphi_1,\varphi_2)$, where
$\varphi_1$ is an $S$-isomorphism $\cX_1\to\cX_2$
preserving all ordinary marked points and $\varphi_2$ is an arrow $\cZ_1\to\cZ_2$ in $\Zzd(S)$
covering $1_S$ so that
$f_2\circ\varphi_1=\varphi_2\circ f_1$. In case $\xi_1\in \Mzdg(S)$ and $\xi_2\in\Mzdg(T)$,
then an arrow $\xi_1\to\xi_2$ consists of a morphism $\sigma\mh S\to T$ and an arrow
$\xi_1\to\sigma\sta\xi_2$ in $\Mzdg(S)$ covering $1_S$.
Let $\mathfrak p\mh\Mzdg\to({\rm Sch})$ be the functor that sends families over $S$ to $S$.
The category $\Mzdg$ coupled with the functor $\mathfrak p$ form a groupoid,
called the groupoid of the relative stable morphisms to $\Zzd$ of type $\Gamma$.
\end{defi}

\begin{theo}
\label{14.10}
The groupoid $\Mzdg$ is an algebraic stack. It is separated and proper over $\kk$.
\end{theo}

\begin{proof}
The proof is parallel to that of $\Mwg$ and will be omitted.
\end{proof}

We now move back to the family $W$ defined before.
We fix an ample line bundle $H$ on $W$ as before and let $H_{Y_i}$ be the
restriction of $H$ to $Y_i$.
We let $(Y_1,D_1)$ and $(Y_2,D_2)$ be the
two relative pairs from decomposing $W_0$. We let $\YYo$ and $\YYt$ be the stack of expanded
relative pairs of $(Y_1,D_1)$ and $(Y_2,D_2)$ respectively.
In the remainder part of this section, we will investigate the relation of the stacks $\Mydo$,
$\Mydt$ and $\Mwg$. We begin with a discussion of gluing a pair of relative stable morphisms
$$(f_1,f_2)\in\Mydo\times\Mydt
$$
to form a stable morphism in $\Mwg$. Here as usual the domain of the new morphism will be the
gluing of the
domains of $f_1$ and $f_2$.
In terms of their associated graphs, this amounts to connecting the roots in the graph $\Gamma_{1}$
to the roots in the graph $\Gamma_{2}$ to form a new graph with edges.\footnote{Namely  we identify the free
end of a root with the free end of another root to form an edge.}

\begin{defi}
\label{14.11}
Let $\Gamma_1$ and $\Gamma_2$ be two admissible graphs with identical numbers of
roots and $k_1$ and $k_2$ legs respectively.
We let $r$ be the number of roots of $\Gamma_1$,
let $k=k_1+k_2$ and let $I\sub[k]$ be a subset of $k_1$ elements.
We say that the triple
$(\Gamma_1,\Gamma_2,I)$ is an admissible triple if the following two holds:
\newline
1. The weight assignments $\mu_{1,i}$ and $\mu_{2,i}$ of the roots of $\Gamma_1$
and $\Gamma_2$ coincides in the sense that $\mu_{1,i}=\mu_{2,i}$ for all $i=1,\cdots,r$.
\newline
2. After connecting the $i$-th root in $\Gamma_1$ with the $i$-th root
in $\Gamma_2$ for all $i$, we obtain a new graph with $k$-edges (but without roots).
We demand that this graph is connected.
\end{defi}

Let $\eta=(\Gamma_1,\Gamma_2,I)$ be an admissible triple. Using the subset $I\sub[k]$
we obtain a unique bijection $[k_1]\cup [k_2]\to[k]$ so that it preserves the orders of
$[k_1]$ and $[k_2]$ and the image of $[k_1]$ is $I\sub [k]$. Using this
bijection, we get a unique ordering of the legs of this new graph. The genus of this graph is
$$g(\eta)=r+1-|V(\Gamma)|+\sum_{v\in\ver(\Gamma_1)\cup\ver(\Gamma_2)}
g(v)
$$
and the degree of this graph is
$$b(\eta)=\sum_{i=1}^{k_1} b_{\Gamma_1}(i)\cdot c_1(H_{Y_1})+\sum_{i=1}^{k_2} b_{\Gamma_2}(i)
\cdot c_1(H_{Y_2}).
$$
For simplicity, we denote by $|\eta|$ the triple
$$|\eta|=(g(\eta),b(\eta),k).
$$

Now let $f_1\mh\cX_1\to\cY_1$ be a relative stable morphism in $\Mydo(S)$,
where $(\cY_1,\cD_1)\in\YYo(S)$, and let $q_{1,i}\mh S\to\cX_1$ be the $i$-th
distinguished section.
Then each $f_1\circ q_{1,i}\mh S\to \cY_1$
factor through $\cD_1\sub\cY_1$, which composed with the first projection $\pr_D$ of
$\cD\cong D\times S$ defines a morphism $S\to D$. This way we obtain an evaluation morphism
\begin{equation}
\label{42.10}
\bq(f_1)=\bigl(\pr_D\circ f_1\circ q_{1,1},\cdots,\pr_D\circ f_1\circ q_{1,r}):S\to D^r.
\end{equation}
Similarly, let $f_2\mh\cX_2\to\cY_2$ be in $\Mydt(S)$ with distinguished sections
$q_{2,i}$. Then we have a similarly defined evaluation morphism
$\bq(f_2)\mh S\to D^r$.
Now we construct a new family of morphisms over $S$.
We glue the $i$-th sections $q_{1,i}(S)\sub\cX_1$ with the $i$-th section $q_{2,i}(S)\sub\cX_2$
for all $i=1,\cdots,r$ to obtain an $S$-family of nodal curves $\cX_1\adj\cX_2$.
Since $\eta=(\Gamma_1,\Gamma_2,I)$ is
admissible, the new family $\cX_1\adj\cX_2/S$ has connected fibers. Also the data
$I\sub [k]$ defines a natural ordering of the $k$ ordinary marked sections of $\cX_1\adj\cX_2$.
Now let $\cY_1^o\adj\cY_2$ be the $S$-scheme constructed in Corollary \ref{7.1}, then the pair
of morphisms $(f_1,f_2)$ defines a morphism
\begin{equation}
\label{14.14}
f_1\adj f_2:\cX_1\adj \cX_2\lra \cY_1^o\adj \cY_2
\end{equation}
if and only if
$$\bq(f_1)\equiv\bq(f_2): S\to D^r.
$$
We state it as a Proposition.

\begin{prop}
\label{14.12} Let $\eta=(\Gamma_1,\Gamma_2,I)$ be an admissible triple and let $f_i$
be relative stable morphisms in $\Mydi(S)$, $i=1,2$. Suppose $\bq(f_1)\equiv\bq(f_2)$,
then the morphism $f_1\adj f_2$ is an $S$-family of stable morphisms in $\Mwg(S)$ with
$\Gamma=|\eta|$.
\end{prop}

\begin{proof} The proof is straightforward and will be omitted.
\end{proof}

Let $f\in\mwg(S)$ be an $S$-family of stable morphisms to $\cW$. In case there is an admissible
triple $\eta=(\Gamma_1,\Gamma_2,I)$ and families $f_i\in\Mydi(S)$ so that $f=f_1\adj f_2$,
then we say $f$ is $\eta$-decomposable or $f$ is $\eta$-decomposed into families $f_1$ and $f_2$.

We let $\Mydi\to D^r$ be the morphisms induced by $\bq$ in (\ref{42.10}). Then this Proposition implies
that we have a natural morphism of stacks
\begin{equation}
\label{42.11}
\Mydo\times_{D^r}\Mydt\lra\Mwg.
\end{equation}
Note that if $r=0$, then $\eta$ is admissible implies that one of $\Gamma_i$ is empty.
Then the fiber product is understood to be either $\Mydt$ in case $\Gamma_1=\emptyset$ or
the otherwise.

Clearly, the above morphism factor through the substack $\Mwg\times_C 0$, and is finite and representable.
We let $\myy$ be the image stack in $\mwg$. We denote the induced morphism by
\begin{equation}
\label{42.74}
\Phi_\eta:\Mydo\times_{D^r}\Mydt\lra\myy.
\end{equation}

We now investigate the degree of $\Phi_\eta$. Let $f\mh X\to\wn_0$ be a stable morphism that lies in
the image of $\Phi_\eta$. By definition, there is an integer $l\in[n+1]$
so that we can $\eta$-decompose $f$ along the divisor $\bD_l\cap \wn_0$. Namely, we have a
pair of relative stable morphisms $f_i\mh X_i\to Y_i[n_i]$, where $n_1=l-1$ and $n_2=n+1-l$,
an ordered finite set $\Sigma$ and embeddings $\Sigma\sub X_i$ so that $X$ is the result of gluing
$X_1$ and $X_2$ along $\Sigma\sub X_1$ and $\Sigma\sub X_2$, and the restriction of $f$
to $X_i$ is $f_i$ using the gluing $Y_1[n_1]^o\adj Y_2[n_2]=\wn_0$.
Further, if we view $\Sigma\sub X_i$ as its ordered
distinguished marked points, then the topological type of $f_i$ with the induced ordering
on its ordinary marked points is exactly the $\Gamma_i$ in $\eta=(\Gamma_1,\Gamma_2,I)$.
Here $I\sub [k]$ is the subset of those marked points that are in $X_1$.
We fix such a ordering on $\Sigma$.

Now let $\sigma\sub S_r$ be any permutation. Then $\sigma$ induces a new ordering
on $\Sigma$, which we denote by $\Sigma^\sigma$. Accordingly, for the $\Gamma_i$ in $\eta=(\Gamma_1,
\Gamma_2,I)$, we define $\Gamma_i^\sigma$ to be the graph $\Gamma_i$ with its roots reordered
according to $\sigma$. Namely, the $j$-th
root of $\Gamma_i^\sigma$ is the $\sigma\upmo(j)$-root of $\Gamma_i$.
We define $\eta^\sigma=(\Gamma_1^\sigma,\Gamma_2^\sigma,I)$.
For any two admissible triples $\eta_1$ and $\eta_2$ we say $\eta_1\sim\eta_2$ (called equivalent)
if there is
a $\sigma\in S_r$ so that $\eta_1$ is isomorphic to $\eta_2^\sigma$.
We define $\Eq(\eta)$ to be the subgroup of $\sigma\in S_r$ so that $\eta\sim\eta^\sigma$.
Note that $\eta\sim\eta^\sigma$
implies that whenever $\sigma(j_1)=j_2$ then the $j_1$-th and the $j_2$-th roots of $\Gamma_i$ have
identical weights $\mu_{j_1}=\mu_{j_2}$ and are attached to the same vertex of $\Gamma_i$.

\begin{prop}
The morphism $\Phi_\eta$ is finite and \'etale. It has pure degree\footnote{Here we say $f\mh X\to Y$
has pure degree $d$ if for any integral $A\sub Y$ the restriction $X\times_Y A\to A$ has degree $d$.}
$|\Eq(\eta)|$.
\end{prop}

\begin{proof}
The proof that the morphism $\Phi_\eta$ is finite and \'etale is straightforward, and will be left
to readers. We now check that the degree of $\Phi_\eta$ is as stated.
Let $\xi\in\image (\Phi_\eta)$ be any element.
It follows from the discussion preceding to the statement of this Proposition that there is
an $l\in[n+1]$ so that the decomposition of $f$ along the divisor $\bD_l\cap \wn_0$
is the $\eta$-decomposition of $f$. Because $f$ is stable, for all other $l\pri\in[n+1]-l$
the decomposition types of $f$ along $\bD_{l\pri}\cap \wn_0$ will be different from $\eta$
(i.e. $\not\sim$).
We now let $\xi_i$ be represented by $f_i\mh X_i\to Y_i[n_i]$,
after fixing an ordering on $\Sigma=f\upmo(\bD_l)$ so that the topological type of
$\xi$ is $\Gamma_i$. Thus $f_i\in\Mydi$ and
$f=f_1\adj f_2$. Now let $\sigma\in S_r$ be any permutation. We let $f_i^\sigma$ be the same
morphism $f_i\mh X_i\to Y_i[n_i]$ except that the distinguished marked points of $X_i$
are reordered according to $\sigma$. This way
$$(f_i^\sigma,f_2^\sigma)\in \MM(\YY_1\urel,\Gamma_1^\sigma)\times_{D^r}\MM(\YY_1\urel,\Gamma_1^\sigma).
$$
Clearly, it is in the domain of $\Phi_\eta$ if and only if $\eta\sim\eta^\sigma$.

To derive the degree formula, we need to investigate the automorphism group $\Aut(\xi)$
of $\xi\in \mwg$ and the automorphism group $\Aut(\xi_i)$ of $\Mydi$.
First, the automorphism group $\Aut((\xi_1,\xi_2))$ is naturally isomorphic to
$\Aut(\xi_1)\times\Aut(\xi_2)$. Because $f$ is derived from gluing $f_1$ and $f_2$
along $\Sigma\sub X_1$ and $\Sigma\sub X_2$,
elements in $\Eq(\eta)$ induce permutations of $\Sigma$. Hence we obtain a natural homomorphism of groups
$h\mh \Aut(\xi)\to S_r$ whose kernel is isomorphic to $\Aut(\xi_1)\times\Aut(\xi_2)$.
Further the image of $h$ lies in the subgroup $\Eq(\eta)\sub S_r$
and the coset $\Eq(\eta)/\image(\eta)$ is exactly the set $\Phi_\eta\upmo(\xi)$. This shows that
$\Phi_\eta\upmo(\xi)$ consists of exactly $|\Eq(\eta)|/|\image(h)|$ distinct elements.
Hence the degree of $\Phi_\eta$ is
$$|\Eq(\eta)|/|\image(\eta)|\cdot |\Aut(\xi)|\cdot 1/|\Aut(\xi_1)\times\Aut(\xi_2)|=|\Eq(\eta)|.
$$
This proves the Proposition.
\end{proof}

\section{Some Comments}
It is interesting to work out the analogous construction for the moduli spaces of stable sheaves
over schemes. Before the discovery of SW-invariants, there were a lot of research work on
using degeneration method to study the Donaldson invariants. The differential geometric approach
to this is well understood. But the algebro-geometric approach to this is still missing. It is clear that
one should consider the moduli of stable sheaves over the stack $\WWc$ and the stack $\Zzd$.
It is interesting to work this out in detail. One technical issue is that after using the expanded degenerations
one should be able to avoid non-locally freeness (of sheaves) along the nodal divisor.
This is the case for curves, following the pioneer
work of \cite{GM}. It is true for surfaces, which is essentially proved by Gieseker and the author in \cite{GL}.

A more challenging degeneration problem is to workout a degeneration of GW-invariants for the families $W\to C$
whose singular fibers are only assumed to have normal crossing singularities. Recall that the family $W/C$ studied in
this paper only deal with the case where $W_0$ has two irreducible components that intersect along a connected
smooth divisor. It is easy to see that this construction can be generalized to the case where $W_0$ is only assumed
to have normal crossing singularity and its singular locus is smooth. When the singular locus of $W_0$ is not smooth,
it is not clear what constitute the class of expanded degenerations of $W/C$.
The progress along this line is important
to the study of degenerations of moduli spaces.

\end{document}